\documentclass[aop,MSNbibl,seceqn,dvips]{arximspdf}
\usepackage{mathbh}
\usepackage{graphicx}

\doi{10.1214/13-AOP855} %kopijuoti is PTS
\volume{42}
\issue{5}
\pubyear{2014}
\firstpage{1980}
\lastpage{2031}

\makeatletter
\newcommand{\widebar}{\overline}

\newcommand{\rrVert}{\Vert}
\newcommand{\rrvert}{\vert}
\newcommand{\llVert}{\Vert}
\newcommand{\llvert}{\vert}
\newcommand{\ii}{\mathrm{i}}
\newcommand{\IE}{\mathbb{IE}}
\newtheorem{theorem}{Theorem}[section]
\newtheorem{lemma}[theorem]{Lemma}
\newtheorem{proposition}[theorem]{Proposition}
\newproclaim{definition}[theorem]{Definition}
\newproclaim{remark}[theorem]{Remark}
\makeatother

\begin{document}
\begin{frontmatter}

\title{The outliers of a deformed Wigner matrix}
\runtitle{The outliers of a deformed Wigner matrix}

\begin{aug}
\author[A]{\fnms{Antti} \snm{Knowles}\corref{}\ead[label=e1]{knowles@math.ethz.ch}\thanksref{t1}}
\and
\author[B]{\fnms{Jun} \snm{Yin}\ead[label=e2]{jyin@math.uwisc.edu}\thanksref{t2}}
\runauthor{A. Knowles and J. Yin}
\affiliation{New York University and University of Wisconsin}
\address[A]{Department of Mathematics\\
ETH Z\"urich\\
R\"{a}mistrasse 101\\
8092 Z\"urich\\
Switzerland\\
\printead{e1}} %adresu isvedimo komanda gale!
\address[B]{Department of Mathematics\\
University of Wisconsin\\
Madison, Wisconsin 53706\\
USA\\
\printead{e2}}
\end{aug}
\thankstext{t1}{Supported in part by NSF Grant DMS-07-57425 Swiss National Science Foundation Grant 144662.}
\thankstext{t2}{Supported in part by NSF Grant DMS-12-07961.}

% HISTORY:
\received{\smonth{7} \syear{2012}}
\revised{\smonth{3} \syear{2013}}

% ABSTRACT
%
\begin{abstract}
We derive the joint asymptotic distribution of the outlier eigenvalues
of an additively deformed Wigner matrix $H$. Our only assumptions on
the deformation are that its rank be fixed and its norm bounded. Our
results extend those of [The isotropic semicircle law and deformation of Wigner
  matrices. Preprint] by admitting overlapping outliers
and by computing the joint distribution of all outliers. In particular,
we give a complete description of the failure of universality first
observed in [\textit{Ann. Probab.} \textbf{37} (2009) 1--47;
\textit{Ann. Inst. Henri Poincar\'e Probab. Stat.} \textbf{48} (1013) 107--133;
Free convolution with a semi-circular
distribution and eigenvalues of spiked deformations of Wigner matrices.
Preprint]. We also show that, under suitable
conditions, outliers may be strongly correlated even if they are far
from each other. Our proof relies on the isotropic local semicircle law
established in [The isotropic semicircle law and deformation of Wigner
  matrices. Preprint]. The main technical achievement of the
current paper is the joint asymptotics of an arbitrary finite family of
random variables of the form $\langle\mathbf{v}, (H - z)^{-1}
\mathbf{w}\rangle$.
\end{abstract}

% KEYWORDS
\begin{keyword}[class=AMS]
\kwd{15B52}
\kwd{60B20}
\kwd{82B44}
\end{keyword}
\begin{keyword}
\kwd{Random matrix}
\kwd{universality}
\kwd{deformation}
\kwd{outliers}
\end{keyword}

\end{frontmatter}

\setcounter{footnote}{2}
%s1 #&#
\section{Introduction}

In this paper, we study a Wigner matrix $H$---a random $N \times N$
matrix whose entries are independent up to symmetry constraints---that
has been deformed by the addition of a finite-rank matrix $A$ belonging
to the same symmetry class as $H$. By Weyl's eigenvalue interlacing
inequalities, such a deformation does not influence the global
statistics of the eigenvalues as $N \to\infty$. Thus, the empirical
eigenvalue densities of the deformed matrix $H + A$ and the undeformed
matrix $H$ have the same large-scale asymptotics, and are governed by
Wigner's famous semicircle law. However, the behavior of individual
eigenvalues may change dramatically under such a deformation. In
particular, deformed Wigner matrices may exhibit
\textit{outliers}---eigenvalues detached from the bulk spectrum. They
were first investigated in \cite{FKoml} for a particular rank-one
deformation. Subsequently, much progress
\cite{SoshPert,SoshPert2,FP,CDMF1,CDMF2,CDMF3,BGN,BGGM1,BGGM2,KY2} has
been made in the understanding of the outliers of deformed Wigner
matrices. We refer to \cite{SoshPert,SoshPert2,KY2} for a more detailed
review of recent developments.

We normalize $H$ so that its spectrum is asymptotically given by the
interval $[-2,2]$. The creation of an outlier is associated with a
sharp transition, where the magnitude of an eigenvalue $d_i$ of $A$
exceeds the threshold $1$. As $d_i$ (resp.,~$-d_i$) becomes larger than
$1$, the largest (resp., smallest) nonoutlier eigenvalue of $H + A$
detaches itself from the bulk spectrum and becomes an outlier. This
transition is conjectured to take place on the scale $\llvert  d_i
\rrvert  - 1 \sim N^{-1/3}$. In fact, this scale was established in
\cite{BBP,Pec,BV1,BV2} for the special cases where $H$ is
Gaussian---the Gaussian Orthogonal Ensemble (GOE) and the Gaussian
Unitary Ensemble (GUE). We sketch the results of \cite{BBP,Pec,BV1,BV2}
in the case of additive deformations of GOE/GUE. For simplicity, we
consider rank-one deformations, although the results of
\mbox{\cite{BBP,Pec,BV1,BV2}} cover arbitrary finite-rank deformations. Let
the eigenvalue $d$ of $A$ be of the form $d = 1 + w N^{-1/3}$ for some
fixed $w \in\mathbb{R}$. In \cite{BBP,Pec,BV1,BV2}, the authors proved
for any fixed $w$ the weak convergence
\[
N^{2/3} \bigl(\lambda_N (H + A ) - 2 \bigr)\quad
\Longrightarrow\quad\Lambda_w,
\]
where $\lambda_N(H+A)$ denotes the largest eigenvalue of $H+A$. In
particular, the largest eigenvalue of $H+A$ fluctuates on the scale
$N^{-2/3}$. Moreover, the asymptotics in $w$ of the law $\Lambda_w$ was
analysed in \cite{BBP,Pec,BV1,BV2,Bthesis}: as $w \to+ \infty$ (and
after an appropriate affine scaling), the law $\Lambda_w$ converges to
a Gaussian; as $w \to- \infty$, the law $\Lambda_w$ converges to the
Tracy--Widom-$\beta$ distribution (where $\beta= 1$ for GOE and $\beta=
2$ for GUE), which famously governs the distribution of the largest
eigenvalue of the underformed matrix $H$ \cite{TW1,TW2}.

The proofs of \cite{BBP,Pec} use an asymptotic analysis of Fredholm
determinants, while those of \cite{BV1,BV2,Bthesis} use an explicit
tridiagonal representation of $H$; both of these approaches rely
heavily on the Gaussian nature of $H$. In order to study the phase
transition for non-Gaussian matrix ensembles, and in particular address
the question of spectral universality, a different approach is needed.
Interestingly, it was observed in \cite{CDMF1,CDMF2,CDMF3} that the
distribution of the outliers is not universal, and may depend on the
law of $H$ as well as the geometry of the eigenvectors of $A$. The
nonuniversality of the outliers was further investigated in
\cite{SoshPert,SoshPert2,KY2}.

In a recent paper \cite{KY2}, we considered finite-rank deformations of
a Wigner matrix whose entries have subexponential decay. The two main
results of \cite{KY2} may be informally summarized as follows.
\begin{longlist}[(a)]
\item[(a)] We proved that the nonoutliers of $H + A$ \textit{stick}
    to the extremal eigenvalues of the original Wigner matrix $H$
    with high precision, provided that each eigenvalue $d_i$ of $A$
    satisfies $\llvert  \llvert  d_i \rrvert  - 1 \rrvert
    \geq(\log N)^{C \log\log N} N^{-1/3}$.

\item[(b)] We identified the asymptotic distribution of a single
    outlier, provided that (i) it is separated from the asymptotic
    bulk spectrum $[-2,2]$ by at least $(\log N)^{C \log\log N}
    N^{-2/3}$ and (ii) it does not overlap with any other outlier
    of $H + A$. Here, two outliers are said to \textit{overlap} if
    their separation is comparable to the scale on which they
    fluctuate; see Section~\ref{secheuristics} below for a precise
    definition.
\end{longlist}

Note that the assumption (i) of (b) is optimal, up to the logarithmic
factor $(\log N)^{C \log\log N}$. Indeed, the extremal bulk eigenvalues
of $H + A$ are known \cite{KY2}, Theorem~2.7,\vadjust{\goodbreak} to fluctuate on the scale
$N^{-2/3}$; for an eigenvalue of $H + A$ to be an outlier, therefore,
we require that its distance from the asymptotic bulk spectrum $[-2,2]$
be much greater than $N^{-2/3}$. See Section~\ref{secheuristics} below
for more details.

The goal of this paper is to extend the result (b) by obtaining a
complete description of the asymptotic distribution of the outliers.
Our only assumptions on the deformation $A \equiv A_N$ are that its
rank be fixed and its norm bounded. (In particular, the eigenvalues of
$A$ may depend on $N$ in an arbitrary fashion, provided they remain
bounded, and its eigenvectors may be an arbitrary orthonormal family.)
Our main result gives the asymptotic joint distribution of all
outliers. Here, an outlier is by definition an eigenvalue of $H+A$
whose classical location [see~(\ref{deftheta}) below] is separated from
the asymptotic bulk spectrum $[-2,2]$ by at least $(\log N)^{C \log\log
N} N^{-2/3}$ for some (large) constant $C$. Our main result is given in
Theorem~\ref{thmmainresultJ} below.

Thus, in this paper we extend the result (b) in two directions: we
allow overlapping outliers, and we derive the joint asymptotic
distribution of all outliers. The distribution of overlapping outliers
is more complicated than that of nonoverlapping outliers, as
overlapping outliers exhibit a level repulsion similar to that among
the bulk eigenvalues of Wigner matrices. This repulsion manifests
itself by the joint distribution of a group of overlapping outliers
being given by the distribution of eigenvalues of a small (explicit)
random matrix [see (\ref{referencematrix}) below]. The mechanism
underlying the repulsion among outliers is therefore the same as that
for the eigenvalues of GUE: the Jacobian relating the
eigenvalue--eigenvector entries to the matrix entries has a Vandermonde
determinant structure, and vanishes if two eigenvalues coincide.
Several special cases of overlapping outliers have already been studied
in the works \cite{SoshPert,SoshPert2,CDMF1,CDMF2,CDMF3}, which in
particular exhibited the level repulsion mechanism described above.

Due to this level repulsion, overlapping outliers are obviously not
asymptotically independent. A novel observation, which follows from our
main result, is that in general nonoverlapping outliers are not
asymptotically independent either; in this case the lack of
independence does not arise from level repulsion, but from a~more
subtle interplay between the distribution of $H$ and the geometry of
the eigenvectors of $A$. In some special cases, such as
GOE/GUE,\vadjust{\goodbreak}
nonoverlapping outliers are, however, asymptotically independent. More
precisely, our main result (Theorem~\ref{thmmainresultJ} below) shows
that two outliers may, under suitable conditions on $H$ and $A$, be
strongly correlated in the limit $N \to\infty$, even if they are far
from each other (e.g., on opposite sides of the bulk spectrum).

Finally, we note that throughout this paper we assume that the entries
of $H$ have subexponential decay. We need this assumption because our
proof relies heavily on the local semicircle law and eigenvalue
rigidity estimates for $H$, proved in \cite{EYY3} under the assumption
of subexponential decay. However, this assumption is not fundamental to
our approach, which may be combined with the recent methods for dealing
with heavy-tailed Wigner matrices developed in \cite{EKYY1,EKYY2,LY}.
Moreover, the assumption that the norm of $A$ be bounded may be easily
removed; in fact, large eigenvalues of $A$ are easier to treat than
small ones.

We remark that recently Pizzo, Renfrew and Soshnikov
\cite{SoshPert,SoshPert2} took a different approach, and derived the
asymptotic distribution of a single group of overlapping outliers under
optimal tail assumptions on $H$. On the other hand, in
\cite{SoshPert,SoshPert2} it is assumed that the eigenvalues of $A$ are
independent of $N$ and that its eigenvectors satisfy a condition which
roughly constrains them to be either strongly localized or delocalized.

%s1.1 #&#
\subsection{Outline of the proof}
As in \cite{KY2}, our proof relies on the \textit{isotropic local
semicircle law}, proved in \cite{KY2}, Theorems 2.2~and~2.3. The
isotropic local semicircle law is an extension of the \textit{local
semicircle law}, whose study was initiated in \mbox{\cite{ESY1,ESY3}}. The
local semicircle law has since become a cornerstone of random matrix
theory, in particular in establishing the universality of Wigner
matrices \cite{ESY4,ESY6,EYY2,EYY3,TV1,TV2}. The strongest versions of
the local semicircle law, proved in \mbox{\cite{EYY3,EKYY1}}, give precise
estimates on the local eigenvalue density, down to scales containing
$N^\varepsilon$ eigenvalues. In fact, as formulated in \cite{EYY3}, the
local semicircle law gives optimal high-probability estimates on the
quantity
%
%e1.1 #&#
\begin{equation}
\label{oldlsc} G_{ij}(z) - \delta_{ij} m(z),
\end{equation}
where $m(z)$ denotes the Stieltjes transform of Wigner's semicircle law
and $G(z):=(H - z)^{-1}$ is the resolvent of $H$.

The isotropic local semicircle law is a generalization of the local
semicircle law, in that it gives optimal high-probability estimates on
the quantity
%
%e1.2 #&#
\begin{equation}
\label{quantityinislc} \bigl\langle\mathbf{v},  \bigl(G(z) - m(z)
\mathbh{1} \bigr) \mathbf{w} \bigr\rangle,
\end{equation}
where $\mathbf{v}$ and $\mathbf{w}$ are arbitrary deterministic
vectors. Clearly, (\ref{oldlsc}) is a special case obtained from
(\ref{quantityinislc}) by setting $\mathbf{v} = \mathbf{e}_i$ and
$\mathbf{w} = \mathbf{e}_j$, where $\mathbf{e}_i$ denotes $i$th
standard basis vector of $\mathbb{C}^N$.

As in the works \cite{SoshPert,SoshPert2,KY2}, a major part of our
proof consists in deriving the asymptotic distribution of the entries
of $G(z)$. The main technical achievement of this paper is to obtain
the \textit{joint} asymptotics of an arbitrary finite family of
variables of the form $\langle\mathbf{v}, {G(z) \mathbf{w}}\rangle$,
whereby the spectral parameters $z$ of different entries may differ,
and are assumed to satisfy $2 + (\log N)^{C \log\log N} N^{-2/3}
\leq\llvert \operatorname{Re}z \rrvert \leq C$ for some positive constant $C$. The
question of the joint asymptotics of the resolvent entries occurs more
generally in several problems on deformed random matrix models, and we
therefore believe that the techniques of this paper are also of
interest for other problems on deformed matrix ensembles.

An important ingredient in our proof is the four-step strategy
introduced in \cite{KY2}. It may be summarized as follows: (i)
reduction to the distribution of the resolvent~$G$, (ii) the case of
Gaussian $H$, (iii) the case of almost Gaussian $H$, (iv) the case of
general $H$. Steps (i)--(iii) in the current paper are substantially
different from their counterparts in \cite{KY2}; this results from
treating an entire overlapping group of outliers simultaneously, as
well as from the need to develop an argument that admits an analysis of
the joint law of different groups. In fact, for pedagogical reasons,
first---in Sections~\mbox{\ref{secreduction2}--\ref{secgencase}}---we give the
proof for the case of a single group of overlapping
outliers,\footnote{In the resolvent language, this means that the
spectral parameters $z$ of all the resolvent entries coincide.} and
then---in Section~\ref{secgendistrJ}---extend it to yield the full
joint distribution. In contrast to the steps (i)--(iii), step (iv)
survives almost unchanged from \cite{KY2}, and in
Section~\ref{secgencase} we give an explanation of the required
modifications.

Another ingredient of our proof is a two-level partitioning of the
outliers combined with near-degenerate perturbation theory for
eigenvalues. Roughly, outliers are partitioned into blocks depending on
whether they overlap. In the finer partition, denoted by $\Pi$ below
(see Definition~\ref{defPi}), we regroup two outliers into the same
block if their mean separation is bounded by some large constant
(denoted by $s$ below) times the magnitude of their fluctuations. Due
to logarithmic error factors of the form $(\log N)^{C \log\log N}$ that
appear naturally in high-probability estimates pervading our proof, we
shall require a second, coarser, partition, denoted by $\Gamma$ below
(see Definition~\ref{defjointcoarse}). In $\Gamma$, we regroup two
outliers into the same block if their mean separation is bounded by
$(\log N)^{C \log\log N}$ times the magnitude of their fluctuations.
The link between $\Gamma$ and $\Pi$ is provided by perturbation theory,
and is performed in Sections~\ref{secconcofproof} (for a single group)
and \ref{secjointdist} (for the full joint distribution).

%s2 #&#
\section{Formulation of results} \label{secresults}

%s2.1 #&#
\subsection{The setup}
Let $H = (h_{ij})_{i,j = 1}^N$ be an $N \times N$ random matrix. We
assume that the upper-triangular entries $(h_{ij}\dvtx i \leq j)$ are
independent complex-valued random variables. The remaining entries of
$H$ are given by imposing $H = H^*$. Here $H^*$ denotes the Hermitian
conjugate of $H$. We assume that all entries are centred, $\mathbb
{E}h_{ij} = 0$. In addition, we assume that one of the two following
conditions holds.
\begin{longlist}[(ii)]
\item[(i)] \textit{Real symmetric Wigner matrix}: $h_{ij} \in
\mathbb{R}$ for all $i,j$ and
\[
\mathbb{E}h_{ii}^2 = \frac{2}{N},\qquad
\mathbb{E}h_{ij}^2 = \frac{1}{N} \qquad(i \neq j).
\]

\item[(ii)] \textit{Complex Hermitian Wigner matrix}:
\[
\mathbb{E}h_{ii}^2 = \frac{1}{N}, \qquad\mathbb{E}
\llvert h_{ij} \rrvert ^2 = \frac{1}{N}, \qquad
\mathbb{E}h_{ij}^2 = 0 \qquad(i \neq j).
\]
\end{longlist}
We introduce the usual index $\beta$ of random matrix theory, defined
to be $1$ in the real symmetric case and $2$ in the complex Hermitian
case. We use the abbreviation GOE/GUE to mean GOE if $H$ is a real
symmetric Wigner matrix with Gaussian entries and GUE if $H$ is a
complex Hermitian Wigner matrix with Gaussian entries. We assume that
the entries of $H$ have uniformly subexponential decay, that is, that
there exists a constant $\vartheta> 0$ such that
%
%e2.1 #&#
\begin{equation}
\label{subexpforh} \mathbb{P} \bigl(\sqrt{N} \llvert h_{ij}
\rrvert \geq x \bigr) \leq\vartheta ^{-1} \exp\bigl(- x^\vartheta
\bigr)
\end{equation}
for all $i$, $j$ and $N$. Note that we do not assume the entries of $H$
to be identically distributed, and we do not require any smoothness in
the distribution of the entries of $H$.

We consider a deformation of fixed, finite rank $r \in\mathbb{N}$. Let
$V \equiv V_N$ be a deterministic $N \times r$ matrix satisfying $V^* V
= \mathbh{1}_r$, and $D \equiv D_N$ be a deterministic $r \times r$
diagonal matrix whose eigenvalues are nonzero. Both $V$ and $D$ depend
on $N$. We sometimes also use the notation $V =
[\mathbf{v}^{(1)},\ldots, \mathbf{v}^{(r)}]$, where
$\mathbf{v}^{(1)},\ldots, \mathbf{v}^{(r)} \in\mathbb{C}^N$ are
orthonormal, as well as $D = \operatorname{diag}(d_1, \dots, d_r)$. We always assume
that the eigenvalues of $D$ satisfy
%
%e2.2 #&#
\begin{equation}
\label{basicconditionsond} -\Sigma+ 1 \leq d_1 \leq
d_2 \leq\cdots\leq d_r \leq\Sigma- 1,
\end{equation}
where $\Sigma$ is some fixed positive constant. We are interested in
the spectrum of the deformed matrix
\[
\widetilde H:= H + V D V^* = H + \sum_{i = 1}^r
d_i \mathbf{v}^{(i)} \bigl(\mathbf{v}^{(i)}
\bigr)^*.
\]

The following definition summarizes our conventions for the spectrum of
a matrix. For our purposes, it is important to allow the matrix entries
and its eigenvalues to be indexed by an arbitrary subset of positive
integers.

%
%de2.1 #&#
\begin{definition} \label{defeigenvalues}
Let $\pi$ be a finite set of positive integers, and let $A =
(A_{ij})_{i,j \in\pi}$ be a $\llvert  \pi\rrvert  \times\llvert
\pi\rrvert $ Hermitian
matrix whose entries are indexed by elements of $\pi$. We denote by
\[
\sigma(A):= \bigl(\lambda_i(A)\bigr)_{i \in\pi} \in\mathbb
{R}^{\pi}
\]
the family of eigenvalues of $A$. We always order the eigenvalues so
that $\lambda_i(A) \leq\lambda_j(A)$ if $i \leq j$.

By a slight abuse of notation, we sometimes identify $\sigma(A)$ with
the set $\{\lambda_i(A)\}_{i \in\pi} \subset\mathbb{R}$. Thus, for
instance, $\operatorname{dist}(\sigma(A), \sigma(B) ):=\min_{i,j}
\llvert \lambda_i(A) - \lambda_j(B) \rrvert $ denotes the distance
between $\sigma(A)$ and $\sigma(B)$ viewed as subsets of $\mathbb{R}$.
\end{definition}

We abbreviate the (random) eigenvalues of $H$ and $\widetilde H$ by
\[
\lambda_\alpha:= \lambda_\alpha(H), \qquad\mu_\alpha:= \lambda_\alpha(\widetilde H).
\]
The following definition introduces a convenient notation for minors of
matrices.

%
%de2.2 #&#
\begin{definition}[(Minors)] \label{defsmallminors}
For an $r \times r$ matrix $A = (A_{ij})_{i,j = 1}^r$ and a subset $\pi
\subset\{1,\ldots,r\}$ of integers, we define the $\llvert  \pi\rrvert
\times \llvert  \pi\rrvert $ matrix
\[
A_{[\pi]} = (A_{ij})_{i,j \in\pi}.
\]
\end{definition}

We shall frequently make use of the logarithmic control parameter
%
%e2.3 #&#
\begin{equation}
\varphi\equiv\varphi_N:= (\log N)^{\log\log N}.
\end{equation}
The interpretation of $\varphi$ is that of a slowly growing parameter
[note that $\varphi\leq N^\varepsilon$ for any $\varepsilon> 0$ and
large enough $N \geq N_0(\varepsilon)$]. Throughout this paper, every
quantity that is not explicitly a constant may depend on $N$, with
the sole exception of the rank $r$ of the deformation, which is
required to be fixed. Unless needed, we consistently drop the argument
$N$ from such quantities.

We denote by $C$ a generic positive large constant, whose value may
change from one expression to the next. For two positive quantities
$A_N$ and $B_N$, we use the notation $A_N \asymp B_N$ to mean $C^{-1}
A_N \leq B_N \leq C A_N$ for some positive constant~$C$. Moreover, we
write $A_N \ll B_N$ if $A_N/B_N \to0$ and $A_N \gg B_N$ if $B_N \ll
A_N$. Finally, for $a < b$ we set $[\![a,b ]\!]:= [a,b]
\cap\mathbb{Z}$.

%s2.2 #&#
\subsection{Heuristics of outliers} \label{secheuristics}
Before stating our results, we give a heuristic description of the
behavior of the outliers. An eigenvalue $d_i$ of $D$ satisfying
%
%e2.4 #&#
\begin{equation}
\label{condforoutlier} \llvert d_i \rrvert - 1 \gg
N^{-1/3}
\end{equation}
gives rise to an outlier $\mu_{\alpha(i)}$ located around its classical
location $\theta(d_i)$, where we defined, for $d \in\mathbb
{R}\setminus
(-1,1)$,
%
%e2.5 #&#
\begin{equation}
\label{deftheta} \theta(d):= d + \frac{1}{d}
\end{equation}
and
%
%e2.6 #&#
\begin{equation}
\label{defalphai} \alpha(i):= \cases{ i, &\quad if $d_i < 0$,
\vspace*{2pt}
\cr
N - r + i, &\quad if $d_i > 0$.}
\end{equation}

Condition (\ref{condforoutlier}) may be heuristically understood as
follows; for simplicity set $r = 1$ and $D = d > 1$. The extremal
eigenvalues of $\widetilde H$ that are not outliers fluctuate on the
scale $N^{-2/3}$ (see \cite{KY2}, Theorem~2.7), the same scale as the
extremal eigenvalues of the undeformed matrix $H$. For the largest
eigenvalue $\mu_N$ of $\widetilde H$ to be an outlier, we require that
its separation from the asymptotic bulk spectrum $[-2,2]$, which is of
the order $\theta(d) - 2$, be much greater than $N^{-2/3}$. This leads
to condition (\ref{condforoutlier}) by a simple expansion of
$\theta$ around 1.

\begin{figure}%[ht!]

\includegraphics{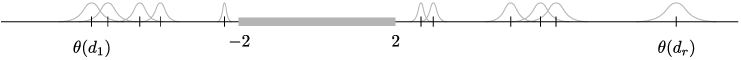}

\caption{A general outlier configuration. We draw the outlier $\mu
_{\alpha(i)}$ associated with $d_i$ using a black line marking its mean location
$\theta(d_i)$ and a grey curve indicating its probability density. The breadth of the
curve associated with $d_i$ is of the order $N^{-1/2} (\llvert  d_i \rrvert  - 1)^{1/2}$. Outliers
whose probability densities overlap satisfy (\protect\ref{overlappingcond}) [or, equivalently,
(\protect\ref{overlappingcond2})]. We do not draw the bulk eigenvalues,
which are contained in the grey bar.}\label{figureoutliers}
\end{figure}

The outlier $\mu_{\alpha(i)}$ associated with $d_i$ fluctuates on the
scale $N^{-1/2} (\llvert  d_i \rrvert  - 1)^{1/2}$. Thus,
$\mu_{\alpha(i)}$ fluctuates on the scale $N^{-1/2}$ if $d_i$ is
well-separated from the critical point $1$, and on the scale $N^{-2/3}$
if $d_i$ is critical, that is, $d_i = 1 + a N^{-1/3}$ for some fixed $a
> 0$. The outliers associated with $d_i$ and $d_j$ \textit{overlap} if
their separation is comparable to or less than the scale on which they
fluctuate. The overlapping condition thus reads
%
%e2.7 #&#
\begin{equation}
\label{overlappingcond} \bigl\llvert \theta(d_i) -
\theta(d_j) \bigr\rrvert \leq C N^{-1/2} \bigl(\llvert
d_i \rrvert - 1\bigr)^{1/2}
\end{equation}
for some (typically large) constant $C > 0$. Note that the factor
$\llvert  d_i \rrvert  - 1$ on the right-hand side could be replaced
with $\llvert  d_j \rrvert  - 1$. Indeed, recalling
(\ref{condforoutlier}), it is not hard to check that
(\ref{overlappingcond}) for some $C > 0$ is equivalent to
(\ref{overlappingcond}) with $d_i$ on the right-hand side replaced with
$d_j$ and the constant $C$ replaced with a constant $C' \asymp C$.
Using~(\ref{deftheta}) and recalling (\ref{condforoutlier}), we may
rewrite the overlapping condition (\ref{overlappingcond}) as
%
%e2.8 #&#
\begin{equation}
\label{overlappingcond2} N^{1/2} \bigl(\llvert d_i
\rrvert - 1\bigr)^{1/2} \llvert d_i - d_j
\rrvert \leq C
\end{equation}
for some $C > 0$. As in (\ref{overlappingcond}), $\llvert  d_i \rrvert
- 1$ may be replaced with $\llvert  d_j \rrvert  - 1$.
Figure~\ref{figureoutliers} summarizes the general picture of outliers.
%
%f1 #&#

%s2.3 #&#
\subsection{The distribution of a single group}
After these preparations, we state our results. We begin by defining a
reference matrix which will describe the distribution of a group of
overlapping outliers.\vspace*{-2pt} Define the moment matrices $\mu^{(3)} =
(\mu^{(3)}_{ij})$ and $\mu^{(4)} = (\mu^{(4)}_{ij})$ of $H$ through
\[
\mu^{(3)}_{ij}:= N^{3/2} \mathbb{E} \bigl(\llvert
h_{ij} \rrvert ^2 h_{ij} \bigr), \qquad
\mu^{(4)}_{ij}:= N^2 \mathbb{E}\llvert h_{ij} \rrvert ^4.
\]
Using the matrices $\mu^{(3)}$ and $\mu^{(4)}$, we define the
deterministic functions
\begin{eqnarray*}
\mathcal P_{ij,kl}(R) &:=& R_{il} R_{kj} +
\mathbf{1}(\beta= 1) R_{ik} R_{jl},
\\
\mathcal Q_{ij,kl}(V) &:=& \frac{1}{\sqrt{N}} \sum
_{a,b} \bigl(\widebar V_{ai} \widebar V_{ak}
V_{al} \mu^{(3)}_{ab} V_{bj} + \widebar V_{ia} \mu ^{(3)}_{ab} V_{bj} \widebar V_{bk} V_{bl}
\\
&&\hspace*{39pt}{} + \widebar V_{ak} \widebar V_{ai}
V_{aj} \mu ^{(3)}_{ab} V_{bl} + \widebar V_{ka} \mu^{(3)}_{ab} V_{bl} \widebar V_{bi} V_{bj} \bigr),
\\
\mathcal R_{ij,kl}(V) &:=& \frac{1}{N} \sum
_{a,b} \bigl(\mu ^{(4)}_{ab} - 4 + \beta
\bigr) \widebar V_{bi} V_{bj} \widebar V_{bk}
V_{bl},
\end{eqnarray*}
where $i,j,k,l \in[\![1,r ]\!]$, $R$ is an $r \times r$
matrix, and $V$ an $N \times r$ matrix. Moreover, we define the
deterministic $r \times r$ matrix
\[
\mathcal S(V):= \frac{1}{N} V^* \mu^{(3)} V.
\]

%
%re2.3 #&#
\begin{remark}\label{remSQRbounded}
Using Cauchy--Schwarz and assumption (\ref{subexpforh}), it is easy
to check that $\mathcal P(V^* V)$, $\mathcal Q(V)$, $\mathcal R(V)$
and $\mathcal S(V)$ are uniformly bounded for $V$ satisfying $0 \leq
V^* V \leq\mathbh{1}$ (in the sense of quadratic forms).
\end{remark}

Next, let $\delta\equiv\delta_N$ be a positive sequence satisfying
$\varphi^{-1} \leq\delta\ll1$. (Our result will be independent of
$\delta$ provided it satisfies this condition; see
Remark~\ref{remchangedelta} below.) The sequence $\delta$ will serve as
a cutoff in the size of the entries of $V$ when computing the law of
$V^* H V$: entries of $V$ smaller than $\delta$ give rise to an
asymptotically Gaussian random variable by the central limit theorem;
the remaining entries are treated separately, and the associated random
variable is in general not Gaussian. Thus, we define the matrix
$V_\delta= (V^{\delta}_{ij})$ through
\[
V^{\delta}_{ij}:= V_{ij} \mathbf{1}\bigl(\llvert
V_{ij} \rrvert > \delta\bigr).
\]

For $\ell\in[\![1,r ]\!]$ satisfying $\llvert  d_\ell\rrvert
> 1$ we define the $r \times r$ matrix
%
%e2.9 #&#
\begin{equation}
\label{defUpsilon} \Upsilon^\ell:= \bigl(\llvert d_\ell
\rrvert + 1\bigr) \bigl(\llvert d_\ell \rrvert - 1\bigr)^{1/2}
\biggl(\frac{N^{1/2} V_\delta^* H V_\delta}{d_\ell^2} + \frac{\mathcal S(V)}{d_\ell^4} \biggr).
\end{equation}
Abbreviate
%
%e2.10 #&#
\begin{equation}
\label{defDelta} \Delta_{ij,kl}:= \mathcal P_{ij,kl}(
\mathbh{1}) = \delta_{il} \delta_{kj} + \mathbf{1}(\beta= 1)
\delta_{ik} \delta_{jl}.
\end{equation}
Note that $\Delta$ is nothing but the covariance matrix of a GOE/GUE
matrix: if $r^{-1/2} \Phi$ is an $r \times r$ GOE/GUE matrix then
$\mathbb{E}
\Phi_{ij} \Phi_{kl} = \Delta_{ij,kl}$. We introduce an $r \times r$
Gaussian matrix $\Psi^\ell$, independent of $H$, which is complex
Hermitian for $\beta= 2$ and real symmetric for $\beta= 1$. The
entries of $\Psi^\ell$ are centred, and their law is determined by the
covariance
%
%e2.11 #&#
\begin{eqnarray}\label{defofcov}
\qquad \mathbb{E}\Psi^\ell_{ij}\Psi^\ell_{kl} &=& \frac{\llvert  d_\ell\rrvert  + 1}{d_\ell^2} \Delta_{ij,kl}
+\bigl(\llvert d_\ell\rrvert + 1\bigr)^2 \bigl(\llvert
d_\ell\rrvert - 1\bigr)
\nonumber\\[-8pt]\\[-8pt]
&&{}\times \biggl(- \frac
{\mathcal P_{ij,kl}(V_\delta^* V_\delta)}{d_\ell^4} +
\frac{\mathcal Q_{ij,kl}(V)}{d_\ell^5}+ \frac{\mathcal
R_{ij,kl}(V)}{d_\ell^6} \biggr)+ E_{ij,kl}.\nonumber
\end{eqnarray}
Here $E_{ij,kl}:=\varphi^{-1} \Delta_{ij,kl}$ is a term, that is,
needed to ensure that the right-hand side of (\ref{defofcov}) is a
nonnegative $r^2 \times r^2$ matrix. This nonnegativity follows as a
by-product of our proof, in which the right-hand side of
(\ref{defofcov}) is obtained from the covariance of an explicit random
matrix; see Proposition~\ref{propalmostgaussian} below for more
details. Note that the term $E_{ij,kl}$ does not influence the
asymptotic distribution of~$\Psi^\ell$.

%
%re2.4 #&#
\begin{remark}\label{remchangedelta}
A different choice of $\delta$, subject to $\varphi^{-1} \leq\delta
\ll1$, leads to the same asymptotic distribution for $\Upsilon^\ell+
\Psi^\ell$. This is an easy consequence of the central limit theorem
and the observation that the matrix entries
\[
\biggl(\bigl(\llvert d_\ell\rrvert + 1\bigr) \bigl(\llvert
d_\ell\rrvert - 1\bigr)^{1/2} \frac{N^{1/2}
V_\delta^* H V_\delta}{d_\ell^2}
\biggr)_{ij}
\]
have covariance matrix $(\llvert  d_\ell\rrvert  + 1)^2 (\llvert
d_\ell\rrvert  - 1) d_\ell^{-4} \mathcal P_{ij,kl}(V_\delta^*
V_\delta)$.
\end{remark}

Before stating our result in full generality, we give a special case
which captures its essence and whose statement is somewhat simpler.

%
%th2.5 #&#
\begin{theorem} \label{thmeasystatement}
For large enough $K$ the following holds. Let $\pi\subset[\![1,r ]\!]$
be a subset of consecutive integers, and fix $\ell\in\pi$. Suppose that
$\llvert  d_\ell\rrvert  \geq1 + \varphi^K N^{-1/3}$. Suppose moreover
that there is a constant $C$ such that
%
%e2.12 #&#
\begin{equation}
\label{samegroup} N^{1/2} \bigl(\llvert d_\ell\rrvert - 1
\bigr)^{1/2} \llvert d_i - d_\ell \rrvert \leq C
\end{equation}
for all $i \in\pi$ and, as $N \to\infty$,
%
%e2.13 #&#
\begin{equation}
\label{differentgroup} N^{1/2} \bigl(\llvert d_\ell\rrvert -
1\bigr)^{1/2} \llvert d_i - d_\ell \rrvert \to
\infty
\end{equation}
for all $i \in[\![1,r ]\!]\setminus\pi$.

Define the rescaled eigenvalues $\bolds\zeta= (\zeta_i)_{i \in\pi}$
through
%
%e2.14 #&#
\begin{equation}
\label{defzeta} \zeta_i:= N^{1/2} \bigl(\llvert
d_\ell\rrvert - 1\bigr)^{-1/2} \bigl(\mu _{\alpha(i)} -
\theta(d_\ell) \bigr),
\end{equation}
where we recall the definition (\ref{defalphai}) of $\alpha(i)$. Let
$\bolds\xi= (\xi_i)_{i \in\pi}$ denote the eigenvalues of the random
$\llvert \pi\rrvert  \times\llvert  \pi\rrvert $ matrix
%
%e2.15 #&#
\begin{equation}
\label{referencematrix} \Upsilon^\ell_{[\pi]} +
\Psi^\ell_{[\pi]} + N^{1/2} \bigl(\llvert
d_\ell\rrvert - 1\bigr)^{1/2} \bigl(\llvert d_\ell
\rrvert + 1\bigr) \bigl(d_\ell ^{-1} - D_{[\pi]}^{-1}
\bigr).
\end{equation}
Then for any bounded and continuous function $f$ we have
\[
\lim_N \bigl(\mathbb{E}f (\bolds\zeta ) - \mathbb {E}f
(\bolds \xi ) \bigr) = 0.
\]
\end{theorem}

The subset $\pi$ indexes outliers that belong to the same group of
overlapping outliers, as required by (\ref{samegroup}) [see also
(\ref{overlappingcond2}) in the preceding discussion]. As required by
(\ref{differentgroup}), the remaining outliers do not overlap with the
outliers indexed by $\pi$.

%
%re2.6 #&#
\begin{remark}
The reference point $\ell$ for the block $\pi$ is arbitrary and
unimportant. See Lemma~\ref{lemreferencepoint} below and the comment
preceding it for a more detailed discussion.
\end{remark}

%
%re2.7 #&#
\begin{remark}
For the special case $\pi= \{\ell\}$, Theorem~\ref{thmeasystatement}
essentially\footnote{In fact, condition of \cite{KY2} analogous to
(\ref{differentgroup}), equation (2.24) in \cite{KY2}, is slightly
stronger than~(\ref{differentgroup}).} reduces to Theorem~2.14 of
\cite{KY2}. In addition, Theorem~\ref{thmeasystatement} corrects a
minor issue in the statement of Theorem~2.14 of \cite{KY2}, where the
variance of $\Upsilon$ was not necessarily positive. Indeed, in the
language of the current paper, in \cite{KY2} the term $V_\delta^* H
V_\delta$ in~(\ref{defUpsilon}) was of the form $V^* H V$, which
amounted to transferring a large Gaussian component from $\Psi$ to
$\Upsilon$. This transfer was ill-advised as it sometimes resulted in
a~negative variance for $\Psi$ (which would however be compensated in the
sum $\Upsilon+ \Psi$ by a large asymptotically Gaussian component in
$\Upsilon$).
\end{remark}

The functions $\mathcal P$, $\mathcal Q$, $\mathcal R$ and $\mathcal S$
in (\ref{defUpsilon}) and (\ref{defofcov}) are in general nonzero in
the limit $N \to\infty$. They encode the \textit{nonuniversality} of
the distribution of the outliers. Thus, the distribution of the
outliers may depend on the law of the entries of $H$ as well as on the
geometry of the eigenvectors $V$.

In the GOE/GUE case, it is easy to check that $\Upsilon^\ell+
\Psi^\ell$ is asymptotically Gaussian with covariance matrix
%
%e2.16 #&#
\begin{equation}
\label{GUEcovariancematrix} \frac{\llvert  d_\ell\rrvert  + 1}{d_\ell^2} \Delta_{ij,kl}.
\end{equation}
Moreover, if $\lim_N \llvert  d_\ell\rrvert  = 1$ then the matrix
$\Upsilon ^\ell+ \Psi^\ell$ converges weakly to a Gaussian matrix with
covariance given by (\ref{GUEcovariancematrix}). In this case,
therefore, the nonuniversality is washed out. Thus, only outliers
separated from the bulk spectrum $[-2,2]$ by a distance of order one
may exhibit nonuniversality.

If $\lim_N \max_{i,j}\llvert  V_{ij} \rrvert  = 0$, then an appropriate
choice of $\delta$ yields $\Upsilon^\ell= (\llvert  d_\ell\rrvert  + 1)
(\llvert  d_\ell\rrvert  - 1)^{1/2} d_\ell^{-4} \mathcal S(V)$ as well
as a matrix $\Psi^\ell$ whose covariance is asymptotically that of the
GOE/GUE case, that is, (\ref{GUEcovariancematrix}). Hence, in this case
the only manifestation of nonuniversality is the deterministic shift
given by $\Upsilon^\ell$.

It is possible to find scenarios in which each term of
(\ref{defUpsilon}) and (\ref{defofcov}) [apart from the trivial error
term $E$ in (\ref{defofcov})] contributes\vspace*{1pt} in the limit $N \to\infty$.
This is, for instance, the case if $\mu^{(3)}_{ij}$ and $\mu^{(4)}_{ij}$
do not depend\vspace*{-1pt} on $i$ and $j$, $\mu^{(4)}_{ij}$ is not asymptotically $4
- \beta$, and an eigenvector $\mathbf{v}^{(i)}$ satisfies $\llVert
\mathbf{v}^{(i)} \rrVert _\infty\geq c$ as well as $\llVert
\mathbf{v}^{(i)} \rrVert _1 \geq c N^{1/2}$ for~some constant $c> 0$.
We refer to \cite{KY2}, Remarks 2.17--2.21, for analogous remarks,
where more details are given for the case $\pi= \{\ell\}$.\looseness=-1

Next, we give the asymptotic distribution of a group of overlapping
outliers in full generality. Thus, Theorem~\ref{thmmainresult} below
holds for arbitrary sequences $V \equiv V_N$ and $D \equiv D_N$
satisfying $V^*V = \mathbh{1}$ and (\ref{basicconditionsond}).\vspace*{-2pt}

%
%de2.8 #&#
\begin{definition} \label{defpis}
Let $N$ and $D$ be given. For $s > 0$ and $\ell\in[\![1,r
]\!]$ satisfying $\llvert  d_\ell\rrvert  > 1$, define $\pi(\ell,
s) \equiv\pi_{N, D}(\ell,s)$ as the smallest subset of $[\![1,r
]\!]$ with the two following properties.
\begin{longlist}[(ii)]
\item[(i)] $\ell\in\pi(\ell,s)$.

\item[(ii)] If for $i,j \in[\![1,r ]\!]$ we have $\llvert  d_i \rrvert
    >  1$  and
%
%e2.17 #&#
\begin{equation}
\label{condonpi} N^{1/2} \bigl(\llvert d_i \rrvert - 1
\bigr)^{1/2} \llvert d_i- d_j \rrvert \leq s,
\end{equation}
then either $i,j \in\pi(\ell, s)$ or $i,j \in[\![1,r
]\!]\setminus\pi(\ell, s)$.\vspace*{-2pt}
\end{longlist}
\end{definition}

The subset $\pi(\ell, s)$ indexes those outliers that belong to the
same group of overlapping outliers as $\ell$, where $s$ is a cutoff
distance used to determine whether two outliers are considered
overlapping. Note that $\pi(\ell,s)$ is a set of consecutive
integers.\vspace*{-2pt}

%
%th2.9 #&#
\begin{theorem}\label{thmmainresult}
For large enough $K$ the following holds. Let $\varepsilon> 0$ be
arbitrary, and let $f_1, \ldots, f_r$ be bounded continuous functions,
where $f_k$ is a function on $\mathbb{R}^k$. Then there exist $N_0 \in
\mathbb{N}$ and $s_0 > 0$ such that for all $N \geq N_0$ and $s \geq
s_0$ the following holds.

Suppose that $\ell\in[\![1,r ]\!]$ satisfies
%
%e2.18 #&#
\begin{equation}
\label{condondell} \llvert d_\ell\rrvert \geq1 +
\varphi^K N^{-1/3}
\end{equation}
and set $\pi:=\pi(\ell,s)$. Then
%
%e2.19 #&#
\begin{equation}
\label{mainresult} \bigl\llvert \mathbb{E}f_{\llvert  \pi\rrvert } (\bolds\zeta ) -
\mathbb{E}f_{\llvert  \pi\rrvert } (\bolds\xi ) \bigr\rrvert \leq\varepsilon,
\end{equation}
where $\bolds\zeta$ and $\bolds\xi$ were defined
Theorem~\ref{thmeasystatement}.\vspace*{-2pt}
\end{theorem}

%s2.4 #&#
\subsection{The joint distribution} \label{secjointdistrresults}

In order to describe the joint distribution of all outliers, we
organize them into groups of overlapping outliers, using a partition
$\Pi$ whose blocks $\pi$ are defined using the subsets $\pi(\ell,s)$
from Definition~\ref{defpis}.

%
%de2.10 #&#
\begin{definition}\label{defPi}
Let $N$ and $D$ be given, and fix $K > 0$ and $s > 0$. We introduce a
partition\footnote{That $\Pi$ is a partition follows from the
observation that $\ell' \in\pi(\ell,s)$ if and only if $\ell\in
\pi(\ell',s)$. Therefore if $\ell$ and $\ell'$ satisfy $\llvert d_\ell
\rrvert  \geq1 + \varphi^K N^{-2/3}$ and $\llvert  d_{\ell'} \rrvert
\geq1 + \varphi^K N^{-2/3}$ then either $\pi(\ell,s) = \pi(\ell',s)$ or
$\pi(\ell,s) \cap\pi(\ell',s) = \varnothing$.} $\Pi\equiv\Pi(N,K,s,D)$
on a subset of $[\![1,r ]\!]$, defined as
\[
\Pi:= \bigl\{\pi(\ell,s)\dvtx \ell\in[\![1,r ]\!], \llvert
d_\ell\rrvert \geq1 + \varphi^K N^{-1/3} \bigr\}.
\]
We also use the notation $\Pi= \{\pi\}_{\pi\in\Pi}$ and $[\Pi]:=
\bigcup_{\pi\in\Pi} \pi$.\vadjust{\goodbreak}
\end{definition}

The indices in $[\Pi]$ give rise to outliers, which are grouped into
the blocks of $\Pi$. Indices in $[\![1,r
]\!]\setminus[\Pi]$ do not give rise to outliers.

For $\pi\in\Pi$, we define
%
%e2.20 #&#
\begin{equation}
\label{defdpi} d _{\pi}:= \min\{d_i\dvtx i\in\pi\}.
\end{equation}
We chose this value for definiteness, although any other choice of
$d_i$ with $i \in\pi$ would do equally well.

Next, in analogy to (\ref{referencematrix}), we define a $\llvert [\Pi]
\rrvert \times\llvert  [\Pi] \rrvert $ reference matrix whose
eigenvalues will have the same asymptotic distribution as the
appropriately rescaled outliers $(\mu_{\alpha(i)})_{i \in[\Pi]}$.
Define the block diagonal $\llvert  [\Pi] \rrvert  \times\llvert  [\Pi]
\rrvert $ matrix $\Upsilon= \bigoplus_{\pi \in \Pi} \Upsilon^\pi$,
where
\[
\Upsilon^\pi:= \bigl(\llvert d_{\pi} \rrvert + 1\bigr)
\bigl(\llvert d_{\pi} \rrvert - 1\bigr)^{1/2} \biggl(
\frac{N^{1/2} V_\delta^* H V_\delta}{d_{\pi}^2} + \frac{\mathcal S(V)}{d_{\pi}^4} \biggr)_{[\pi]}.
\]
In addition, we introduce a Hermitian, Gaussian $\llvert  [\Pi] \rrvert
\times \llvert  [\Pi] \rrvert $ matrix $\Psi$, that is, independent of
$H$ and whose entries have mean zero. It is block diagonal, \mbox{$\Psi=
\bigoplus_{\pi\in\Pi} \Psi^\pi$,} where the block $\Psi^\pi=
(\Psi^\pi_{ij})_{i,j \in \pi}$ is a $\llvert  \pi\rrvert  \times\llvert
\pi\rrvert $ matrix. The law of $\Psi$ is determined by the covariance
%
%e2.21 #&#
%e2.22 #&#
\begin{eqnarray}\label{defPP}
\mathbb{E}\Psi^\pi_{ij} \Psi^{\pi'}
_{kl} &=& \frac{\llvert
d_{\pi} \rrvert  + 1
}{d_{\pi}^{2} }\delta_{\pi\pi'} \Delta_{ij, kl}
+ \delta_{\pi\pi'} E_{ij,kl}\nonumber
\\
&&{}+ \biggl(\prod_{p={\pi}, {\pi'}}\frac{(\llvert  d_p \rrvert  -
1)^{1/2} (|d_p|+1 )}{d_p^2} \biggr)
\nonumber\\[-8pt]\\[-8pt]
&&\hspace*{11pt}{}\times \biggl(- \mathcal P_{ij,kl}\bigl(V_\delta^* V_\delta
\bigr) + \frac{1}{d_{{\pi}} d_{{\pi'}} } \mathcal R_{ij,kl}(V)\nonumber
\\
&&\hspace*{61pt}{} + \frac{\mathcal W_{ij,kl}(V)}{d_{\pi'}} +
\frac{\mathcal W_{kl,ij}(V)}{d_{\pi}} \biggr),\nonumber
\end{eqnarray}
where we defined
\[
\mathcal W_{ij,kl}(V):= \frac{1}{\sqrt{N}} \sum
_{a,b} \bigl(\widebar V_{ai} \widebar V_{ak}
V_{al} \mu^{(3)}_{ab} V_{bj} + \widebar V_{ia} \mu ^{(3)}_{ab} V_{bj} \widebar V_{bk} V_{bl} \bigr).
\]
(Note that $\mathcal Q_{ij,kl} = \mathcal W_{ij,kl} + \mathcal
W_{kl,ij}$.) As in (\ref{defofcov}), the factor $E_{ij,kl} =
\varphi^{-1} \Delta_{ij,kl}$, whose contribution vanishes in the limit
$N \to\infty$, simply ensures that the right-hand side of (\ref{defPP})
defines a nonnegative matrix; this nonnegativity is an immediate
corollary of our proof in Section~\ref{secgendistrJ}.

Next, in analogy to (\ref{defzeta}), we introduce the rescaled family
of outliers $\bolds\zeta= (\zeta_i^\pi\dvtx \pi\in\Pi, i \in\pi)
\in\mathbb{R}^{[\Pi]}$ whose entries are defined by
%
%e2.23 #&#
\begin{equation}
\label{defzetaJ} \zeta^\pi_{i}:= N^{1/2}
\bigl(\llvert d_{\pi} \rrvert - 1\bigr)^{-1/2} \bigl(
\mu_{\alpha(i)} - \theta(d_{\pi} ) \bigr),
\end{equation}
where we recall the definition (\ref{defalphai}) of $\alpha(i)$.
Moreover, for $\pi\in\Pi$ let $\bolds\xi^\pi= (\xi_i^\pi\dvtx i \in
\pi)$ denote the eigenvalues of the random $\llvert  \pi\rrvert
\times\llvert \pi\rrvert $ matrix
\[
\Upsilon^\pi +\Psi^\pi+ N^{1/2} \bigl(\llvert
d_{\pi} \rrvert - 1\bigr)^{1/2} (|d_{\pi}|+1 )
\bigl(d_{\pi}^{-1} - D_{[\pi]}^{-1} \bigr)
\]
and write $\bolds\xi= (\bolds\xi^\pi\dvtx \pi\in\Pi) = (\xi_i^\pi\dvtx
\pi\in\Pi, i \in\pi) \in\mathbb{R}^{[\Pi]}$. We may now state our main
result in its greatest generality.

%
%th2.11 #&#
\begin{theorem}\label{thmmainresultJ}
For large enough $K$ the following holds. Let $\varepsilon> 0$ be
arbitrary, and let $f_1, \ldots, f_r$ be bounded continuous functions,
where $f_k$ is a function on $\mathbb{R}^k$. Then there exist $N_0 \in
\mathbb{N}$ and $s_0 > 0$ such that for all $N \geq N_0$ and $s \geq
s_0$ we have
\[
\bigl\llvert \mathbb{E}f_{|[\Pi]|} (\bolds\zeta ) - \mathbb
{E}f_{|[\Pi]|} (\bolds\xi ) \bigr\rrvert \leq \varepsilon.
\]
\end{theorem}

We conclude this section by drawing some consequences from
Theorem~\ref{thmmainresultJ}. In the GOE/GUE case, it is easy to see
that the law of the block matrix $\Upsilon+ \Psi$ is asymptotically
Gaussian with covariance
\[
\frac{\llvert  d_{\pi} \rrvert  + 1 }{d_{\pi}^{2} }\delta_{\pi\pi
'} \Delta _{ij, kl}.
\]
In particular, we find that overlapping outliers repel each other
according to the usual random matrix level repulsion, while
nonoverlapping outliers are asymptotically independent.

In general outliers are not asymptotically independent, even if they do
not overlap. Such correlations arise from correlations between
different blocks of $\Upsilon+ \Psi$. There are two possible sources
for these correlations: the term $V_\delta^* H V_\delta$ in the
definition of $\Upsilon$, and the terms $\mathcal R$ and $\mathcal W$
in the covariance (\ref{defPP}) of the Gaussian matrix $\Psi$. Thus,
two outliers may be strongly correlated even if they are located on
opposite sides of the bulk spectrum.

%s3 #&#
\section{Tools} \label{sectools}
The rest of this paper is devoted to the proofs of Theorems
\ref{thmeasystatement}, \ref{thmmainresult} and \ref{thmmainresultJ}.
Sections~\ref{sectools}--\ref{secconcofproof} are devoted to the proof
of Theorem~\ref{thmmainresult}; Theorem~\ref{thmeasystatement} is an
easy corollary of Theorem~\ref{thmmainresult}. Finally,
Theorem~\ref{thmmainresultJ} is proved in Section~\ref{secjointdist} by
an extension of the arguments of
Sections~\ref{sectools}--\ref{secconcofproof}.

We begin with a preliminary section that collects tools we shall use in
the proof. We introduce the spectral parameter
\[
z = E + \ii\eta,
\]
which will be used as the argument of Stieltjes transforms and
resolvents. In the following, we often use the notation $E = \operatorname{Re}z$ and
$\eta= \operatorname{Im}z$ without further comment. Let
\[
\varrho(x):= \frac{1}{2 \pi} \sqrt{\bigl[4 - x^2\bigr]_+}
\qquad(x \in\mathbb{R})
\]
denote the density of the local semicircle law, and
%
%e3.1 #&#
\begin{equation}
\label{definitionofmsc} m(z):= \int\frac{\varrho(x)}{x - z}\,\mathrm{d}x \qquad
\bigl(z \notin[-2,2]\bigr)
\end{equation}
its Stieltjes transform. It is well known that the Stieltjes transform
$m$ satisfies the identity
%
%e3.2 #&#
\begin{equation}
\label{identityformsc} m(z) + \frac{1}{m(z)} + z = 0.
\end{equation}
It is easy to see that (\ref{identityformsc}) and the definition
(\ref{deftheta}) imply
%
%e3.3 #&#
\begin{equation}
\label{mtheta} m\bigl(\theta(d)\bigr) = -\frac{1}{d}.
\end{equation}
For $E \in\mathbb{R}$, define
%
%e3.4 #&#
\begin{equation}
\label{defkappa} \kappa_E:= \bigl\llvert \llvert E \rrvert - 2
\bigr\rrvert,
\end{equation}
the distance from $E$ to the spectral edges $\pm2$. We have the simple
estimate
%
%e3.5 #&#
\begin{equation}
\label{kappad} \kappa_{\theta(d)} \asymp\bigl(\llvert d \rrvert - 1
\bigr)^2
\end{equation}
for $\llvert  d \rrvert  > 1$. The following lemma collects some useful
properties
of $m$.

%
%le3.1 #&#
\begin{lemma} \label{lemmamsc}
For $\llvert  z \rrvert  \leq2 \Sigma$, we have
%
%e3.6 #&#
\begin{equation}
\label{boundsonmsc} \bigl\llvert m(z) \bigr\rrvert \asymp1, \qquad\bigl
\llvert 1 - m(z)^2 \bigr\rrvert \asymp \sqrt{\kappa+ \eta}.
\end{equation}
Moreover,
\[
\operatorname{Im}m(z) \asymp \cases{ \displaystyle\sqrt{\kappa+ \eta}, &\quad if $\llvert E
\rrvert \leq2$, \vspace*{2pt}
\cr
\displaystyle\frac{\eta}{\sqrt{\kappa+ \eta}}, & \quad if $\llvert E \rrvert \geq2$.}
\]
(Here the implicit constants depend on $\Sigma$.)
\end{lemma}

\begin{pf}
The proof is an elementary calculation; see Lemma~4.2 in \cite{EYY2}.
\end{pf}

The following definition introduces a notion of high probability that
is suitable for our needs.

%
%de3.2 #&#
\begin{definition}[(High probability events)]
We say that an $N$-dependent event $\Xi$ holds with \textit{high
probability} if there is some constant $C$ such that
%
%e3.7 #&#
\begin{equation}
\label{highprob} \mathbb{P}\bigl(\Xi^c\bigr) \leq N^C
\exp(-\varphi)
\end{equation}
for large enough $N$.
\end{definition}

Next, we give the key tool behind the proof of
Theorem~\ref{thmmainresult}: the \textit{Isotropic local semicircle
law}. We use the notation $\mathbf{v} = (v_i)_{i = 1}^N
\in\mathbb{C}^N$ for the components of a~vector. We introduce the
standard scalar product $\langle\mathbf{v},
{\mathbf{w}}\rangle:=\sum_i \bar v_i w_i$. For $\eta> 0$, we define
the resolvent of $H$ through
\[
G(z):= (H - z)^{-1}.
\]
The following result was proved in \cite{KY2}, Theorem~2.3.

%
%th3.3 #&#
\begin{theorem}[(Isotropic local semicircle law outside of the spectrum)]
\label{theoremstrongestimate} Fix $\Sigma\geq3$. There exists a
constant $C$ such that for large enough $K$ and any deterministic
$\mathbf{v}, \mathbf{w} \in\mathbb{C}^N$ we have with high probability
%
%e3.8 #&#
\begin{equation}
\label{strongestimate} \bigl\llvert \bigl\langle\mathbf{v}, {G(z) \mathbf{w}}
\bigr\rangle- m(z) \langle {\mathbf{v} }, {\mathbf{w}}\rangle \bigr\rrvert \leq
\varphi^{C} \sqrt{\frac{\operatorname{Im}m(z)}{N \eta}} \llVert \mathbf{v} \rrVert
\llVert \mathbf{w} \rrVert
\end{equation}
for all
\[
E \in \bigl[{-\Sigma, -2 - \varphi^{K} N^{-2/3}} \bigr] \cup
\bigl[{2 + \varphi^{K} N^{-2/3}, \Sigma} \bigr]\quad\mbox{and}\quad\eta\in(0, \Sigma].
\]
\end{theorem}

Using (\ref{kappad}) and Lemma~\ref{lemmamsc}, we find that the control
parameter in (\ref{strongestimate}) may be written as
%
%e3.9 #&#
\begin{equation}
\label{controlparam} \sqrt{\frac{\operatorname{Im}m(z)}{N \eta}} \asymp N^{-1/2} (
\kappa_E + \eta)^{-1/4} \leq N^{-1/2}
\kappa_E^{-1/4}.
\end{equation}
The following result provides sharp (up to logarithmic factors) large
deviations bounds on the locations of the outliers.

%
%th3.4 #&#
\begin{theorem}[(Locations of the deformed eigenvalues)]\label{theoremrank-klde}
There exists a constant $C$ such that, for large enough $K$ and under
condition (\ref{basicconditionsond}), we have
%
%e3.10 #&#
\begin{equation}
\label{LDEforoutliers} \bigl\llvert \mu_{\alpha(i)} -
\theta(d_i) \bigr\rrvert \leq\varphi ^{C} N^{-1/2}
\bigl(\llvert d_i \rrvert - 1\bigr)^{1/2}
\end{equation}
with high probability provided that $\llvert  d_i \rrvert  \geq1 +
\varphi^K N^{-1/3}$.
\end{theorem}

\begin{pf}
This was essentially proved in \cite{KY2}, Theorem~2.7, by setting
$\psi = 1$ there; see equation (2.20) of \cite{KY2}. Note that
Theorem~2.7 of \cite{KY2} has slightly stronger assumptions than
Theorem~\ref{theoremrank-klde}, requiring in addition that there be no
eigenvalues $d_j$ of $D$ satisfying $\llvert  \llvert  d_j \rrvert  - 1
\rrvert  < \varphi^K N^{-1/3}$. However, this assumption was only
needed for equation (2.21) of \cite{KY2}, and the proof from Section~6
of \cite{KY2} may be applied verbatim to (\ref{LDEforoutliers}) under
the assumptions of Theorem~\ref{theoremrank-klde}.
\end{pf}

We shall often need to consider minors of $H$, which are the content of
the following definition. It is a convenient extension of
Definition~\ref{defsmallminors}.

%
%de3.5 #&#
\begin{definition}[(Minors and partial expectation)]
\textup{(i)} For $U \subset[\![1,N ]\!]$, we define
\[
H^{(U)}:= H_{[U^c]} = (h_{ij})_{i,j \in U^c},
\]
where $U^c:=[\![1,N ]\!]\setminus U$. Moreover, we define the
resolvent of $H^{(U)}$ through
\[
G^{(U)}(z):= \bigl(H^{(U)} - z\bigr)^{-1}.
\]

(ii) Set
\[
\sum_i^{(U)}:= \sum
_{i\dvtx i \notin U}.
\]
When $U = \{a\}$, we abbreviate $(\{a\})$ by $(a)$ in the above
definitions; similarly, we write $(ab)$ instead of $(\{a,b\})$.

(iii) For $U \subset[\![1,N ]\!]$ define the partial expectation
    $\mathbb{E}_{U}(X):=\mathbb{E}(X | H^{(U)})$.
\end{definition}

Next, we record some basic large deviations estimates from \cite{KY2},
Lemma~3.5.

%
%le3.6 #&#
\begin{lemma}[(Large deviations estimates)] \label{lemmaLDE}
Let $a_1, \ldots, a_N$, $b_1, \ldots, b_M$ be independent random
variables with zero mean and unit variance. Assume that there is a
constant $\vartheta> 0$ such that
%
%e3.11 #&#
\begin{eqnarray}\label{subexponentialdecayforLDE}
\mathbb{P}\bigl(\llvert a_i \rrvert \geq x\bigr) &\leq&
\vartheta^{-1} \exp\bigl(-x^\vartheta\bigr) \qquad(i = 1, \ldots, N),
\nonumber\\[-8pt]\\[-8pt]
\mathbb{P}\bigl(\llvert b_i
\rrvert \geq x\bigr) &\leq&\vartheta^{-1} \exp\bigl(-x^\vartheta
\bigr) \qquad(i = 1, \ldots, M).\nonumber
\end{eqnarray}
Then there exists a constant $\rho\equiv\rho(\vartheta) > 1$ such
that, for any $\xi> 0$ and any deterministic complex numbers $A_i$ and
$B_{ij}$, we have with high probability
%
%e3.12 #&#
%e3.13 #&#
%e3.14 #&#
\begin{eqnarray}
\label{diagLDE} \biggl\llvert \sum_i
A_i \llvert a_i \rrvert ^2 - \sum
_i A_i \biggr\rrvert &\leq&\varphi
^{\rho\xi} \biggl(\sum_i \llvert
A_i \rrvert ^2 \biggr)^{1/2},
\\
\label{offdiagLDE} \biggl\llvert \sum_{i \neq j} \bar a_i B_{ij} a_j \biggr\rrvert &\leq&\varphi
^{\rho
\xi} \biggl(\sum_{i \neq j} \llvert
B_{ij} \rrvert ^2 \biggr)^{1/2},
\\
\label{two-setLDE} \biggl\llvert \sum_{i,j}
a_i B_{ij} b_j \biggr\rrvert &\leq&
\varphi^{\rho
\xi} \biggl(\sum_{i,j} \llvert
B_{ij} \rrvert ^2 \biggr)^{1/2}.
\end{eqnarray}
\end{lemma}

We conclude this preliminary section by quoting a result on the
eigenvalue rigidity of $H$. Denote by $\gamma_1 \leq\gamma_2 \leq
\cdots\leq\gamma_N$ the classical locations of the eigenvalues of
$H$, defined through
%
%e3.15 #&#
\begin{equation}
\label{classicallocation} N \int_{-\infty}^{\gamma_\alpha}
\varrho(x)\,\mathrm{d}x = \alpha \qquad(1 \leq\alpha\leq N).
\end{equation}
The following result was proved in \cite{EYY3}, Theorem~2.2.

%
%th3.7 #&#
\begin{theorem}[(Rigidity of eigenvalues)] \label{theoremrigidity}
There exists a constant $C$ such that we have with high probability
\[
\llvert \lambda_\alpha- \gamma_\alpha\rrvert \leq\varphi
^{C} \bigl(\min\{\alpha, N + 1 - \alpha\} \bigr)^{-1/3}
N^{-2/3}
\]
for all $\alpha\in[\![1, N ]\!]$.
\end{theorem}

%s4 #&#
\section{Coarser grouping of outliers and reduction to the law of
$G$}\label{secreduction2} For the following, we fix the sequences
$(V_N)_N$ and $(D_N)_N$. It will sometimes be convenient to assume that
%
%e4.1 #&#
\begin{equation}
\label{limdexists} \lim_N d_i^{(N)}
\qquad\mbox{exists for all } i \in[\![1,r ]\!].
\end{equation}
To that end, we invoke the following elementary result.

%
%le4.1 #&#
\begin{lemma} \label{lemsequencetrick}
Let $(a_N)_{N}$ be a sequence of nonnegative numbers and $\varepsilon>
0$. The following statements are equivalent.
\begin{longlist}[(ii)]
\item[(i)] $a_N \leq\varepsilon$ for large enough $N$.

\item[(ii)] Each subsequence has a further subsequence along which $a_N
    \leq\varepsilon$.
\end{longlist}
\end{lemma}

We use Lemma~\ref{lemsequencetrick} by setting $a_N$ to be the
left-hand side of (\ref{mainresult}). Using
Lemma~\ref{lemsequencetrick}, we therefore find that
Theorem~\ref{thmmainresult} holds for arbitrary~$D$ if it holds for~$D$
satisfying (\ref{limdexists}). From now on, we therefore assume without
loss of generality that (\ref{limdexists}) holds.

For the proof of Theorem~\ref{thmmainresult}, we need a new subset of
$[\![1,r ]\!]$, denoted by $\gamma(\ell)$, which is larger than or
equal to the subset $\pi(\ell,s)$ from Definition~\ref{defpis}.

%
%de4.2 #&#
\begin{definition} \label{defgamma}
For $\ell\in[\![1,r ]\!]$ satisfying (\ref{condondell}), define
$\gamma(\ell) \equiv\gamma_{N,D,K}(\ell)$ as the smallest subset of
$[\![1,r ]\!]$ with the two following properties.
\begin{longlist}[(ii)]
\item[(i)] $\ell\in\gamma(\ell)$.

\item[(ii)] If for $i,j \in[\![1,r ]\!]$ we have $\llvert  d_i \rrvert
    > 1$ and
%
%e4.2 #&#
\begin{equation}
\label{condforgamma} N^{1/2} \bigl(\llvert d_i \rrvert -
1\bigr)^{1/2} \llvert d_i- d_j \rrvert \leq
\varphi ^{K/2},
\end{equation}
then either $i,j \in\gamma(\ell)$ or $i,j \in\bar\gamma(\ell)$.
\end{longlist}
Here we use the notation $\bar\gamma(\ell):=[\![1,r ]\!]
\setminus\gamma(\ell)$.
\end{definition}

Note that $\gamma(\ell)$ is a set of consecutive integers. Similar to
$\pi(\ell,s)$, the set $\gamma(\ell)$ indexes outliers that are close
to that indexed by $\ell$, except that now the threshold used to
determine whether two outliers overlap is larger ($\varphi^{K/2}$
instead of the \mbox{$N$-}independent $s$). This need to regroup outliers into
larger subsets arises from the perturbation theory argument in
Proposition~\ref{propositionreduction} below. At the end of the proof,
in Section~\ref{secconcofproof}, we shall use perturbation theory a
second time to obtain a statement involving outliers in $\pi(\ell,s)$
instead of $\gamma(\ell)$.

For the following, we introduce the abbreviation
\[
\delta_\rho(d):= \varphi^\rho N^{-1/2} \bigl(
\llvert d \rrvert - 1\bigr)^{-1/2},
\]
so that (\ref{condforgamma}) reads $\llvert  d_i - d_j \rrvert  \leq
\delta_{K/2}(d_i)$. We have the following elementary result.

%
%le4.3 #&#
\begin{lemma} \label{lemmaddprime}
Let $\rho> 0$. If $\llvert  d \rrvert  \geq1 + \varphi^\rho N^{-1/3}$
and $\llvert  d - d' \rrvert  \leq\delta_{\rho}(d)$, then
\[
\bigl\llvert d' \bigr\rrvert - 1 = \bigl(\llvert d \rrvert - 1
\bigr) \bigl(1 + O\bigl(\varphi^{-\rho
/2}\bigr) \bigr).
\]
\end{lemma}

For brevity, we fix $\ell$ satisfying (\ref{condondell}), and
abbreviate $\gamma\equiv\gamma(\ell)$ and
$\bar\gamma\equiv\bar\gamma(\ell)$ when there is no risk of confusion.
The indices of $\gamma$ and $\bar\gamma$ are separated in the following
sense.

%
%le4.4 #&#
\begin{lemma}
If $i \in\gamma$ and $j \in\bar\gamma$, then
%
%e4.3 #&#
\begin{equation}
\label{differentblock} \llvert d_i - d_j \rrvert >
\delta_{K/2}(d_i).
\end{equation}
If $i,j \in\gamma$, then
%
%e4.4 #&#
\begin{equation}
\label{sameblock} \llvert d_i - d_j \rrvert \leq2 r
\delta_{K/2}(d_i).
\end{equation}
\end{lemma}

\begin{pf}
The bound (\ref{differentblock}) follows immediately from the
definition of $\gamma$. The bound (\ref{sameblock}) follows immediately
from Lemma~\ref{lemmaddprime} and the fact that $\gamma$ is a set of at
most $r$ consecutive integers.
\end{pf}

Since $D$ is diagonal, we may write
\[
D = D_{[\gamma]} \oplus D_{[\bar\gamma]}.
\]
The matrix $D_{[\gamma]}$ has dimensions $\llvert  \gamma\rrvert
\times
\llvert  \gamma\rrvert $ and eigenvalues $(d_i)_{i \in\gamma}$.
Define the region
%
%e4.5 #&#
\begin{equation}
\label{defB} \mathcal B:= \Bigl[\min_{i \in\gamma}
\bigl(d_i - \delta _{K/4}(d_i) \bigr), \max
_{i \in\gamma} \bigl(d_i + \delta _{K/4}(d_i)
\bigr) \Bigr].
\end{equation}
From (\ref{condondell}), (\ref{sameblock}) and Lemma~\ref{lemmaddprime}
we get, for any $i \in\gamma$, that
\begin{eqnarray*}
\llvert d_i \rrvert - \delta_{K/4}(d_i) &
\geq&\llvert d_\ell \rrvert - \llvert d_i - d_\ell
\rrvert - 2 \varphi^{K/4} N^{-1/2} \bigl(\llvert d_\ell
\rrvert - 1\bigr)^{-1/2}
\\
& \geq&1 + \varphi^K N^{-1/3} - (2r + 2)
\varphi^{K/2}N^{-1/2} \bigl(\llvert d_\ell\rrvert - 1
\bigr)^{-1/2}
\\
& \geq&1 + \varphi^K N^{-1/3} - (2 r + 2) N^{-1/3}
\\
& > &1.
\end{eqnarray*}
We therefore conclude that $\mathcal B
\subset\mathbb{R}\setminus[-1,1]$. For large enough $K$ a simple
estimate using the definition of $\theta$ and the bound
(\ref{LDEforoutliers}) yields for all $i \in\gamma$
%
%e4.6 #&#
\begin{equation}
\label{roughlocationofq-group} \sigma(\widetilde H) \cap\theta(\mathcal B) = \{
\mu_{\alpha(i)}\}_{i
\in\gamma}
\end{equation}
with high probability. In other words, $\theta(\mathcal B)$ houses with
high probability all of the outliers indexed by $\gamma$, and no other
eigenvalues of $\widetilde H$. Moreover, from
Theorem~\ref{theoremrigidity} we find that for large enough $K$ the
region $\theta(\mathcal B)$ contains with high probability no
eigenvalues of $H$.

We may now state the main result of this section. Introduce the $r
\times r$ matrix
\[
M(z):= V^* G(z) V.
\]
To shorten notation, for $i$ satisfying $\llvert  d_i \rrvert  > 1$ we often
abbreviate
\[
\theta_i:= \theta(d_i).
\]

%
%pr4.5 #&#
\begin{proposition} \label{propositionreduction}
The following holds for large enough $K$. Let $\ell\in[\![1,r ]\!]$
satisfy (\ref{condondell}), and write $\gamma\equiv\gamma(\ell)$. Then
for all $i \in\gamma$ we have
%
%e4.7 #&#
\begin{equation}
\label{claimofreduction}
\quad\qquad\biggl\llvert \mu_{\alpha(i)} -
\lambda_i \biggl(\theta_{\ell} - \frac
{1}{m'(\theta_{\ell})}
\bigl(M(\theta_{\ell}) + D^{-1} \bigr)_{[\gamma]}
\biggr) \biggr\rrvert \leq\varphi^{-1} N^{-1/2} \bigl(\llvert
d_{\ell} \rrvert - 1\bigr)^{1/2}
\end{equation}
with high probability. [Recall Definitions \ref{defeigenvalues} and
\ref{defsmallminors} for the meaning of $\lambda_i(\cdot)$ on the
left-hand side.]
\end{proposition}

\begin{pf}
%We have to introduce some additional randomness in order to (almost
%surely) avoid pathological coincidences.
Our strategy for locating the outliers is based on the well-known fact
that $x \notin\sigma(H)$ is an eigenvalue of $\widetilde H$ if and only
if $M(x) + D^{-1}$ has a zero eigenvalue (see, e.g., Lemma~6.1 of
\cite{KY2}). Below, we develop a counting argument that finds the
eigenvalues of $\widetilde H$ by analysing the behavior of each
eigenvalue of $M(x) + D^{-1}$ as $x$ varies. For our argument to work,
it is important that no two eigenvalues of $M(x) + D^{-1}$
simultaneously cross the origin. [This condition is made precise in the
claim (\ref{eqzvai}) below.] In order to rule out such coincidences, we
introduce additional randomness, by adding a small perturbation
$\varepsilon\Delta$, where $\Delta$ has an absolutely continuous law.
The sole purpose of this perturbation is to exclude these coincidences
almost surely in the randomness of $\Delta$. This perturbation is
purely \textit{qualitative} in the sense that $\varepsilon> 0$ may be
arbitrarily small; once the counting argument is concluded, we may
easily take $\varepsilon\to0$ and recover the claim for $\varepsilon=
0$ by a trivial continuity argument.

Thus, let $\Delta$ be an $r \times r$ Hermitian random matrix whose
upper-triangular entries are independent and have an absolutely
continuous law supported in the unit disc. Moreover, let $\Delta$ be
independent of $H$. Let $\varepsilon> 0$. We shall prove the claim of
Proposition~\ref{propositionreduction} for the matrix $\widetilde
H{}^\varepsilon:=H + V (D^{-1} + \varepsilon\Delta)^{-1} V^*$ for small
enough $\varepsilon$ (depending on $N$), instead of $\widetilde H = H +
V D V^*$.
%Having done this, the claim for $\widetilde H$ follows easily by taking the
%limit $\varepsilon\to0$.

Define the $r \times r$ matrix
%
%e4.8 #&#
\begin{equation}
\label{defA} A^\varepsilon(x):= M(x) - m(x) + D^{-1} +
\varepsilon\Delta.
\end{equation}
From \cite{KY2}, Lemma~6.1, we get that $x \notin\sigma(H)$ is an
eigenvalue of $\widetilde H{}^\varepsilon$ if and only if
$A^\varepsilon(x) + m(x)$ has a zero eigenvalue. Similar to
Proposition~7.1 in \cite{KY2}, we use perturbation theory to compare
the eigenvalues of $A^\varepsilon(x)$ with those of the block matrix
\[
\tilde A^\varepsilon(x):= A^\varepsilon_{[\gamma]}(x) \oplus
A^\varepsilon_{[\bar\gamma]}(x).
\]

In order to apply perturbation theory, we must establish a lower bound
on the spectral gap
\[
\operatorname{dist} \bigl(\sigma \bigl(A_{[\gamma]}^\varepsilon(\theta_{\ell
})
\bigr), \sigma \bigl(A_{[\bar\gamma]}^\varepsilon(\theta _{\ell})
\bigr) \bigr).
\]
We find, for large enough $K$ and small enough $\varepsilon$ (depending
on $N$), that with high probability
%
%e4.9 #&#
\begin{eqnarray} \label{lowerboundforpert}
&& \operatorname{dist}
\bigl(\sigma \bigl(A_{[\gamma]}^\varepsilon( \theta_{\ell }) \bigr),
\sigma \bigl(A_{[\bar\gamma]}^\varepsilon( \theta _{\ell}) \bigr)
\bigr)\nonumber
\\
&&\qquad  \geq \operatorname{dist} \bigl(\sigma
\bigl(D_{[\gamma]}^{-1}\bigr), \sigma\bigl(D_{[\bar\gamma ]}^{-1}
\bigr) \bigr) - \delta_C(d_{\ell}) - r \varepsilon
\\
&&\qquad \geq c \delta_{K/2}(d_{\ell}) - \delta_C(d_{\ell}) \geq\delta_{K/2 -
1}(d_{\ell});\nonumber
\end{eqnarray}
in the first step we used Lemma~\ref{lemmagrowthofspectrum}, $\llVert
\varepsilon\Delta\rrVert \leq r \varepsilon$ and
%
%e4.10 #&#
\begin{equation}
\label{M-mnorm} \bigl\llVert M(\theta_\ell) - m(
\theta_\ell) \bigr\rrVert \leq\delta _C(d_{\ell})
\end{equation}
by Theorem~\ref{theoremstrongestimate}, (\ref{kappad}),
(\ref{controlparam}) and (\ref{condondell}); in the second step we used
(\ref{differentblock}) and chose~$\varepsilon$ to be small enough
(depending on $N$); in the last step we chose $K$ to be large enough
(depending on $C$).

Next, Theorem~\ref{theoremstrongestimate}, (\ref{kappad}) and
(\ref{controlparam}) yield, with high probability,
%
%e4.11 #&#
\begin{equation}
\label{upperboundforpert} \bigl\llVert A^\varepsilon(\theta_{\ell})
- \tilde A^\varepsilon(\theta _{\ell}) \bigr\rrVert \leq
\delta_{K/4 - 2}(d_{\ell})
\end{equation}
for large enough $K$ and small enough $\varepsilon$ (depending on $N$).

Define the regions
\begin{eqnarray*}
\mathcal D &:=& \bigcup_{i \in\gamma} \bigl[{d_i^{-1}
- \delta _{K/ 4}(d_{\ell}), d_i^{-1} +
\delta_{K / 4}(d_{\ell})} \bigr],
\\
\widebar{\mathcal D} &:=&
\bigcup_{i \in\bar\gamma} \bigl[{d_i^{-1}
- \delta_{K / 4}(d_{\ell}), d_i^{-1} +
\delta_{K / 4}(d_{\ell
})} \bigr],
\end{eqnarray*}
which are disjoint by (\ref{differentblock}). Using (\ref{M-mnorm}), we
find that for large enough $K$ and small enough $\varepsilon$
(depending on $N$) we have, with high probability,
\[
\sigma \bigl(A_{[\gamma]}^\varepsilon(\theta_{\ell}) \bigr)
\subset\mathcal D, \qquad\sigma \bigl(A^\varepsilon(\theta _{\ell})
\bigr) \subset\mathcal D \cup\widebar{\mathcal D}.
\]
Moreover, both $A^\varepsilon(\theta_{\ell})$ and
$A_{[\gamma]}^\varepsilon(\theta_{\ell})$ have exactly $\llvert
\gamma\rrvert $ eigenvalues in $\mathcal D$; we denote these
eigenvalues by $(a_i^\varepsilon)_{i \in\gamma}$ and $(\tilde
a_i^\varepsilon)_{i \in \gamma}$, respectively.

We may now apply perturbation theory. Invoking
Proposition~\ref{propperturbation} using (\ref{lowerboundforpert}) and
(\ref{upperboundforpert}) yields with high probability
%
%e4.12 #&#
\begin{equation}
\label{boundfromperturbation} a_i^\varepsilon= \tilde a^\varepsilon_i + O \biggl(\frac{\delta
_{K/4 - 2}(d_{\ell})^2}{\delta_{K/2 - 1}(d_{\ell})} \biggr) = \tilde a^\varepsilon_i + O \bigl(\delta_{-3}(d_{\ell})\bigr)
\end{equation}
for $i \in\gamma$.

Next, we allow the argument $x$ of $A^\varepsilon(x)$ to vary in order
to locate the eigenvalues of $\widetilde H{}^\varepsilon$. We recall
the following derivative bound from \cite{KY2}, Lem\-ma~7.2: there is a
constant $C$ such that for large enough $K$ we have for all
$\ell^2$-normalised $\mathbf{v}, \mathbf{w} \in\mathbb{C}^N$, with high
probability,
%
%e4.13 #&#
\begin{eqnarray}\label{derivativeofG}
\bigl\llvert \partial_x
G_{\mathbf{v} \mathbf{w}}(x) - \partial_x m(x) \langle\mathbf{v},
\mathbf{w}\rangle\bigr\rrvert \leq\varphi^C N^{-1/3} \kappa
_x^{-1}
\nonumber\\[-8pt]\\[-8pt]
\eqntext{\mbox{for } x \in \bigl[{-\Sigma, -2 -
\varphi^{K/2} N^{-2/3}} \bigr] \cup \bigl[{2 +
\varphi^{K/2} N^{-1/3}, \Sigma} \bigr].}
\end{eqnarray}
By the definition (\ref{defB}) of $\mathcal B$, we find from
Lemma~\ref{lemmaddprime}, (\ref{condondell}) and (\ref{sameblock}) that
%
%e4.14 #&#
\begin{equation}
\label{positionofxinBq}
\qquad x \in\theta(\mathcal B) \quad\Longrightarrow\quad
\theta \bigl(d_{\ell} - 3 r \delta_{K/2}(d_{\ell})
\bigr) \leq x \leq \theta \bigl(d_{\ell} + 3 r \delta_{K/2}(d_{\ell})
\bigr).
\end{equation}
We deduce using Lemma~\ref{lemmaddprime}, (\ref{condondell}) and
(\ref{kappad}) that
%
%e4.15 #&#
\begin{equation}
\label{kappaxequivd-1} \kappa_x \asymp\bigl(\llvert
d_{\ell} \rrvert - 1\bigr)^2\qquad\mbox {for } x \in\theta(
\mathcal B).
\end{equation}
Therefore from Theorem~\ref{theoremstrongestimate}, we conclude with
high probability
%
%e4.16 #&#
\begin{equation}
\label{Mm} M(x) = m(x) + O \bigl(\delta_C(d_{\ell})
\bigr) \qquad\mbox {for } x \in\theta(\mathcal B).
\end{equation}
Similarly, from (\ref{derivativeofG}) we get with high probability
%
%e4.17 #&#
\begin{equation}
\label{Mmprime} M'(x) = m'(x) + O \bigl(
\varphi^C N^{-1/3} \bigl(\llvert d_{\ell} \rrvert -
1\bigr)^{-2} \bigr) \qquad\mbox{for } x \in\theta(\mathcal B).
\end{equation}

With these preliminary bounds, we may vary $x \in\theta(\mathcal B)$.
Let $(a_i(x))_{i \in\gamma}$ denote the continuous family of
eigenvalues of $A^\varepsilon(x)$ satisfying
$a_i^\varepsilon(\theta_{\ell}) = a_i^\varepsilon$ for $i \in \gamma$.
For the following argument, it is helpful to keep
Figure~\ref{figureevs} in mind. We make the following claim:
{\renewcommand{\theequation}{\mbox{$\ast$}}
\begin{eqnarray}\label{eqzvai}
&& \mbox{Almost surely, for all }x \in\theta(\mathcal B)\mbox{ we have that }
\nonumber\\[-10pt]\\[-10pt]
&&\qquad a_i^\varepsilon(x) = - m(x)\mbox{ for at most one }i \in
\gamma.\nonumber
\end{eqnarray}\setcounter{equation}{17}}%
We omit the details of the proof\,\footnote{The claim (\ref{eqzvai})
reduces to the following statement. Let $B(x)$ with $x \in I$ and
$\Delta$ be Hermitian matrices such that $B(x)$ is deterministic and
depends smoothly on $x$, and $\Delta$ has an absolutely continuous law;
then, almost surely in $\Delta$, for all $x \in I$ the matrix $B(x) +
\Delta$ has at most one zero eigenvalue. Let $S$ denote the subset of
matrices with multiple eigenvalues at zero, so that $S$ is an algebraic
variety of codimension two. The claim therefore reduces to the
statement that the path $\{B(x)\}_{x \in I} + \Delta$ almost surely
does not intersect $S$, which is standard.} of (\ref{eqzvai}). Note
that the necessity for (\ref{eqzvai}) to hold is the only reason we had
to introduce the additional randomness $\Delta$ into $\widetilde
H{}^\varepsilon$.
%
%f2 #&#
\begin{figure}%[ht!]

\includegraphics{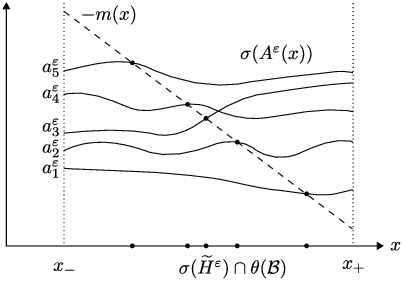}

\caption{The spectrum of $A^\varepsilon(x)$ for $x \in\theta
(\mathcal B)$. For definiteness, we chose $\gamma= [\![1,5 ]\!]$. The
region $x \in\theta(\mathcal B)$ is delimited by dotted lines. The
eigenvalues of $\widetilde H{}^\varepsilon$ are labelled by black dots on the
$x$-axis.}\label{figureevs}
\end{figure}

For definiteness, suppose for the following that $d_\ell> 1$. We claim
that for all $i \in\gamma$ we have with high probability
%
%e4.18 #&#
\begin{equation}
\label{amaorder} a_i^\varepsilon(x_-) \leq-m(x_-), \qquad-m(x_+) \leq a_i^\varepsilon(x_+),
\end{equation}
where $x_\pm$ denote the endpoints of the interval $\theta(\mathcal
B)$. Let us focus on the first estimate; the second one is proved
similarly. Let $i:=\min\gamma$. Since $d \mapsto d - \delta_{K/4}(d)$
is increasing, we find that the left endpoint of $\mathcal B$ is $d_i -
\delta_{K/4}(d_i)$. From (\ref{Mm}) and
Lemma~\ref{lemmagrowthofspectrum}, we find with high probability
\begin{eqnarray*}
\max_{x \in\theta(\mathcal B)} \max_{j \in\gamma} a_j^\varepsilon
(x) &\leq& \frac{1}{d_i} + \delta_C(d_\ell) + r
\varepsilon
\\
&\leq&\frac{1}{d_i - \delta_{K/4}(d_i)} - c \delta_{K/4}(d_i) +
\delta_C(d_\ell) + r \varepsilon
\\
&\leq& \frac{1}{d_i - \delta _{K/4}(d_i)} = -m(x_-);
\end{eqnarray*}
in the second step we used $1 \leq d_i \leq\Sigma- 1$; the third step
holds for large enough $K$ and small enough $\varepsilon$ (depending on
$N$), by Lemma~\ref{lemmaddprime}; the last step follows from~(\ref{mtheta}). This concludes the proof of (\ref{amaorder}).

Recall that $\widetilde H{}^\varepsilon$ has with high probability
exactly $\llvert  \gamma\rrvert $ eigenvalues in $\theta(\mathcal B)$.
By continuity of $a_i^\varepsilon(x)$ the property (\ref{eqzvai}) and
(\ref{amaorder}), we therefore get that the function~$-m(x)$ intersects
each function $a^\varepsilon_i(x)$, $i \in\gamma$, exactly once in
$\theta(\mathcal B)$. Let $i \in\gamma$ and denote by $x_i^\varepsilon$
the unique point (with high probability) in $\theta(\mathcal B)$ at
which $a^\varepsilon_i(x_i^\varepsilon) = -m(x_i^\varepsilon)$.

From the definition of $A^\varepsilon$ and (\ref{Mmprime}) we get, with
high probability,
%
%e4.19 #&#
\begin{eqnarray}\label{mxaeps}
-m\bigl(x_i^\varepsilon\bigr) &=&
a^\varepsilon_i(\theta_{\ell}) + O \bigl(
\varphi^C N^{-1/3} \bigl(\llvert d_{\ell} \rrvert -
1\bigr)^{-2} \bigl\llvert x_i^\varepsilon-
\theta_{\ell} \bigr\rrvert \bigr)
\nonumber\\[-9pt]\\[-9pt]
&=& a_i^\varepsilon+ O
\bigl(\varphi^{K/2 + C} N^{-5/6} \bigl(\llvert d_{\ell}
\rrvert - 1\bigr)^{-3/2} \bigr),\nonumber
\end{eqnarray}
where in the second step we used (\ref{positionofxinBq}), the fact that
$x_i^\varepsilon\in\theta(\mathcal B)$, and the elementary bound
$\llvert  \theta'(d) \rrvert  \asymp\llvert  d \rrvert  - 1$. [Recall
that by definition $a_i^\varepsilon(\theta_{\ell}) = a_i^\varepsilon$.]
Now we may use (\ref{boundfromperturbation}) and (\ref{mxaeps}) to get
%
%e4.20 #&#
\begin{equation}
\label{mxatilde} -m\bigl(x_i^\varepsilon\bigr) = \tilde a_i^\varepsilon+ O \bigl(\delta _{-3}(d_{\ell})
+ \varphi^{K/2 + C} N^{-5/6} \bigl(\llvert d_{\ell} \rrvert
- 1\bigr)^{-3/2} \bigr)
\end{equation}
with high probability. Now we expand the left-hand side using the identity
%
%e4.21 #&#
\begin{equation}
\label{mprime} m' = \frac{m^2}{1 - m^2} \asymp
\kappa_x^{-1/2},
\end{equation}
which follows easily from (\ref{identityformsc}); in the second step we
used Lemma~\ref{lemmamsc}. Differentiating again, we get $m''(x)
\asymp\kappa_x^{-3/2}$. From (\ref{kappaxequivd-1}), we therefore get
%
%e4.22 #&#
\begin{eqnarray}\label{expandingm}
m\bigl(x_i^\varepsilon\bigr) & =& m(\theta_{\ell}) +
m'(\theta_{\ell}) \bigl(x_i^\varepsilon-
\theta_{\ell}\bigr)\nonumber
\\[-1pt]
&&{} + O \bigl(\bigl(\llvert d_{\ell} \rrvert -
1\bigr)^{-3} \bigl(\bigl(\llvert d_{\ell} \rrvert - 1\bigr)
\delta _{K/2}(d_{\ell}) \bigr)^2 \bigr)
\\[-1pt]
& =& m(\theta_{\ell}) + m'(
\theta_{\ell}) \bigl(x_i^\varepsilon-
\theta_{\ell}\bigr) + O \bigl(\varphi^{K} \bigl(\llvert
d_{\ell} \rrvert - 1\bigr)^{-2} N^{-1} \bigr)\nonumber
\end{eqnarray}
with high probability. Solving $x_i^\varepsilon$ from
(\ref{expandingm}) and $-m(x_i^\varepsilon)$ from (\ref{mxatilde}), we
find for large enough $K$ with high probability
\begin{eqnarray*}
x_i^\varepsilon& =& \theta_{\ell} - \frac{1}{m'(\theta_{\ell})}
\bigl(\tilde a_i^\varepsilon+ m(\theta_{\ell}) \bigr)
\\
&&{}+ O \bigl(\varphi^{-3} N^{-1/2} \bigl(\llvert
d_{\ell} \rrvert - 1\bigr)^{1/2}
\\[-1pt]
&&\hspace*{24pt}{} + \varphi^{K/2 + C}
N^{-5/6} \bigl(\llvert d_{\ell} \rrvert - 1\bigr)^{-1/2}
+ \varphi^{K} N^{-1} \bigl(\llvert d_{\ell} \rrvert
- 1\bigr)^{-1} \bigr)
\\
& =& \theta_{\ell} - \frac{1}{m'(\theta_{\ell})} \bigl(\tilde a_i^\varepsilon+
m(\theta_{\ell}) \bigr)
\\
&&{}+ O \bigl(\varphi^{-3} N^{-1/2} \bigl(\llvert
d_{\ell} \rrvert - 1\bigr)^{1/2}
\\[-1pt]
&&\hspace*{24pt}{}+ \varphi^{-K/2 + C}
N^{-1/2} \bigl(\llvert d_{\ell} \rrvert - 1\bigr)^{1/2}
+ \varphi^{-K/2} N^{-1/2} \bigl(\llvert d_{\ell} \rrvert
- 1\bigr)^{1/2} \bigr)
\\
& =& \theta_{\ell} - \frac{1}{m'(\theta_{\ell})} \bigl(\tilde a_i^\varepsilon+
m(\theta_{\ell}) \bigr) + O \bigl(\varphi^{-2}
N^{-1/2} \bigl(\llvert d_{\ell} \rrvert - 1\bigr)^{1/2}
\bigr);
\end{eqnarray*}
in the first step we estimated the error terms using $m'(\theta_\ell)
\asymp(\llvert  d_\ell\rrvert  - 1)^{-1}$ by (\ref{mprime}) and
(\ref{kappaxequivd-1}); in the second step we used (\ref{condondell});
the last step follows by choosing $K$ large enough. Thus, we conclude
that
\[
x_i^\varepsilon= \lambda_i \biggl(
\theta_{\ell} - \frac
{1}{m'(\theta_{\ell})} \bigl(M_{[\gamma]}(
\theta_{\ell}) + D_{[\gamma]}^{-1} + \varepsilon
\Delta_{[\gamma]} \bigr) \biggr) + O \bigl(\varphi^{-2}
N^{-1/2} \bigl(\llvert d_{\ell} \rrvert - 1\bigr)^{1/2}
\bigr)
\]
with high probability for small enough $\varepsilon$ (depending on
$N$). Taking $\varepsilon\to0$ completes the proof.
\end{pf}

We conclude this section with a remark on the choice of the reference
point~$\theta_\ell$ in Proposition~\ref{propositionreduction}. By
definition of $\gamma$, if $i \in\gamma(\ell)$ then $\gamma(i) =
\gamma(\ell)$. Obviously, the distribution of the overlapping group of
outliers $(\mu_{\alpha(i)})_{i \in\gamma}$ cannot depend on the
particular choice of $\ell\in\gamma$. Nevertheless, the reference
matrix $\theta_{\ell} - \frac{1}{m'(\theta_{\ell})}
(M_{[\gamma]}(\theta_{\ell}) + D_{[\gamma]}^{-1} )$ in
(\ref{claimofreduction}) depends explicitly on $\ell\in\gamma$ via
$\theta_{\ell}$. This is not a contradiction, however, since a
different choice of $\ell$ leads to a reference matrix which only
differs from the original one by an error term of order $O
(\varphi^{-1} N^{-1/2} (\llvert d_{\ell} \rrvert  - 1)^{1/2} )$; this
difference may be absorbed into the error term on the right-hand side
of (\ref{claimofreduction}). We shall need this fact in
Section~\ref{secjointdist}. The precise statement is as follows. (To
simplify notation, we state it without loss of generality for the case
$\gamma= [\![1,r ]\!]$.)

%
%le4.6 #&#
\begin{lemma}\label{lemreferencepoint}
Suppose that $\gamma(1) = [\![1,r ]\!]$ and that $\llvert
d_1 \rrvert  \geq1 + \varphi^K N^{-1/3}$. Let
\[
d,\tilde d \in \bigl[{d_1 - \delta_{K/2 + 1}(d_1), d_1 + \delta_{K/2 + 1}(d_1)} \bigr].
\]
Then for large enough $K$ we have
\begin{eqnarray*}
&& \biggl\llVert \biggl(\theta- \frac{1}{m'(\theta)} \bigl(M(\theta) +
D^{-1} \bigr) \biggr) - \biggl(\tilde\theta- \frac
{1}{m'(\tilde\theta)}
\bigl(M(\tilde\theta) + D^{-1} \bigr) \biggr) \biggr\rrVert
\\
&&\qquad \leq \varphi^{-1} N^{-1/2} \bigl(\llvert d_1 \rrvert -
1\bigr)^{1/2},
\end{eqnarray*}
where we abbreviated $\theta\equiv\theta(d)$ and $\tilde\theta
\equiv\theta(\tilde d)$.
\end{lemma}

\begin{pf}
We write
\begin{eqnarray*}
&& \biggl(\theta_i - \frac{1}{m'(\theta)} \bigl(M(\theta) +
D^{-1} \bigr) \biggr) - \biggl(\tilde\theta- \frac
{1}{m'(\tilde\theta)}
\bigl(M(\tilde\theta) + D^{-1} \bigr) \biggr)
\\
&&\qquad = \theta- \tilde\theta+ \frac{1}{m'(\theta)} \bigl(M(\tilde\theta) - M(
\theta) \bigr) + \biggl(\frac
{1}{m'(\tilde\theta)} - \frac{1}{m'(\theta)} \biggr)
\bigl(M(\tilde\theta) + D^{-1} \bigr)
\\
&&\qquad = \theta- \tilde\theta+ \frac{1}{m'(\theta)} \bigl(m(\tilde\theta) - m(
\theta) \bigr) + \biggl(\frac
{1}{m'(\tilde\theta)} - \frac{1}{m'(\theta)} \biggr)
\bigl(m(\tilde\theta) + \tilde d^{-1} \bigr)
\\
&&\quad\qquad{}+ O \bigl(\varphi^{K/2 + C} N^{-5/6} \bigl(\llvert
d_1 \rrvert - 1\bigr)^{-1/2} + \varphi^{K/2 + C}
N^{-1} \bigl(\llvert d_1 \rrvert - 1\bigr)^{-1}
\bigr)
\\
&&\qquad = d + \frac{1}{d} - \tilde d - \frac{1}{\tilde d} +
\bigl(d^2 - 1\bigr) \biggl(\frac{1}{d} - \frac{1}{\tilde d}
\biggr) + O \bigl(\varphi^{-2} N^{-1/2} \bigl(\llvert
d_1 \rrvert - 1\bigr)^{1/2} \bigr)
\\
&&\qquad = O \bigl(\varphi^{-2} N^{-1/2} \bigl(\llvert
d_1 \rrvert - 1\bigr)^{1/2} \bigr)
\end{eqnarray*}
with high probability; in the second step we wrote $M(\tilde\theta
) - M(\theta) = \int_\theta^{\tilde\theta} M'(\xi)\,\mathrm {d}\xi$
and used (\ref{Mmprime}) and Lemma~\ref{lemmaddprime}, as well as
Theorem~\ref{theoremstrongestimate}, (\ref{controlparam}),
(\ref{kappad}), (\ref{mprime}), and the~fact that $m''(x)
\asymp\kappa_x^{-3/2}$; in the third step we used (\ref{deftheta}),
(\ref{mtheta}), and the assumption that $K$ is large enough; in the
last step we used that $(d - \tilde d)^2 \leq4 \varphi^{K + 1}
N^{-1} (\llvert  d_1 \rrvert  - 1)^{-1}$.
\end{pf}

%s5 #&#
\section{The Gaussian case} \label{secGaussian}
Suppose that $\ell$ satisfies (\ref{condondell}). By
Proposition~\ref{propositionreduction}, in order to analyse the joint
distribution of the outliers $(\mu_{\alpha(i)})_{i \in\gamma}$ with
$\gamma\equiv\gamma(\ell)$, it suffices to analyse the distribution of
the eigenvalues of the $\llvert \gamma\rrvert  \times\llvert  \gamma
\rrvert $ matrix $M_{[\gamma]}(\theta_\ell)$. In this section, we do
this under the assumption that the entries of $H$ are Gaussian, that
is, that $H$ is a GOE/GUE matrix.

Recall that $\gamma$ may depend on $N$. To simplify notation, in
Sections~\ref{secGaussian}--\ref{secgencase} we take $\gamma= [\![1,r
]\!]$, which allows us to drop subscripts $[\gamma]$ and avoid minor
nuisances arising from the fact that $\gamma$ may depend on $N$. In
fact, this special case will easily imply the case of general $\gamma$;
see Section~\ref{secconcofproof}.

The following definition is a convenient shorthand for the equivalence
relation defined by two random matrices of fixed size having the same
asymptotic distribution.

%
%de5.1 #&#
\begin{definition}
For two sequences $X_N$ and $Y_N$ of random $k \times k$ matrices,
where $k \in\mathbb{N}$ is fixed, we write $X \stackrel{d}{\sim}Y$ if
\[
\lim_N \bigl(\mathbb{E}f(X_N) -
\mathbb{E}f(Y_N) \bigr) = 0
\]
for all continuous and bounded $f$.
\end{definition}

Let $\Phi= (\Phi_{ij})_{i,j = 1}^r$ be an $r \times r$ GOE/GUE matrix
multiplied by $\sqrt{r}$. In other words, the covariances of $\Phi$ are
given by
%
%e5.1 #&#
\begin{equation}
\label{GUEcovariances} \mathbb{E}\Phi_{ij} \Phi_{kl} =
\Delta_{ij,kl},
\end{equation}
where $\Delta_{ij,kl}$ was defined in (\ref{defDelta}). The following
proposition is the main result of this section. It provides the joint
distribution of the eigenvalues of $M(\theta)$, which, by
Proposition~\ref{propositionreduction}, immediately yields the
distribution of the $\gamma$-group of outliers under the assumption
that $H$ is a GOE/GUE matrix. However, since we are ultimately
interested in non-Gaussian $H$, we shall not combine it
Proposition~\ref{propositionreduction} directly, but instead use it as
an input for the more general case covered in Section~\ref{secAGcase}.

%
%pr5.2 #&#
\begin{proposition} \label{propositionGaussiancase}
The following\vspace*{1pt} holds for large enough $K$. Let $\theta\equiv\theta(d)$
for some $d$ satisfying $\llvert  d \rrvert  \geq1 + \varphi^{K}
N^{-1/3}$. Suppose moreover that $H$ is a GOE/GUE matrix. Then
\[
N^{1/2} \bigl(\llvert d \rrvert - 1\bigr)^{1/2} \bigl(M(
\theta) - m(\theta ) \bigr) \stackrel{d} {\sim} \frac{1}{\llvert  d \rrvert  \sqrt{\llvert  d
\rrvert  + 1}} \Phi.
\]
\end{proposition}

\begin{pf}
Throughout the proof, we drop the spectral parameter $z = \theta$ from
quantities such as $M(\theta)$. By unitary invariance of $H$, we may
assume that \mbox{$V_{ij} = \delta_{ij}$,} that is, $\mathbf{v}^{(i)}$ is the
$i$th standard basis vector of $\mathbb{C}^N$. By Schur's complement
formula, we therefore\vadjust{\goodbreak} get $M = B^{-1}$ where $B = (B_{ij})_{i,j = 1}^r$
is the Hermitian $r \times r$ matrix defined by
\[
B_{ij}:= h_{ij} - \theta- \sum
_{a,b}^{(1 \cdots r)} h_{ia} G^{(1 \cdots r)}_{ab}
h_{bj}.
\]

We now claim that
%
%e5.2 #&#
\begin{equation}
\label{traceGestimate} \Biggl\llvert \frac{1}{N} \sum
_a^{(1 \cdots r)} G^{(1 \cdots r)}_{aa} - m
\Biggr\rrvert \leq\varphi^C N^{-1} \kappa^{-1}_\theta.
\end{equation}
Bearing later applications in mind, we in fact prove, for any $\ell\in
\mathbb{N}$, that
%
%e5.3 #&#
\begin{equation}
\label{Gkm} \biggl\llvert \operatorname{Tr}G^\ell- N \int\frac{\varrho(x)}{(x - \theta)^\ell}\,
\mathrm{d}x \biggr\rrvert \leq\varphi^C \kappa_\theta^{-\ell}
\end{equation}
with high probability. Applying (\ref{Gkm}) with $\ell= 1$ to the minor
$H^{(1 \cdots r)}$ immediately yields (\ref{traceGestimate}). In order
to prove (\ref{Gkm}), we use Theorem~\ref{theoremrigidity} to get with
high probability
%
%e5.4 #&#
%e5.5 #&#
\begin{eqnarray}\label{traceGstep1}
&& \biggl\llvert \sum_\alpha
\frac{1}{(\lambda_\alpha- \theta)^\ell} - \sum_\alpha\frac{1}{(\gamma_\alpha- \theta)^\ell}
\biggr\rrvert\nonumber
\\
&&\qquad \leq \varphi^C \sum_{\alpha= 1}^{N/2}
\frac{\alpha^{-1/3} N^{-2/3}}{(\llvert
\theta\rrvert  - \llvert  \gamma_\alpha\rrvert )^{\ell+1}}
\leq\frac{\varphi^C}{N} \sum_{\alpha= 1}^{N/2}
\frac{(\alpha/N)^{-1/3}}{((\alpha/N)^{2/3} + \kappa_\theta)^{\ell+1}}
\\
&&\qquad \leq\varphi^C \int_0^\infty
\frac{x^{-1/3}}{(x^{2/3} +
\kappa_\theta)^{\ell+1}} \,\mathrm{d}x \leq \frac{\varphi^C}{\kappa_\theta^\ell};\nonumber
\end{eqnarray}
in the first step we estimated the contribution of $\alpha> N/2$ by the
contribution of $N + 1 - \alpha$, and used that $\llvert  \lambda
_\alpha- \gamma_\alpha\rrvert  \ll\llvert  \theta\rrvert  - \llvert
\gamma _\alpha \rrvert $ with high probability by
Theorem~\ref{theoremrigidity} and the assumption on $\theta$ (for large
enough $K$); in the second step we used the estimate
%
%e5.6 #&#
\begin{equation}
\label{behaviourofgamma} 2 - \llvert \gamma_\alpha\rrvert \asymp
\alpha^{2/3} N^{-2/3}
\end{equation}
for $\alpha\leq N/2$, as follows from the definition of
$\gamma_\alpha$. Similarly, setting $\gamma_0:=-2$, we find
%
%e5.7 #&#
%e5.8 #&#
\begin{eqnarray}
\label{traceGstep2}
N \int\frac{\varrho(x)}{(x - \theta)^\ell} \,\mathrm{d}x &=& N \sum
_{\alpha= 1}^N \int_{\gamma_{\alpha- 1}}^{\gamma_\alpha}
\frac{\varrho(x)}{(x - \theta)^\ell} \,\mathrm{d}x\nonumber
\\
&=& \sum_{\alpha= 1}^N \frac{1}{(\gamma_\alpha- \theta)^\ell} + O
\Biggl(\sum_{\alpha= 1}^{N/2} \frac{\alpha^{-1/3}
N^{-2/3}}{(\llvert  \theta\rrvert  - \llvert  \gamma_\alpha\rrvert )^{\ell+ 1}}
\Biggr)
\\
&=& \sum_{\alpha= 1}^N \frac{1}{(\gamma_\alpha- \theta)^\ell}
+ O \biggl(\frac{1}{\kappa_\theta^\ell} \biggr).\nonumber
\end{eqnarray}
Now (\ref{Gkm}) follows from (\ref{traceGstep1}) and\vadjust{\goodbreak}
(\ref{traceGstep2}).

Using $\mathbb{E}h_{ia} h_{bj} = \delta_{ij} \delta_{ab} N^{-1}$ and
(\ref{kappad}), we therefore get from (\ref{traceGestimate})
\begin{eqnarray*}
&& \sum_{a,b}^{(1 \cdots r)} h_{ia}
G^{(1 \cdots r)}_{ab} h_{bj} - \delta _{ij} m
\\
&&\qquad = (\mathbh{1}- \mathbb{E}_{1 \cdots r}) \sum_{a,b}^{(1 \cdots r)}
h_{ia} G^{(1
\cdots r)}_{ab} h_{bj} + O \bigl(
\varphi^C N^{-1} (d - 1)^{-2} \bigr)
\end{eqnarray*}
with high probability. We may therefore write
\[
B_{ij} = - \theta- m - (-h_{ij} + W_{ij} +
R_{ij} ),
\]
where
\[
W_{ij}:= (\mathbh{1}- \mathbb{E}_{1 \cdots r}) \sum
_{a,b}^{(1
\cdots r)} h_{ia} G^{(1 \cdots r)}_{ab}
h_{bj}\quad\mbox{and}\quad R_{ij} = O \bigl(
\varphi^C N^{-1} \bigl(\llvert d \rrvert - 1
\bigr)^{-2} \bigr)
\]
with high probability.

Next, we claim that
%
%e5.9 #&#
\begin{equation}
\label{estimateonW} W_{ij} = O \bigl(\varphi^C
N^{-1/2} \bigl(\llvert d \rrvert - 1\bigr)^{-1/2} \bigr)
\end{equation}
with high probability. Indeed, using Lemma~\ref{lemmaLDE} we get
\begin{eqnarray*}
\llvert W_{ij} \rrvert &\leq& \varphi^C \Biggl(
\frac{1}{N^2} \sum_{a,b}^{(1 \cdots r)} \bigl
\llvert G_{ab}^{(1 \cdots r)} \bigr\rrvert ^2
\Biggr)^{1/2}
\\
&=& \varphi^C \biggl(\frac{1}{N^2} \operatorname{Tr}
\bigl(G^{(1
\cdots r) *} G^{(1 \cdots r)} \bigr) \biggr)^{1/2}
\leq \varphi ^C N^{-1/2} \bigl(\llvert d \rrvert - 1 \bigr)^{-1/2}
\end{eqnarray*}
with high probability. In the last step we used (\ref{Gkm}), (\ref
{derivativeofG}), and $G = G^*$ to get (dropping the upper indices to
simplify notation)
\begin{eqnarray*}
\frac{1}{N^2} \operatorname{Tr}\bigl(G^* G\bigr) &=& N^{-1} m' + O
\bigl(\varphi^C N^{-2} \kappa_\theta^{-2}
\bigr)
\\
&=& O \bigl(N^{-1} \kappa_\theta^{-1/2} +
\varphi^C N^{-2} \kappa_\theta^{-2}
\bigr)
= O\bigl(N^{-1} \bigl(\llvert d \rrvert - 1\bigr)^{-1}\bigr)
\end{eqnarray*}
with high probability.

Using the bounds (\ref{estimateonW}) and $\llvert  h_{ij} \rrvert
\leq\varphi^C N^{-1/2}$ with high probability [as follows from
(\ref{subexpforh})], we may expand with (\ref{identityformsc}) to get
\[
M_{ij} = m \delta_{ij} + m^2
(-h_{ij} + W_{ij}) + O \bigl(\varphi^C
N^{-1} \bigl(\llvert d \rrvert - 1\bigr)^{-2} \bigr)
\]
with high probability. Let $H_{[1 \cdots r]} = H^{(r+1 \cdots N)}$
denote the upper
$r\times r$ block of $H$. Thus we get
%
%e5.10 #&#
\begin{eqnarray}
\label{schurforGaussiancase}
\qquad && N^{1/2} \bigl(\llvert d \rrvert - 1
\bigr)^{1/2} (M - m)
\nonumber\\[-8pt]\\[-8pt]
&&\qquad = m^2 N^{1/2} \bigl(\llvert d
\rrvert - 1\bigr)^{1/2} (-H_{[1 \cdots r]} + W) + O \bigl(
\varphi^C N^{-1/2} \bigl(\llvert d \rrvert - 1
\bigr)^{-3/2} \bigr)\nonumber
\end{eqnarray}
with high probability. In particular, for large enough $K$ we get
%
%e5.11 #&#
\begin{equation}
\label{HHW} N^{1/2} \bigl(\llvert d \rrvert - 1\bigr)^{1/2}
(M - m) \stackrel{d} {\sim} m^2 N^{1/2} \bigl(\llvert d
\rrvert - 1\bigr)^{1/2} (-H_{[1 \cdots r]} + W).
\end{equation}

By definition, $H_{[1 \cdots r]}$ and $W$ are independent. What
therefore remains is to compute the asymptotic distribution of $W$. We
claim that $W$ converges in law to an $r \times r$ Gaussian matrix:
%
%e5.12 #&#
\begin{equation}
\label{convergenceofW} N^{1/2} \bigl(\llvert d \rrvert - 1
\bigr)^{1/2} W \stackrel{d} {\sim} \frac
{1}{\sqrt{\llvert  d \rrvert  +
1}} \Phi.
\end{equation}
By the Cram\'er--Wold device, it suffices to show that
\[
N^{1/2} \bigl(\llvert d \rrvert - 1\bigr)^{1/2} \sum
_{i,j} Q_{ij} W_{ij} \stackrel{d} {\sim}
\frac
{1}{\sqrt{\llvert  d \rrvert  + 1}} \sum_{i,j} Q_{ij}
\Phi_{ij}
\]
for any deterministic matrix $Q = (Q_{ij})$ satisfying $Q = Q^*$ and
$Q_{ij} \in\mathbb{R}$ if $\beta= 1$. To that end, we diagonalize
$G^{(1\cdots r)}$ by writing
\[
N^{-1/2} \bigl(\llvert d \rrvert - 1\bigr)^{1/2}
G^{(1\cdots r)} = U^* \Lambda U,
\]
where $U$ is a unitary $(N - r) \times(N - r)$ matrix and $\Lambda=
\operatorname{diag}(\Lambda_{r + 1}, \ldots, \Lambda_N)$. Moreover, we
introduce the $r \times(N - r)$ matrix $h:=(h_{ia}\dvtx i \leq r, a
\geq r + 1)$. Since the entries of $h$ are i.i.d. Gaussians, $U$ is\vspace*{-1pt}
orthogonal/unitary, and $H$ is independent of $(\Lambda, U)$, we find
that $ (\Lambda, U h ) \stackrel{d}{=} (\Lambda, h )$. We conclude that
\begin{eqnarray*}
N^{1/2} \bigl(\llvert d \rrvert - 1\bigr)^{1/2} \sum
_{i,j = 1}^r Q_{ij} W_{ij} &=& N (1
- \mathbb{E}_{1 \cdots r}) \operatorname{Tr} \bigl(Q h^* U^* \Lambda U h \bigr)
\\
&\stackrel{d} {=}& N (1 - \mathbb{E}_{1 \cdots r}) \operatorname{Tr} \bigl(Q h^* \Lambda h
\bigr)
\\
&=& \sum_{a}^{(1 \cdots r)}
\Lambda_a \sum_{i,j = 1}^r
Q_{ij} N (h_{ia} h_{aj} - \mathbb{E}h_{ia}
h_{aj} )
\\
&=:& X.
\end{eqnarray*}
%
%Let $(\Lambda_a)_{a = r + 1}^N$ denote the eigenvalues of $N^{-1/2} (
%the vector $(h_{ai})_{a = r+1}^N$ for all $i \in\ququ{1,r}$ and the
%fact that $H^{(1 \cdots r)}$ is independent of the family $(h_{ia}
%N^{1/2} (\abs{d} - 1)^{1/2} \sum_{i,j = 1}^r Q_{ij} W_{ij} \eqdist
%h_{aj} - \E h_{ia} h_{aj}} \eqd X.
Note that $ (\sum_{i,j} Q_{ij} N  (h_{ia} h_{aj} - \mathbb {E}h_{ia}
h_{aj} ) )_{a = r+1}^N$ is a family of i.i.d. random variables,
independent of $\Lambda$, with variance $2 \beta^{-1} \operatorname{Tr}Q^2$.
Therefore,
\begin{eqnarray*}
\mathbb{E}X^2 &=& \frac{2}{\beta} \operatorname{Tr}Q^2 \sum
_{a}^{(1 \cdots r)} \Lambda_a^2
\\
&=& \frac{2}{\beta} \operatorname{Tr}Q^2 N^{-1} \bigl(\llvert d \rrvert - 1
\bigr) \operatorname{Tr} \bigl(G^{(1 \cdots r)}\bigr)^2
\\
&=& \frac{2}{\beta} \operatorname{Tr}Q^2 \bigl(\bigl(\llvert d \rrvert - 1\bigr)
m' + O\bigl(\varphi^C N^{-1} \bigl(\llvert d
\rrvert - 1\bigr)^{-3}\bigr) \bigr)
\\
&=& \frac{2}{\beta}\operatorname{Tr}Q^2 \bigl(\bigl(\llvert d \rrvert - 1\bigr) m' + O
\bigl(\varphi ^{-1}\bigr) \bigr)
\end{eqnarray*}
with high probability for large enough $K$, where we used (\ref{Gkm}).
Moreover, we have
\begin{eqnarray*}
\sum_a^{(1 \cdots r)} \Lambda_a^4
&=& N^{-2} \bigl(\llvert d \rrvert - 1\bigr)^2 \operatorname{Tr}
\bigl(G^{(1 \cdots r)}\bigr)^4
\\
&=& N^{-2} \bigl(\llvert d
\rrvert - 1\bigr)^2 \bigl(N m'''
/ 6 + O\bigl(\varphi^C \bigl(\llvert d \rrvert - 1
\bigr)^{-8}\bigr) \bigr)
\\
&=& O \bigl(N^{-1} \bigl(\llvert d \rrvert - 1\bigr)^{-3} +
N^{-2} \bigl(\llvert d \rrvert - 1\bigr)^{-6} \bigr)
= O
\bigl(\varphi^{-1}\bigr)
\end{eqnarray*}
with high probability for large enough $K$, where in the second step we
used (\ref{Gkm}) and in the third step the estimate $m''' \asymp
\kappa_\theta^{-5/2}$ as follows by differentiating (\ref{mprime})
twice and from Lemma~\ref{lemmamsc}.

We conclude from the central limit theorem that
\[
X \stackrel{d} {\sim} \mathcal N \biggl(0, \frac{2}{\beta
(\llvert  d \rrvert  + 1)}
\operatorname{Tr}Q^2 \biggr),
\]
where we used the identity
\[
\bigl(\llvert d \rrvert - 1\bigr) m' = \frac{1}{\llvert  d \rrvert  + 1}
\]
as follows from (\ref{mprime}) and (\ref{mtheta}). Thus,
(\ref{convergenceofW}) follows the identity
\[
\frac{1}{\sqrt{\llvert  d \rrvert  + 1}} \sum_{i,j} Q_{ij}
\Phi_{ij} \stackrel{d} {=} \mathcal N \biggl(0, \frac{2}{\beta(\llvert  d \rrvert  + 1)} \operatorname{Tr}
Q^2 \biggr)
\]
as follows from a from a simple variance calculation.

Next, by definition of $H_{[1 \cdots r]}$ we have
\[
-N^{1/2} \bigl(\llvert d \rrvert - 1\bigr)^{1/2}
H_{[1 \cdots r]} \stackrel {d} {=} \bigl(\llvert d \rrvert - 1
\bigr)^{1/2} \Phi.
\]
Thus, we find
\[
N^{1/2} \bigl(\llvert d \rrvert - 1\bigr)^{1/2}
(-H_{[1 \cdots r]} + W) \stackrel{d} {\sim} \frac
{\llvert  d \rrvert }{\sqrt{\llvert  d \rrvert  + 1}} \Phi.
\]
The claim now follows from (\ref{HHW}) and (\ref{mtheta}).
\end{pf}

%s6 #&#
\section{The almost Gaussian case}\label{secAGcase}

The next step of the proof is to consider the case where most entries
of $H$ are Gaussian. The exponent $\rho\geq2$ is used to define a
cutoff scale in the entries of $V$, below which the corresponding
entries of $H$ are assumed to be Gaussian.
Proposition~\ref{propalmostgaussian} will ultimately be fed into
Lemma~\ref{lemmaGreenfunctioncompwithshift} below, at which time we
shall choose $\rho$ to be large enough.

%
%pr6.1 #&#
\begin{proposition} \label{propalmostgaussian}
The following holds for large enough $K$. Let $\theta\equiv\theta(d)$
for some $d$ satisfying $\llvert  d \rrvert  \geq1 + \varphi^{K}
N^{-1/3}$. Let
$\rho\geq2$. Suppose that the Wigner matrix $H$ satisfies
%
%e6.1 #&#
\begin{equation}
\label{condforpropsition} \max_{1 \leq l \leq r} \max\bigl\{\llvert
V_{il} \rrvert, \llvert V_{jl} \rrvert \bigr\} \leq
\varphi^{-\rho} \quad\Longrightarrow\quad h_{ij}\mbox{ is
Gaussian}.
\end{equation}
Then
\[
N^{1/2} \bigl(\llvert d \rrvert - 1\bigr)^{1/2} \bigl(M(
\theta) - m(\theta ) \bigr) \stackrel{d} {\sim} - N^{1/2} \bigl(\llvert
d \rrvert - 1\bigr)^{1/2} d^{-2} V_\delta^* H
V_\delta+ \Psi_0,
\]
where $\Psi_0 = \Psi_0^*$ is a Gaussian matrix, independent of $H$,
with centred entries and covariance
\begin{eqnarray*}
\mathbb{E}(\Psi_0)_{ij} (\Psi_0)_{kl}
&=& \frac{\llvert  d \rrvert  - 1}{d^4} \bigl(\Delta_{ij,kl} - \mathcal P_{ij,kl}
\bigl(V_\delta^* V_\delta\bigr) \bigr) + \frac
{1}{d^4 (\llvert  d \rrvert  + 1)}
\Delta_{ij,kl}
\\
&&{} + \frac{\llvert  d\rrvert  - 1}{d^5} \mathcal Q_{ij,kl}(V) +
\frac{\llvert  d \rrvert  - 1}{d^6} \mathcal R_{ij,kl}(V).
\end{eqnarray*}
\end{proposition}

\begin{pf}
Throughout the proof, we drop the spectral parameter $z = \theta$ from
our notation.
\begin{longlist}[\textit{Step}~2.]
\item[\textit{Step} 1.] We start with some linear algebra in order to
    write the matrix $M$ in a form amenable to analysis. Since $\llVert
    \mathbf{v}^{(l)} \rrVert  = 1$ for all $l$ we find that
\[
\bigl\llvert \bigl\{i\dvtx \llvert V_{il} \rrvert >
\varphi^{-\rho} \bigr\} \bigr\rrvert \leq \varphi^{2\rho}.
\]
We shall permute the rows of $V$ by using an $N \times N$ permutation
matrix $O$ according to $M = V^* G V = (O V)^* OGO^* OV$. It is
easy to see that we may permute the rows of $V$ by setting $V \mapsto O
V$ so that after the permutation we have
\[
V = \pmatrix{ U \cr W },
\]
where:
\begin{longlist}[(iii)]
\item[(i)] $U$ is a $\mu\times r$ matrix and $W$ an $(N - \mu) \times
    r$ matrix,

\item[(ii)] $\llvert  W_{il} \rrvert  \leq\varphi^{-\rho}$ for all $i$
    and $l$,

\item[(iii)] $\mu\leq r \varphi^{2\rho}$.
\end{longlist}
After the permutation $H \mapsto O H O^*$, we may write $H$ as
\[
H = \pmatrix{ A & B^* \cr B & H_0 },
\]
where $A$ is a $\mu\times\mu$ matrix, $B$ an $(N - \mu) \times\mu$
matrix, and $H_0$ an $(N - \mu) \times(N - \mu)$ matrix with Gaussian
entries [as follows from (\ref{condforpropsition})].

Next, we rotate the rows of $W$ by choosing a unitary $(N - \mu) \times
(N - \mu)$ matrix $\widetilde S$ such that
\[
\widetilde S W = \pmatrix{ \widetilde W \cr
0 },
\]
where $\widetilde W$ is an $r \times r$ matrix that satisfies
%
%e6.2 #&#
\begin{equation}
\label{UUWW} U^* U + W^* W = U^* U + \widetilde W^* \widetilde W = \mathbh{1}_r.
\end{equation}
Thus, we get
\begin{eqnarray*}
M & =& V^* \pmatrix{ 1 & 0 \cr
0 & \widetilde S^* } \pmatrix{ 1 & 0
\cr
0 & \widetilde S } \pmatrix{ A - \theta& B^* \cr
B &
H_0 - \theta }^{-1} \pmatrix{ 1 & 0 \cr
0 &
\widetilde S^* } \pmatrix{ 1 & 0 \cr
0 & \widetilde S } V
\\
& \stackrel{d} {=}& \pmatrix{ U \vspace*{2pt}\cr
\widetilde W \cr 0 }^*
\pmatrix{ A - \theta& B^* \widetilde S^* \cr
\widetilde S B & H_0 -
\theta }^{-1} \pmatrix{ U \vspace*{2pt}\cr
\widetilde W \cr
0 },
\end{eqnarray*}
where $\stackrel{d}{=}$ denotes equality in distribution. Here we used the
unitary invariance of the Gaussian matrix $H_0$.

Next, we decompose
\[
H_0 = \pmatrix{ H_1 & Z^* \cr
Z &
H_2 }, \qquad \widetilde S = \pmatrix{ R \cr
S },
\]
where $H_1$ is an $r \times r$ Gaussian matrix, $Z$ an $(N - \mu- r)
\times r$ Gaussian matrix, and $H_2$ an $(N - \mu- r) \times(N - \mu
- r)$ Gaussian matrix. Moreover, $R$ is an $r \times(N - \mu)$ matrix
and we have
\[
R R^* = \mathbh{1}_r, \qquad S S^* = \mathbh{1}_{N - \mu-
r},
\qquad R S^* = 0, \qquad R^* R + S^* S = \mathbh{1}_{N - \mu
}.
\]
Thus, we find
\begin{eqnarray*}
M &\stackrel{d} {=}& \pmatrix{ U \vspace*{2pt}\cr
\widetilde W \cr
0 }^*
\pmatrix{ A - \theta& B^* R^* & B^* S^* \cr
RB & H_1 -
\theta& Z^* \cr
SB & Z & H_2 - \theta }^{-1}
\pmatrix{ U \vspace*{2pt}\cr
\widetilde W \cr
0 }
\\
&=& \pmatrix{ Y \cr
0 }^* \pmatrix{ \tilde A - \theta& F^* \cr
F & H_2 -
\theta }^{-1} \pmatrix{ Y \cr
0 } =: \Theta,
\end{eqnarray*}
where
\[
Y:= \pmatrix{ U \cr
\widetilde W }, \qquad F:= (SB, Z), \qquad\tilde A:=
\pmatrix{ A & B^* R^* \cr
RB & H_1 }.
\]
Here $Y$ is a $(\mu+ r) \times r$ matrix satisfying $Y^* Y = \mathbh{1}_r$,
and $F$ is an $(N - \mu- r) \times(\mu+r)$ matrix.

\item[\textit{Step} 2.] We claim that
%
%e6.3 #&#
\begin{equation}
\label{FFstar} F^* F = \mathbh{1}_{\mu+r} + O\bigl(
\varphi^C N^{-1/2}\bigr)
\end{equation}
with high probability (in the sense of matrix entries). In order to
prove (\ref{FFstar}), we write
\[
F^* F = \pmatrix{ B^* S^* S B & B^* S^* Z \cr
Z^* SB & Z^* Z }
\]
and consider each block separately. For $i \neq j$, we get using
(\ref{two-setLDE})
\begin{eqnarray*}
\bigl\llvert \bigl(B^* S^* S B \bigr)_{ij} \bigr\rrvert &=& \biggl
\llvert \sum_{k,l} \widebar B_{k i} \bigl(S^* S
\bigr)_{kl} B_{lj} \biggr\rrvert
\\
&\leq&\frac
{\varphi^C}{N}
\biggl(\sum_{k,l} \bigl\llvert \bigl(S^* S
\bigr)_{kl} \bigr\rrvert ^2 \biggr)^{1/2} =
\frac
{\varphi^C}{N} \bigl(\operatorname{Tr}\bigl(S^* S\bigr)^2 \bigr)^{1/2}
\leq\varphi ^C N^{-1/2}
\end{eqnarray*}
with high probability. Similarly, (\ref{diagLDE}) and
(\ref{offdiagLDE}) yield
\[
\bigl(B^* S^* S B \bigr)_{ii} = \sum_{k}
\bigl(S^* S\bigr)_{kk} \llvert B_{ki} \rrvert
^2 + \sum_{k \neq l} \widebar B_{k i}
\bigl(S^* S\bigr)_{kl} B_{li} = 1 + O\bigl(
\varphi^C N^{-1/2}\bigr)
\]
with high probability, where we used that $N^{-1} \operatorname{Tr}S^* S = 1 - (\mu+
r) N^{-1}$. Next, from (\ref{diagLDE}), (\ref{offdiagLDE}) and
(\ref{two-setLDE}) we easily get
%
%e6.4 #&#
\begin{equation}
\label{Z*Z} Z^* Z = \mathbh{1}_r + O\bigl(\varphi^C
N^{-1/2}\bigr)
\end{equation}
with high probability. Finally, (\ref{two-setLDE}) yields
\begin{eqnarray*}
\bigl\llvert \bigl(B^* S^* Z\bigr)_{ij} \bigr\rrvert &=& \biggl\llvert
\sum_{k,l} \widebar B_{ki} S^*_{kl}
Z_{lj} \biggr\rrvert\\
&\leq&\frac{\varphi^C}{N} \biggl(\sum
_{k,l} \bigl\llvert S^*_{kl} \bigr\rrvert
\biggr)^{1/2} = \frac{\varphi
^C}{N} \bigl(\operatorname{Tr}S^* S \bigr)^{1/2}
\leq\varphi^C N^{-1/2}
\end{eqnarray*}
with high probability. This concludes the proof of (\ref{FFstar}).

Next, we define
\[
G_2:= (H_2 - \theta)^{-1}
\]
and claim that
%
%e6.5 #&#
\begin{equation}
\label{FGF} F^* G_2 F = m + O \bigl(\varphi^C
N^{-1/2} \bigl(\llvert d \rrvert - 1\bigr)^{-1/2} \bigr)
\end{equation}
with high probability (in the sense of matrix entries). Since $N^{1/2}
(N - \mu- r)^{-1/2} H_2$ is an $(N - \mu- r) \times(N - \mu- r)$
GOE/GUE matrix that is independent of $F$, (\ref{FGF}) follows from
Theorem~\ref{theoremstrongestimate}, (\ref{kappad}),
(\ref{controlparam}) and (\ref{FFstar}).

\item[\textit{Step} 3.] For the following, we use the letter $\mathcal
    E$ to denote any (random) error term satisfying $\llvert  \mathcal E
    \rrvert \leq\varphi^C N^{-1} (\llvert  d \rrvert  - 1)^{-1}$ with
    high probability for some constant $C$. We apply Schur's complement
    formula to get
\begin{eqnarray*}
\Theta& =& Y^* \bigl(-\theta- m - \bigl(-\tilde A + F^* G_2 F - m
\bigr) \bigr)^{-1} Y
\\
& =& m Y^* Y - m^2 Y^* \tilde A Y + m^2 \bigl(Y^* F^*
G_2 F Y - m Y^* Y \bigr) + \mathcal E
\\
& =& m - m^2 Y^* \tilde A Y + m^2 \bigl(Y^* F^*
G_2 F Y - m \bigr) + \mathcal E,
\end{eqnarray*}
where in the second step we expanded using (\ref{identityformsc}) and
estimated the error term using (\ref{FGF}), $\mu\leq\varphi^C$, and
$\llVert \tilde A \rrVert \leq\varphi^C N^{-1/2}$ with high
probability. Using $R^* \widetilde W = W$, we get
\begin{eqnarray*}
\Theta &=& m - m^2 \bigl(U^* A U + U^* B^* W + W^* B U + \widetilde W^*
H_1 \widetilde W \bigr)
\\
&&{} + m^2 Y^* F^* (G_2 - m) F Y + m^3 \bigl(Y^*
F^* F Y - \mathbh{1}\bigr) + \mathcal E.
\end{eqnarray*}
Next, we rewrite the term $Y^* F^* (G_2 - m) F Y$ so as to decouple the
randomness of~$H_2$ from that of $F$. From (\ref{FFstar}), we find
%
%e6.6 #&#
\begin{equation}
\label{YFFY} Y^* F^* F Y = \mathbh{1}_r + O\bigl(
\varphi^C N^{-1/2}\bigr)
\end{equation}
with high probability. Define the deterministic $(N - \mu- r) \times
r$ matrix
\[
E_1:= \pmatrix{ \mathbh{1}_r \cr
0_{(N - \mu- 2r) \times r} }.
\]
Next, we claim that there is a unitary $(N - \mu- r) \times(N - \mu-
r)$ matrix $O_1$, which is $F$-measurable, such that
%
%e6.7 #&#
\begin{equation}
\label{OFY-E} \llVert O_1 FY - E_1 \rrVert \leq
\varphi^C N^{-1/2}
\end{equation}
with high probability. In order to prove (\ref{OFY-E}), write
$(\mathbf{x}_1, \ldots, \mathbf{x}_r):=FY$. Then (\ref{YFFY}) simply
states that the vectors $\mathbf{x}_1, \ldots, \mathbf{x}_r$ form a
basis of an $r$-dimensional subspace, which is orthonormal up to errors
of order $\varphi^C N^{-1/2}$ with high probability. More precisely, we
choose a unitary matrix $U_1$ such that $U_1 \mathbf{x}_1$ lies in the
direction of~$\mathbf{e}_1$. Hence, by (\ref{YFFY}), we have $U_1 FY =
(\mathbf{e}_1, U_1 \mathbf{x}_2, \ldots, U_1 \mathbf{x}_r) +
O(\varphi^C N^{-1/2})$ with high probability. Note moreover that by
(\ref{YFFY}) we have $\langle{\mathbf{e}_1}, {U_1 \mathbf{x}_i}\rangle=
O(\varphi^C N^{-1/2})$ with high probability for $i \geq2$. Next, we
choose a unitary matrix $U_2$ that leaves $\mathbf{e}_1$ invariant and
maximizes $\langle{\mathbf{e}_2}, {U_2 U_1 \mathbf{x}_2}\rangle$.
Hence, again by (\ref{YFFY}), we have $U_2 U_1 FY = (\mathbf{e}_1,
\mathbf{e}_2, U_2 U_1 \mathbf{e}_3,\break  \ldots, U_2 U_1 \mathbf{e}_r) +
O(\varphi^C N^{-1/2})$ with high probability. We continue in this
manner, at the $k$th step choosing a unitary matrix $U_k$ that leaves
$\mathbf{e}_1, \ldots, \mathbf{e}_{k - 1}$ invariant and maximizes
$\langle{\mathbf{e}_k}, {U_k \cdots U_1 \mathbf{x}_k}\rangle$. Finally,
we define $O_1:=U_r \cdots U_1$. By construction, the estimate in
(\ref{OFY-E}) holds. Moreover, since $Y$ is deterministic, $O_1$ is
clearly $F$-measurable. This concludes the proof of (\ref{OFY-E}).

Using Theorem~\ref{theoremstrongestimate} and the fact that $F$ and
$H_2$ are independent, we therefore get from (\ref{OFY-E})
\[
(O_1 FY )^* (G_2 - m) O_1 F Y =
E_1^* (G_2 - m) E_1 + \mathcal E.
\]
We conclude that
\begin{eqnarray*}
M &\stackrel{d} {=}& m - m^2 \bigl(U^* A U + U^* B^* W + W^* B U +
\widetilde W^* H_1 \widetilde W \bigr)
\\
&&{}+ m^2 E_1^*
(G_2 - m) E_1 + m^3 \bigl(Y^* F^* F Y -
\mathbh{1}\bigr) + \mathcal E,
\end{eqnarray*}
where we used that $O_1 G_2 O_1^* \stackrel{d}{=}G_2$ and that all
terms apart from $m^2 E_1^* (G_2 - m) E_1$ are independent of $H_2$.

Next, we compute
\begin{eqnarray*}
&& Y^* F^* F Y
\\
&&\qquad  = U^* B^* S^* S BU + U^* B^* S^* Z \widetilde W + \widetilde W^* Z^* SBU + \widetilde W^* Z^* Z \widetilde W
\\
&&\qquad  = U^* B^* B U - U^* B^* R^* R B U + U^* B^* S^* Z \widetilde W + \widetilde W^* Z^* SBU
+ \widetilde W^* Z^* Z \widetilde W
\\
&&\qquad  = U^* B^* B U + U^* B^* S^* Z \widetilde W + \widetilde W^* Z^* SBU + \widetilde W^* Z^* Z \widetilde W +
O\bigl(\varphi^C N^{-1}\bigr)
\end{eqnarray*}
with high probability, where in the last step we used
Lemma~\ref{lemmaLDE} and\break  \mbox{$\operatorname{Tr} (R^* R)^2 = r$.} Using (\ref{UUWW}), we
rewrite
\begin{eqnarray*}
&& U^* B^* B U + \widetilde W{}^* Z^* Z \widetilde W - \mathbh{1}
\\
&&\qquad= \IE \bigl(U^* B^* B U + \widetilde W^* Z^* Z \widetilde W \bigr) - \frac{\mu}{N} U^* U - \frac{\mu+
r}{N} \widetilde W^*
\widetilde W,
\end{eqnarray*}
where we introduced the notation $\IE X:=X - \mathbb{E}X$.

Thus, we conclude that
%
%e6.8 #&#
\begin{equation}
\label{M-msum} M - m \stackrel{d} {=} \Theta_1 +
\Theta_2 + \Theta_3 + \Theta _4 + \mathcal E,
\end{equation}
where
\begin{eqnarray*}
\Theta_1 &:=& m^2 E_1^* (G_2
- m) E_1,
\\
\Theta_2 &:=& -m^2 U^* A U,
\\
\Theta_3 &:=& -m^2 \widetilde W^* H_1 \widetilde W,
\\
\Theta_4 &:=& -m^2 \bigl(U^* B^* W + W^* B U \bigr)
\\
&&{}  + m^3 \IE \bigl(U^* B^* B U + U^* B^* S^* Z \widetilde W + \widetilde W^* Z^* SBU +
\widetilde W^* Z^* Z \widetilde W \bigr).
\end{eqnarray*}
By definition, the random variables $\Theta_1$, $\Theta_2$, $\Theta_3$
and $\Theta_4$ are independent.

\item[\textit{Step} 4.] We compute the asymptotics of $\Theta_1$,
    $\Theta_2$, and $\Theta_3$. We begin with $\Theta_1$. We shall
    apply Proposition~\ref{propositionGaussiancase} to the $(N - \mu-
    r) \times(N - \mu - r)$ Gaussian matrix~$H_2$. Thus, in
    Proposition~\ref{propositionGaussiancase} we replace $N$ with $N -
    \mu- r$, $H$ with $H_2$, and $M(\theta) = V^* (H - \theta)^{-1} V$
    by $V^* (H_2 - \theta)^{-1} V$ with $V:=E_1$. Since $\mu+ r
    \leq\varphi^C$ we find that $N - \mu - r \asymp N$. We therefore
    conclude from Proposition~\ref{propositionGaussiancase} that
\[
N^{1/2} \bigl(\llvert d \rrvert - 1\bigr)^{1/2}
\Theta_1 \stackrel{d} {\sim } \frac{1}{\llvert  d \rrvert ^3
\sqrt{\llvert  d \rrvert  + 1}} \Phi.
\]
Here we used (\ref{mtheta}). Recall that $\Phi$ is the rescaled GOE/GUE
matrix satisfying (\ref{GUEcovariances}).

In order to deal with $\Theta_2$, we introduce, in analogy to
$V_\delta$, the matrix $U_\delta= (U^\delta_{il})$ whose entries are
defined by $U_{il}^\delta:=U_{il} \mathbf{1}(\llvert  U_{il} \rrvert
> \delta)$. In particular, since $\delta\geq\varphi^{-1} \geq\varphi^{-\rho}$, we have
$V_\delta= {U_\delta\choose0}$. Writing $\widehat U_\delta= (\widehat
U_{il}^\delta):=U - U_\delta$, we get
\[
U^* A U = U_\delta^* A U_\delta+ \widehat U_\delta^* A
U_\delta+ U_\delta^* A \widehat U_\delta+ \widehat
U_\delta^* A \widehat U_\delta.
\]
Next, we define the matrices
\[
\Psi_1:= \widehat U_\delta^* A U_\delta+
U_\delta^* A \widehat U_\delta, \qquad\Psi_2:=\widehat
U_\delta^* A \widehat U_\delta.
\]
Note that, by definition, $\Psi_1$, $\Psi_2$ and $U_\delta^* A
U_\delta$ are independent. We now compute the covariances of the
matrices $\Psi_1$ and $\Psi_2$. A simple calculation yields
%
%e6.9 #&#
\begin{eqnarray}\label{covariancePsi}
N \mathbb{E}(\Psi_1)_{ij}
(\Psi_1)_{kl} &=& 2 \mathcal T_{ij,kl}
\bigl(U_\delta ^* U_\delta, \widehat U^*_\delta\widehat
U_\delta\bigr),
\nonumber\\[-8pt]\\[-8pt]
N \mathbb {E}(\Psi _2)_{ij}
(\Psi_2)_{kl} &=& \mathcal T_{ij,kl}\bigl(\widehat
U_\delta^* \widehat U_\delta,\widehat U_\delta^* \widehat
U_\delta\bigr),\nonumber
\end{eqnarray}
where we defined
\[
\mathcal T_{ij,kl}(R,T):= \tfrac{1}{2} \bigl(R_{il}
T_{kj} + R_{kj} T_{il} + \mathbf{1}(\beta= 1)
(R_{ik} T_{jl} + R_{jl} T_{ik} )
\bigr).
\]
For example, let us prove the second identity for the case $\beta= 2$.
Using $N \mathbb{E}h_{ab} h_{cd} = \delta_{ad} \delta_{bc}$ we find
\begin{eqnarray*}
N \mathbb{E}(\Psi_2)_{ij} (\Psi_2)_{kl}
&=& N \mathbb{E}\sum_{a,b,c,d = 1}^\mu \widehat
U^{\delta*}_{ia} h_{ab} \widehat U^\delta_{bj}
\widehat U^{\delta
*}_{kc} h_{cd} \widehat U^\delta_{dl}
\\
&=& \sum_{a,b = 1}^\mu\widehat U^{\delta*}_{ia}
\widehat U^\delta_{bj} \widehat U^{\delta*}_{kb} \widehat
U^\delta_{al} = \bigl(\widehat U^*_\delta \widehat
U_\delta \bigr)_{il} \bigl(\widehat U^*_\delta\widehat
U_\delta  \bigr)_{kj}.
\end{eqnarray*}
The other cases are handled similarly. Moreover, since by definition we
have $\llvert  \widehat U_{il}^\delta\rrvert  \leq\delta\ll1$, the
central limit theorem implies that $N^{1/2} \Psi_1$ and $N^{1/2}
\Psi_2$ converge to a Gaussian random matrix. Hence, the asymptotics of
$\Psi_1$ and $\Psi_2$ are governed entirely by their covariances
(\ref{covariancePsi}).

Similarly, $\Theta_3$ is Gaussian with covariance
\[
N \mathbb{E}(\Theta_3)_{ij} (\Theta_3)_{kl}
= d^{-4} \mathcal T_{ij,kl}\bigl(W^*W,W^*W\bigr),
\]
where we used (\ref{mtheta}). Using $U_\delta^* A U_\delta= V_\delta^*
H V_\delta$, we therefore conclude that
%
%e6.10 #&#
\begin{eqnarray}\label{sum3}
&& N^{1/2} \bigl(\llvert d \rrvert - 1\bigr)^{1/2}
(\Theta_1 + \Theta_2 + \Theta _3 )
\\
&&\qquad \stackrel{d} {\sim} \frac{1}{\llvert  d \rrvert ^3
\sqrt{\llvert  d \rrvert  + 1}} \Phi- N^{1/2} \frac{(\llvert  d \rrvert  - 1)^{1/2}}{d^2}
V_\delta ^* H V_\delta+ \Psi_3,\nonumber
\end{eqnarray}
where $\Psi_3$ is Gaussian with covariance
%
%e6.11 #&#
\begin{eqnarray}\label{varPsi3}
&& \mathbb{E}(\Psi_3)_{ij} (\Psi_3)_{kl}\nonumber
\\
&&\qquad = \frac{\llvert  d \rrvert  - 1}{d^{4}} \bigl(2\mathcal
T_{ij,kl}\bigl(U_\delta^* U_\delta, \widehat
U^*_\delta\widehat U_\delta\bigr)+ \mathcal T_{ij,kl}\bigl(
\widehat U_\delta^* \widehat U_\delta,\widehat U_\delta^* \widehat
U_\delta\bigr)
\\
&&\hspace*{163pt}
{} + \mathcal T_{ij,kl}\bigl(W^*W,W^*W\bigr) \bigr).\nonumber
\end{eqnarray}

\item[\textit{Step} 5.] Next, we compute the asymptotics of
    $\Theta_4$. We shall prove that $N^{1/2} (\llvert  d \rrvert  - 1)^{1/2}
    \Theta_4$ is asymptotically Gaussian, and compute its covariance
    matrix.

Using Lemma~\ref{lemmaLDE}, we find
%
%e6.12 #&#
\begin{eqnarray}\label{BSSB}
\bigl(B^* S^* S B\bigr)_{ij} & =& \sum_{k \neq l}
\bigl(S^* S\bigr)_{kl} \widebar B_{ki} B_{lj}
+ \sum _{k} \bigl(S^* S\bigr)_{kk} \widebar B_{ki} B_{kj}\nonumber
\\
& =& \sum_{k \neq l} \bigl(S^* S\bigr)_{kl}
\widebar B_{ki} B_{lj} + \sum_{k}
\bigl(S^* S\bigr)_{kk} \biggl(\widebar B_{ki} B_{kj}
- \frac{\delta_{ij}}{N} \biggr)
\nonumber\\[-8pt]\\[-8pt]
&&{}+ \frac{N - \mu- r}{N} \delta_{ij}\nonumber
\\
& =& \delta_{ij} + O\bigl(\varphi^C N^{-1/2}\bigr)\nonumber
\end{eqnarray}
with high probability. Define the deterministic $(N - \mu- r) \times
\mu$ matrix
\[
E_2:= \pmatrix{ \mathbh{1}_\mu\cr
0_{(N - 2 \mu- r) \times\mu} }.
\]
Exactly as after (\ref{BSSB}) we find that (\ref{BSSB}) and Gaussian
elimination imply that there is a unitary $(N - \mu- r) \times(N - \mu-
r)$ matrix $O_2$, which is $B$-measurable, such that
\[
\llVert O_2 SB - E_2 \rrVert \leq\varphi^C
N^{-1/2}
\]
with high probability. Thus, we get
\begin{eqnarray*}
\bigl\llvert \bigl(\widetilde W^* Z^* (O_2 SB - E_2) U
\bigr)_{ij} \bigr\rrvert & =& \biggl\llvert \sum
_k \widetilde W^*_{ik} \sum
_l \bigl((O_2 SB - E_2) U
\bigr)_{lj} \widebar Z_{lk} \biggr\rrvert
\\
& \leq&\varphi^C N^{-1/2} \bigl(U^* (O_2 S B -
E_2)^* (O_2 S B - E_2) U
\bigr)_{ii}^{1/2}
\\
& \leq&\varphi^C N^{-1}
\end{eqnarray*}
with high probability. Using that $Z$ is independent of $B$ and $O_2$,
we therefore find
\begin{eqnarray*}
\Theta_4 &\stackrel{d} {=}& -m^2 \bigl(U^* B^* W + W^* B
U \bigr)
\\
&&{} + m^3 \IE \bigl(U^* B^* B U + U^* E_2^* Z \widetilde W + \widetilde W^* Z^* E_2U + \widetilde W^* Z^* Z \widetilde W \bigr) + \mathcal E.
\end{eqnarray*}
Defining the $(N - \mu- r) \times r$ matrix
\[
\widetilde U:= E_2 U = \pmatrix{ U \cr
0_{(N - \mu- 2r) \times r} },
\]
we therefore have
\[
\Theta_4 \stackrel{d} {=} \Theta_4' +
\Theta_4'' + \mathcal E,
\]
where
\begin{eqnarray*}
\Theta_4' &:=& -m^2 \bigl(U^* B^* W + W^*
BU \bigr) + m^3 \IE \bigl(U^* B^* B U\bigr),
\\
\Theta_4'' &:=& m^3 \bigl(
\widetilde U^* Z \widetilde W + \widetilde W^* Z^* \widetilde U + \IE\bigl(\widetilde W Z^* Z \widetilde W\bigr) \bigr).
\end{eqnarray*}
By definition, $\Theta_4'$ and $\Theta_4''$ are independent. Recalling
that $\llvert  W_{il} \rrvert  \leq\varphi^{-\rho}$, we find from the central
limit theorem that $N^{1/2} \Theta_4'$ and $N^{1/2} \Theta_4''$ are
each asymptotically Gaussian. Hence, it suffices to compute their
covariances. A straightforward computation yields
\begin{eqnarray*}
N \mathbb{E}\bigl(\Theta'_4\bigr)_{ij}
\bigl(\Theta'_4\bigr)_{kl} &=& 2
m^4 \mathcal T_{ij,kl}\bigl(U^* U,W^* W\bigr) -
m^5 \mathcal Q_{ij,kl}(U,W)
\\
&&{} + m^6 \bigl(\mathcal
T_{ij,kl}\bigl(U^* U, U^* U\bigr) + \mathcal R_{ij,kl}(U) \bigr),
\end{eqnarray*}
where we defined
\begin{eqnarray*}
\mathcal Q_{ij,kl}(U,W)&:=& N^{-1/2} \sum
_{a,b} \bigl( \widebar U_{ai} \widebar U_{ak}
U_{al} \mu^{(3)}_{ab} W_{bj} + \widebar
W_{ia} \mu^{(3)}_{ab} U_{bj} \widebar U_{bk} U_{bl}
\\
&&\hspace*{47pt}{}+ \widebar U_{ak} \widebar U_{ai} U_{aj}
\mu^{(3)}_{ab} W_{bl} + \widebar W_{ka}
\mu^{(3)}_{ab} U_{bl} \widebar U_{bi}
U_{bj} \bigr).
\end{eqnarray*}
[By a slight abuse of notation, we write $\mathcal R_{ij,kl}(U)$ by
identifying $U$ with the $N \times r$ vector ${U\choose0}$.]

We may similarly deal with $\Theta_4''$. Using $\widetilde U^*
\widetilde U = U^* U$ and $\widetilde W^* \widetilde W = W^* W$ we find
\[
N \mathbb{E}\bigl(\Theta_4''
\bigr)_{ij} \bigl(\Theta_4''
\bigr)_{kl} = 2 m^6 \mathcal T_{ij,kl}\bigl(U^*
U,W^* W\bigr) + m^6 \mathcal T_{ij,kl}\bigl(W^* W,W^* W\bigr).
\]
Combining $\Theta_4'$ and $\Theta_4''$, and recalling (\ref{mtheta}),
we find
%
%e6.13 #&#
\begin{eqnarray}\label{theta4cov}
N \mathbb{E}(\Theta_4)_{ij} (
\Theta_4)_{kl} &=& 2 d^{-4} \mathcal
T_{ij,kl}\bigl(U^* U,W^* W\bigr) + d^{-5} \mathcal
Q_{ij,kl}(U,W)
\nonumber\\[-8pt]\\[-8pt]
&&{} + d^{-6} \bigl(\Delta_{ij,kl} +\mathcal R_{ij,kl}(U) \bigr),\nonumber
\end{eqnarray}
where we used that
\begin{eqnarray*}
&& \mathcal T_{ij,kl}\bigl(U^* U,U^* U\bigr) + \mathcal T_{ij,kl}
\bigl(W^*W,W^*W\bigr) + 2 \mathcal T_{ij,kl}\bigl(U^* U,W^* W\bigr)
\\
&&\qquad =\mathcal T_{ij,kl}(\mathbh {1},\mathbh{1}) = \Delta_{ij,kl}
\end{eqnarray*}
as follows from $W^* W + U^* U = \mathbh{1}$.

\item[\textit{Step} 6.] We may now consider the sum $\Theta_1 +
    \Theta_2 + \Theta_3 + \Theta_4$. From (\ref{M-msum}), (\ref{sum3}), (\ref
    {varPsi3}), (\ref{theta4cov}), and the definition of $\mathcal E$,
    we get
\[
N^{1/2} \bigl(\llvert d \rrvert - 1\bigr)^{1/2} (M - m)
\stackrel{d} {\sim} - N^{1/2} \bigl(\llvert d \rrvert - 1
\bigr)^{1/2} d^{-2} V_\delta^* H V_\delta+
\Psi_4,
\]
where $\Psi_4 = \Psi_4^*$ is a Gaussian matrix, independent of $H$,
with covariance
\begin{eqnarray*}
&& \mathbb{E}(\Psi_4)_{ij} (\Psi_4)_{kl}
\\
&&\qquad= \frac{\llvert  d \rrvert  - 1}{d^4} \bigl(\Delta_{ij,kl} - \mathcal P_{ij,kl}
\bigl(V_\delta^* V_\delta\bigr) \bigr) + \frac{\llvert  d
\rrvert  - 1}{d^5}
\mathcal Q_{ij,kl}(U,W)
\\
&&\quad\qquad{} + \frac{\llvert  d \rrvert  - 1}{d^6} \bigl(\Delta
_{ij,kl} + \mathcal R_{ij,kl}(U) \bigr) + \frac{1}{d^6 (\llvert  d
\rrvert  + 1)}
\Delta_{ij,kl}.
\end{eqnarray*}
Here we used that
\begin{eqnarray*}
&& 2 \mathcal T_{ij,kl}\bigl(U_\delta^* U_\delta, \widehat
U_\delta^* \widehat U_\delta\bigr) + \mathcal T_{ij,kl}\bigl(
\widehat U_\delta^* \widehat U_\delta, \widehat U_\delta^* \widehat
U_\delta\bigr)
\\
&&\quad{} + \mathcal T_{ij,kl}\bigl(W^* W, W^* W\bigr) + 2
\mathcal T_{ij,kl}\bigl(U^* U, W^* W\bigr)
\\
&&\qquad = \Delta_{ij,kl} - \mathcal T_{ij,kl}\bigl(U_\delta^*
U_\delta, U_\delta^* U_\delta\bigr) =
\Delta_{ij,kl} - \mathcal P_{ij,kl}\bigl(V_\delta^*
V_\delta\bigr)
\end{eqnarray*}
as follows from the bilinearity of $\mathcal T_{ij,kl}(\cdot, \cdot)$
as well as the identities $\mathcal T_{ij,kl}(\mathbh{1},\mathbh{1}) =
\Delta_{ij,kl}$, $\mathbh{1}= U_\delta^* U_\delta+ \widehat U_\delta^*
\widehat U_\delta+ W^* W$ and $U_\delta^* U_\delta= V_\delta^*
V_\delta$.

Using that $U$ is a $\mu\times r$ matrix with $\mu\leq r
\varphi^{2\rho}$ and $\llvert  W_{il} \rrvert  \leq\varphi^{-\rho
}$, we easily
find that
%
%e6.14 #&#
\begin{eqnarray}\label{errorincov}
\mathcal Q_{ij,kl}(U,W) &=& \mathcal Q_{ij,kl}(V) + O\bigl(\varphi^{-\rho}\bigr),
\nonumber\\[-8pt]\\[-8pt]
\mathcal R_{ij,kl}(U) &=& \mathcal R_{ij,kl}(V) + O\bigl(
\varphi^{-2\rho}\bigr).\nonumber
\end{eqnarray}
Since $\rho\geq2$, it is not hard to see that the errors on the
right-hand side of (\ref{errorincov}) are bounded from above (in the
sense of matrices) by the matrix $E_{ij,kl} = \varphi^{-1}
\Delta_{ij,kl}$. In particular, from (\ref{theta4cov}) we get that the
matrix
\[
2 d^{-4} \mathcal T_{ij,kl}\bigl(U^* U,W^* W\bigr) +
d^{-5} \mathcal Q_{ij,kl}(V) + d^{-6} \bigl(
\Delta_{ij,kl} + \mathcal R_{ij,kl}(V) \bigr) + E_{ij,kl}
\]
is nonnegative, from which we conclude that the right-hand side of
(\ref{defofcov}) is nonnegative. This completes the proof.\quad\qed
\end{longlist}\noqed
\end{pf}

%s7 #&#
\section{The general case} \label{secgencase}
The general case follows from Proposition~\ref{propalmostgaussian} and
Green function comparison. The argument is almost identical to that of
Section~7.4 in \cite{KY2}, and we only sketch the differences.

Let $H = (N^{-1/2} X_{ij})$ be an arbitrary real symmetric/complex
Hermitian Wigner matrix and $(N^{-1/2} Y_{ij})$ a GOE/GUE matrix
independent of $H$. For $\rho> 0$, define the subset
\[
I_\rho:= \bigl\{i \in[\![1,N ]\!]\dvtx \llvert V_{il}
\rrvert \leq\varphi^{-\rho}\mbox{ for all } l \in[\![1,r ]\!
] \bigr\}.
\]
Define a new Wigner matrix $\widehat H = (N^{-1/2} \widehat X_{ij})$
through
\[
\widehat X_{ij}:= \cases{ Y_{ij}, &\quad if $i \in
I_\rho$ and $j \in I_\rho$, \vspace*{2pt}
\cr
X_{ij}, &\quad otherwise.}
\]
Thus, $\widehat H$ satisfies the assumptions of
Proposition~\ref{propalmostgaussian}. Let
\[
J_\rho:= \bigl\{(i,j)\dvtx 1 \leq i \leq j \leq N, i \in
I_\rho \mbox{ and } j \in I_\rho \bigr\}.
\]
Choose a bijective map $\phi\dvtx J_\rho\to\{1, \ldots, \llvert
J_\rho\rrvert \}$. For $1 \leq\tau\leq\llvert  J_\rho \rrvert $ denote
by $H_\tau = (h_{ij}^\tau)$ the Hermitian matrix defined by
\[
h_{ij}^\tau:= \cases{ N^{-1/2} X_{ij}, &
\quad if $\phi(i,j)\leq\tau$ \vspace*{2pt}
\cr
N^{-1/2} \widehat X_{ij}, &\quad otherwise} \qquad(i \leq j).
\]
In particular, $H_0 = \widehat H$ and $H_{\llvert  J_\rho\rrvert } =
H$. Let now $(a,b) \in J_\rho$ satisfy $\phi(a,b) = \tau$. We write
\[
H_{\tau-1} = Q + N^{-1/2} \bigl(Y_{ab}
E^{(ab)} + \mathbf{1}(a \neq b) Y_{ba} E^{(ba)}
\bigr)
\]
and
\[
H_\tau= Q + N^{-1/2} \bigl(X_{ab}
E^{(ab)}+ \mathbf{1}(a \neq b) X_{ba} E^{(ba)} \bigr).
\]
Here $E^{(ab)}$ denotes the matrix with entries $E^{(ab)}_{ij}:=
\delta_{ai} \delta_{bj}$. Hence we have $Q_{ab} = Q_{ba} = 0$, and the
matrices $H_{\tau- 1}$ and $H_{\tau}$ differ only in the entries
$(a,b)$ and $(b,a)$.

Next, we introduce the resolvents
\[
R(z):= \frac{1}{Q-z}, \qquad S(z):= \frac
{1}{H_{\tau-1}-z}, \qquad T(z):=
\frac{1}{H_{\tau}-z}.
\]
Let $\llvert  d \rrvert  \geq1 + \varphi^K N^{-1/2}$. Set $z:=\theta
(d) + \ii N^{-4}$ (as in \cite{KY2}, Section~7.4, we add a~small
imaginary part to $z$ to ensure weak control on low-probability events)
and define
%
%e7.1 #&#
\begin{equation}
\label{definitionofx} x_R:= N^{1/2} \bigl(\llvert d
\rrvert - 1\bigr)^{1/2} \bigl(V^* R(z) V - m(z) \bigr).
\end{equation}
The quantities $x_S$ and $x_T$ are defined analogously with $R$
replaced by $S$ and $T$, respectively.

The following estimate is the main comparison estimate. It is very
similar to Lemma~7.13 of \cite{KY2}.

%
%le7.1 #&#
\begin{lemma} \label{lemmaGreenfunctioncompwithshift}
Provided $\rho$ is a large enough constant, the following holds. Let $f
\in C^3(\mathbb{C}^{r \times r})$ be bounded with bounded derivatives
and $q
\equiv q_N$ be an arbitrary deterministic sequence of $r \times r$
matrices. Then
%
%e7.2 #&#
%e7.3 #&#
\begin{eqnarray}
\mathbb{E}f(x_T + q) & =& \mathbb{E}f(x_R + q)+ \sum_{i,j = 1}^r
Z^{(ab)}_{ij} \mathbb{E}\frac{\partial f}{\partial x_{ij}}(x_R +q)
\nonumber\\[-8pt]\label{comparisonforT}\\[-8pt]
&&{}  + A_{ab} + O \bigl(\varphi^{-1} \mathcal E_{ab} \bigr),\nonumber
\\
\label{comparisonforS} \mathbb{E}f(x_S + q) & =&
\mathbb{E}f(x_R + q) + A_{ab} + O \bigl(
\varphi^{-1} \mathcal E_{ab} \bigr),
\end{eqnarray}
where $A_{ab}$ satisfies $\llvert  A_{ab} \rrvert  \leq\varphi^{-1}$,
\[
Z^{(ab)}_{ij}:= - N^{-1} \bigl(\llvert d \rrvert -
1\bigr)^{1/2} \bigl(m^4 \mu^{(3)}_{ab}
\widebar V_{ai} V_{bj} + m^4 \mu^{(3)}_{ba}
\widebar V_{bi} V_{aj} \bigr)
\]
and
%
%e7.4 #&#
\begin{equation}\label{defofcalE}
\mathcal E_{ab}:= \sum
_{i,j = 1}^r \sum_{\sigma, \tau= 0}^2
N^{-2 + \sigma/2 + \tau/2} \llvert V_{ai} \rrvert ^\sigma\llvert
V_{bj} \rrvert ^\tau+ \delta_{ab} \sum
_{i = 1}^r \sum_{\sigma= 0}^2
N^{-1 + \sigma/2} \llvert V_{ai} \rrvert ^\sigma.\hspace*{-30pt}
\end{equation}
\end{lemma}

\begin{pf}
The proof follows the proof of Lemma~7.13 of \cite{KY2} with cosmetic
modifications whose details we omit.
\end{pf}

Using Lemma~\ref{lemmaGreenfunctioncompwithshift}, we may now complete
the proof in the general case. The following proposition is the main
result of this section, and is the conclusion of the arguments from
Sections~\ref{secGaussian}--\ref{secgencase}.

%
%pr7.2 #&#
\begin{proposition} \label{propgeneraldistribution}
The following holds for large enough $K$. Let $\theta\equiv\theta(d)$
for some $d$ satisfying $\llvert  d \rrvert  \geq1 + \varphi^{K}
N^{-1/3}$. Then
\begin{eqnarray*}
&& N^{1/2} \bigl(\llvert d \rrvert - 1\bigr)^{1/2} \bigl(M(
\theta) - m(\theta ) \bigr)
\\
&&\qquad \stackrel{d} {\sim} - \frac{N^{1/2} (\llvert  d \rrvert  -
1)^{1/2}}{d^2}
V_\delta^* H V_\delta- \frac{(\llvert  d \rrvert  - 1)^{1/2} \mathcal S(V)}{d^4} + \Psi_0,
\end{eqnarray*}
where $\Psi_0$ is the Gaussian matrix from
Proposition~\ref{propalmostgaussian}.
\end{proposition}

\begin{pf}
The proof follows the proof of Theorem~2.14 in Section~7.4 of
\cite{KY2} with cosmetic modifications whose details we omit. The main
inputs are Proposition~\ref{propalmostgaussian} and
Lemma~\ref{lemmaGreenfunctioncompwithshift}. The imaginary part of the
spectral parameter $z = \theta(d) + \ii N^{-4}$ is easily removed using
the estimate $m(z) = -d + O(N^{-3})$. The condition $f \in C^3$ in
Lemma~\ref{lemmaGreenfunctioncompwithshift} can be relaxed to $f \in C$
by standard properties of weak convergence of measures.
\end{pf}

%s8 #&#
\section{Conclusion of the proof of Theorems \texorpdfstring{\protect\ref{thmeasystatement}}{2.5} and \texorpdfstring{\protect\ref{thmmainresult}}{2.9}} \label{secconcofproof}

We may now conclude the proof of Theorems \ref{thmeasystatement} and
\ref{thmmainresult}. First, we note that Theorem~\ref{thmeasystatement}
is an easy corollary of Theorem~\ref{thmmainresult}. We focus therefore
on the proof of Theorem~\ref{thmmainresult}.

Fix $K$ to be the constant from
Proposition~\ref{propgeneraldistribution}. Fix $\ell\in[\![1,r ]\!]$
and define the subset
\[
\Lambda:= \bigl\{N\in\mathbb{N}\dvtx \bigl\llvert d_\ell^{(N)}
\bigr\rrvert \geq1 + \varphi^K N^{-1/3} \bigr\}.
\]
We assume that $\Lambda$ is a subsequence (i.e., infinite), for
otherwise the claim of Theorem~\ref{thmmainresult} is vacuous. For
given $s > 0$, we introduce the partition
%
%e8.1 #&#
\begin{equation}
\label{Lambdapartition} \Lambda= \bigcup_{\gamma,\pi}
\Lambda_{\pi, \gamma}(s),
\end{equation}
where the union ranges over subsets $\pi,\gamma$ of $[\![1,r
]\!]$ satisfying $\ell\in\pi\subset\gamma\subset[\![1,r
]\!]$, and
\[
\Lambda_{\pi, \gamma}(s):= \bigl\{N \in\Lambda\dvtx \gamma _N(\ell) =
\gamma, \pi_N(\ell, s) = \pi \bigr\},
\]
where $\pi_N(\ell,s) \equiv\pi(\ell,s)$ and $\gamma_N(\ell) \equiv
\gamma(\ell)$ are the subsets from Definitions \ref{defpis} and
\ref{defgamma}.

We shall prove the following result.

%
%pr8.1 #&#
\begin{proposition} \label{propconclusion}
Fix $\ell$, $\pi$ and $\gamma$ satisfying $\ell\in\pi\subset \gamma
\subset[\![1,r ]\!]$. Let $\varepsilon> 0$ be given, and let $f_1,
\ldots, f_r$ be bounded continuous functions,\vadjust{\goodbreak} where $f_k$ is a function
on $\mathbb{R}^k$ satisfying $\llVert  f_k \rrVert _\infty\leq1$. Then
there exist constants $N_0$ and $s_0$, both depending on $\varepsilon$
and $f_1, \ldots, f_r$, such that (\ref{mainresult}) holds for all $s
\geq s_0$ and all $N \geq N_0$ satisfying $N \in\Lambda_{\pi,
\gamma}(s)$.
\end{proposition}

Before proving Proposition~\ref{propconclusion}, we note that it
immediately implies Theorem~\ref{thmmainresult}, since the partition
(\ref{Lambdapartition}) ranges over a finite family containing $O(1)$
elements.

\begin{pf*}{Proof of Proposition~\ref{propconclusion}}
From (\ref{roughlocationofq-group}), we know that $\theta(\mathcal B)$
contains with high probability precisely $\llvert  \gamma\rrvert $
outliers, namely $(\mu_{\alpha(i)})_{i \in\gamma}$. Following
(\ref{defzeta}), for $i \in\gamma$ we introduce the rescaled
eigenvalues
\[
\zeta_i = N^{1/2} \bigl(\llvert d_\ell\rrvert -
1\bigr)^{-1/2} (\mu _{\alpha(i)} - \theta_\ell).
\]
In order to identify the asymptotics of $\zeta_i$, we introduce the
$\llvert  \gamma\rrvert  \times\llvert  \gamma\rrvert $ matrices
\begin{eqnarray*}
X & \equiv& X_N:= -N^{1/2} \bigl(\llvert d_\ell
\rrvert - 1\bigr)^{1/2} \bigl(\llvert d_\ell\rrvert + 1\bigr)
\bigl(M_{[\gamma]}(\theta_\ell) - m(\theta_\ell )
\bigr),
\\
Y & \equiv& Y_N:= -N^{1/2} \bigl(\llvert d_\ell
\rrvert - 1\bigr)^{1/2} \bigl(\llvert d_\ell\rrvert + 1\bigr)
\bigl(D_{[\gamma]}^{-1} - d_\ell ^{-1} \bigr).
\end{eqnarray*}
Note that $X$ is random and $Y$ deterministic. From
(\ref{claimofreduction}), (\ref{mtheta}) and (\ref{mprime}), we get for
all $i \in \gamma$ that
%
%e8.2 #&#
\begin{equation}
\label{maininputforproof} \bigl\llvert \zeta_i -
\lambda_i(X + Y) \bigr\rrvert \leq\varphi^{-1}
\end{equation}
with high probability. By Proposition~\ref{propgeneraldistribution} and
Remark~\ref{remSQRbounded}, the family $(X_N)_N$ is tight.

By definition of $\pi$ and Lemma~\ref{lemmaddprime}, if $i \in \pi$ and
$j \in\gamma\setminus\pi$ then
%
%e8.3 #&#
\begin{equation}
\label{gapforpi} \llvert d_i - d_j \rrvert > s
N^{-1/2} \bigl(\llvert d_\ell\rrvert - 1\bigr)^{-1/2} /
2.
\end{equation}
We have the splitting
\[
D_{[\gamma]} = D_{[\pi]} \oplus D_{[\gamma\setminus\pi]}.
\]
We shall apply perturbation theory to the matrix $X+Y$. In order to do
so, we truncate $X$ by defining $X^t:=X \mathbf{1}(\llVert  X \rrVert
\leq t)$ for $t > 0$. Then by tightness of $X$ there exists a $t \equiv
t(\varepsilon)
> 0$ such that
%
%e8.4 #&#
\begin{equation}
\label{proboftruncation} \mathbb{P}\bigl(X_N \neq
X_N^t\bigr) \leq\frac{\varepsilon}{5}
\end{equation}
for all $N$. For the truncated matrices, we find the spectral gap
\begin{eqnarray*}
&& \operatorname{dist} \bigl(\sigma \bigl(X^t_{[\pi]} + Y_{[\pi]}
\bigr), \sigma \bigl(X^t_{[\gamma\setminus\pi]} + Y_{[\gamma\setminus
\pi]}\bigr) \bigr)
\\
&&\qquad \geq \operatorname{dist} \bigl(\sigma (Y_{[\pi]} ), \sigma
(Y_{[\gamma\setminus\pi]} ) \bigr) - 2 t \geq c s - 2t,
\end{eqnarray*}
where the constant $c$ only depends on $\Sigma$ in
(\ref{basicconditionsond}); here in the last step we used~(\ref{gapforpi}). Proposition~\ref{propperturbation} therefore yields
%
%e8.5 #&#
\begin{equation}
\label{perttheoryforpi} \bigl\llvert \lambda_i
\bigl(X^t + Y\bigr) - \lambda_i\bigl(X^t_{[\pi]}
+ Y_{[\pi]}\bigr) \bigr\rrvert \leq\frac{t^2}{cs - 2t - 2t^2}.
\end{equation}

We conclude that for there exists an $s_0$ and an $N_0$, both depending
on $\varepsilon$ and $f_{\llvert  \pi\rrvert }$, such that for $s
\geq s_0$ and $N
\geq N_0$ satisfying $N \in\Lambda_{\pi,\gamma}(s)$ we have
\begin{eqnarray*}
&& \bigl\llvert \mathbb{E}f_{\llvert  \pi\rrvert } \bigl((\zeta_i)_{i \in \pi}
\bigr) - \mathbb{E} f_{\llvert  \pi\rrvert } \bigl( \bigl(\lambda_i(X_{[\pi]}
+ Y_{[\pi]}) \bigr)_{i \in\pi} \bigr) \bigr\rrvert
\\
&&\qquad \leq \bigl \llvert \mathbb{E}f_{\llvert  \pi\rrvert } \bigl((\zeta _i)_{i \in\pi
}\bigr) - \mathbb{E}f_{\llvert  \pi\rrvert } \bigl( \bigl(\lambda_i
\bigl(X_{[\pi]}^t + Y_{[\pi]}\bigr)
\bigr)_{i \in\pi} \bigr) \bigr\rrvert + \frac{\varepsilon}{5}
\\
&&\qquad \leq\bigl\llvert \mathbb{E}f_{\llvert  \pi\rrvert } \bigl((\zeta
_i)_{i \in\pi
} \bigr) - \mathbb{E}f_{\llvert  \pi\rrvert } \bigl(
\bigl(\lambda_i\bigl(X^t + Y\bigr)
\bigr)_{i \in\pi} \bigr) \bigr\rrvert + \frac
{2\varepsilon}{5}
\\
&&\qquad \leq\bigl\llvert \mathbb{E}f_{\llvert  \pi\rrvert } \bigl((\zeta
_i)_{i \in\pi
} \bigr) - \mathbb{E}f_{\llvert  \pi\rrvert } \bigl(
\bigl(\lambda_i(X + Y) \bigr)_{i \in\pi} \bigr) \bigr\rrvert
+ \frac {3\varepsilon}{5}
\\
&&\qquad \leq\frac{4\varepsilon}{5},
\end{eqnarray*}
where in the first step we used (\ref{proboftruncation}), in the second
step (\ref{perttheoryforpi}) and dominated convergence, in the third
step (\ref{proboftruncation}) again, and in the last step
(\ref{maininputforproof}) and dominated convergence.
Proposition~\ref{propconclusion} now follows from
Proposition~\ref{propgeneraldistribution} applied to the $\llvert
\pi\rrvert \times\llvert  \pi \rrvert $ matrix
\[
N^{1/2} \bigl(\llvert d_\ell- 1 \rrvert \bigr)^{1/2}
\bigl(M(\theta_\ell) - m(\theta _\ell)
\bigr)_{[\pi]} = - \bigl(\llvert d_\ell\rrvert + 1
\bigr)^{-1} X_{[\pi]}.
\]\upqed
\end{pf*}

%s9 #&#
\section{The joint distribution: Proof of Theorem~\texorpdfstring{\protect\ref{thmmainresultJ}}{2.11}} \label{secjointdist}

In this final section, we extend the arguments of
Sections~\ref{secreduction2}--\ref{secconcofproof} to cover the joint
distribution of all outliers, and hence prove
Theorem~\ref{thmmainresultJ}.

We begin by introducing a coarser partition $\Gamma$, defined
analogously to $\Pi$ from Definition~\ref{defPi}, except that
$\pi(\ell,s)$ is replaced with $\gamma(\ell)$ from
Definition~\ref{defgamma}.

%
%de9.1 #&#
\begin{definition}\label{defjointcoarse}
Let $N$ and $D$ be given, and fix $K > 0$. We introduce a
partition\footnote{As in the footnote to Definition~\ref{defPi}, it is
easy to see that $\Gamma$ is a partition.} $\Gamma\equiv \Gamma(N,K,D)$
on a subset of $[\![1,r ]\!]$, defined as
\[
\Gamma:= \bigl\{\gamma(\ell)\dvtx \ell\in[\![1,r ]\!], \llvert
d_\ell\rrvert \geq1 + \varphi^K N^{-1/3} \bigr\}.
\]
We also use the notation $\Gamma= \{\gamma\}_{\gamma\in\Gamma}$.
\end{definition}

It is immediate from Definitions \ref{defPi} and \ref{defjointcoarse}
that $[\Pi] \subset\bigcup_{\gamma\in\Gamma} \gamma$ and that for each
$\pi\in\Pi$ there is a (unique) $\gamma\in\Gamma$ such that
$\pi\subset\gamma$. In analogy to (\ref{defdpi}), we set for
definiteness
\[
d_\gamma:= \min\{d_i\dvtx i\in\gamma\}, \qquad\theta
_\gamma:= \theta(d_\gamma).
\]
Note that for $\pi\in\gamma$ we have
%
%e9.1 #&#
\begin{equation}
\label{ddpiga} \frac{d_\pi}{d_\gamma} = 1+o(1), \qquad\frac{\llvert  d_\pi
\rrvert  -
1}{\llvert  d_\gamma\rrvert  - 1} = 1 + o(1).
\end{equation}

The following result follows from
Proposition~\ref{propositionreduction} and (\ref{mprime}).

%
%pr9.2 #&#
\begin{proposition} \label{propositionreductionJ}
The following holds for large enough $K$. For any $\gamma\in\Gamma$
and $i \in\gamma$ we have
%
%e9.2 #&#
\begin{equation}
\label{claimofreductionJ}
\quad\qquad \bigl\llvert \mu_{\alpha(i)} -
\lambda_i \bigl(\theta_{\gamma} - \bigl(d_{\gamma}^2-1
\bigr) \bigl(M(\theta_{\gamma} ) + D^{-1} \bigr)_{[\gamma]}
\bigr) \bigr\rrvert \leq\varphi^{-1} N^{-1/2} \bigl(\llvert
d_{\gamma} \rrvert - 1\bigr)^{1/2}
\end{equation}
with high probability.
\end{proposition}

As in Section~\ref{secconcofproof}, we may assume without loss of
generality that the partitions $\Pi$~and~$\Gamma$ are independent of
$N$. [Otherwise partition
\begin{eqnarray*}
\mathbb{N} &=& \bigcup_{\Gamma,\Pi} \Lambda_{\Pi,\Gamma}(s),
\\
\Lambda_{\Pi,\Gamma}(s) &:=& \bigl\{N \in\mathbb{N}\dvtx \Gamma (N,K,D) =
\Gamma, \Pi(N,K,s,D) = \Pi \bigr\}.
\end{eqnarray*}
Since the union is over a finite family of $O(1)$ subsets of $\mathbb
{N}$, we may first fix $\Gamma$~and~$\Pi$ and then restrict ourselves
to $N \in \Lambda_{\Gamma,\Pi}(s)$.] As in the proof of
Proposition~\ref{propconclusion}, we define for each $\pi\in\Gamma$ the
$\llvert \pi \rrvert \times \llvert  \pi\rrvert $ matrix
\[
X^{\pi}:= -N^{1/2} \bigl(\llvert d_{\pi} \rrvert - 1
\bigr)^{1/2} \bigl(\llvert d_{\pi} \rrvert + 1\bigr) \bigl(M(
\theta_{\pi}) - m(\theta_{\pi
}) \bigr)_{[\pi]}.
\]
The joint distribution of $(X^\pi)_{\pi\in\Pi}$ is described by the
following result, which is analogous to
Proposition~\ref{propgeneraldistribution}.

%
%pr9.3 #&#
\begin{proposition} \label{propgeneraldistributionJ}
For large enough $K$, we have
%
%e9.3 #&#
\begin{equation}
\label{XUP} \bigoplus_{\pi\in\Pi} X^\pi
\stackrel{d} {\sim} \bigoplus_{\pi\in\Pi} \bigl(
\Upsilon^\pi+\Psi^\pi \bigr),
\end{equation}
where $\Upsilon^\pi$ and $\Psi^\pi$ were defined in
Section~\ref{secjointdistrresults}.
\end{proposition}

We postpone the proof of Proposition~\ref{propgeneraldistributionJ} to
the next section, and finish the proof of Theorem~\ref{thmmainresultJ}
first. In order to identify the location of $\zeta_i^\pi$, we invoke
Proposition~\ref{propositionreductionJ} and make use of the freedom
provided by Lemma~\ref{lemreferencepoint} in order to change the
reference point $\theta_\gamma$ at will. Thus,
Proposition~\ref{propositionreductionJ} and
Lemma~\ref{lemreferencepoint} yield, for any $\pi\in\Pi$, $i \in\pi$,
and $\gamma\in\Gamma$ containing $\pi$, that
%
%e9.4 #&#
\begin{eqnarray}\label{xs1j}
\qquad\zeta^\pi_{i} &=& N^{1/2} \bigl(
\llvert d_{\pi} \rrvert - 1\bigr)^{-1/2} (\mu _{\alpha(i)} -
\theta_\pi )
\nonumber\\[-8pt]\\[-8pt]
&=& -N^{1/2} \bigl(\llvert d_{\pi}
\rrvert - 1\bigr)^{1/2} \bigl(\llvert d_{\pi} \rrvert + 1\bigr)
\lambda_i \bigl(\bigl(M(\theta _\pi) + D^{-1}
\bigr)_{[\gamma]} \bigr) + O\bigl(\varphi^{-1}\bigr)\nonumber
\end{eqnarray}
with high probability, where we used (\ref{mprime}), (\ref{ddpiga})
and Lemma~\ref{lemmagrowthofspectrum}.

Next, for $\pi\in\Pi$ let $\gamma(\pi)$ denote the unique element of
$\Gamma$ that contains $\pi$. For each $\pi\in\Pi$, we introduce the
$\llvert  \gamma(\pi) \rrvert  \times\llvert  \gamma(\pi) \rrvert $
matrices
\begin{eqnarray*}
\widetilde X{}^{\pi} &:=& -N^{1/2} \bigl(\llvert d_{\pi}
\rrvert - 1\bigr)^{1/2} \bigl(\llvert d_{\pi} \rrvert + 1\bigr)
\bigl(M(\theta_{\pi}) - m(\theta _{\pi}) \bigr)_{[\gamma(\pi)]},
\\
\widetilde Y{}^\pi&:=& -N^{1/2} \bigl(\llvert d_{\pi}
\rrvert - 1\bigr)^{1/2} \bigl(\llvert d_{\pi} \rrvert + 1\bigr)
\bigl(D^{-1} -d_{\pi}^{-1} \bigr)_{[\gamma
(\pi)]}.
\end{eqnarray*}
Thus, (\ref{xs1j}) reads
\[
\zeta^\pi_{i} = \lambda_i\bigl(\widetilde X{}^\pi+ \widetilde Y{}^\pi\bigr) + O\bigl(\varphi^{-1}
\bigr)
\]
with high probability. By Proposition~\ref{propgeneraldistribution} and
Remark~\ref{remSQRbounded}, $\widetilde X{}^\pi$ is tight (in $N$). We
may now repeat verbatim the truncation and perturbation theory argument
from the proof of Proposition~\ref{propconclusion}, following
(\ref{gapforpi}). The conclusion is that there exists an $s_0$ and an
$N_0$, both depending on $\varepsilon$ and $f_{\llvert  [\Pi] \rrvert
}$, such that for $s \geq s_0$ and $N \geq N_0$ we have
\[
\bigl\llvert \mathbb{E}f_{\llvert  [\Pi] \rrvert } \bigl(\bigl(\zeta_i^\pi
\bigr)_{\pi\in\Pi, i
\in\pi} \bigr) - \mathbb{E}f_{\llvert  [\Pi] \rrvert } \bigl( \bigl(
\lambda_i \bigl[{\bigl(\widetilde X{}^\pi+ \widetilde Y{}^\pi
\bigr)_{[\pi]}} \bigr] \bigr)_{\pi\in\Pi, i \in\pi
} \bigr) \bigr\rrvert \leq
\frac{\varepsilon}{2}.
\]
The claim now follows from Proposition~\ref{propgeneraldistributionJ}
and the observation that\break  $(\widetilde X{}^\pi)_{[\pi]} = X^\pi$. This
concludes the proof of Theorem~\ref{thmmainresultJ}.

%s9.1 #&#
\subsection{Proof of Proposition~\texorpdfstring{\protect\ref{propgeneraldistributionJ}}{9.3}} \label{secgendistrJ}
What remains is to prove Proposition~\ref{propgeneraldistributionJ}.
Clearly, it is a generalization of
Proposition~\ref{propgeneraldistribution}. The proof of
Proposition~\ref{propgeneraldistributionJ} relies on the same
three-step strategy as that of
Proposition~\ref{propgeneraldistribution}: the Gaussian case, the
almost Gaussian case and the general case.

We begin with the Gaussian case (generalization of
Section~\ref{secGaussian}).

%
%pr9.4 #&#
\begin{proposition} \label{propositionGaussiancase2}
Suppose that $H$ is a GOE/GUE matrix. Then for large enough $K$ we have
\[
\bigoplus_{\pi\in\Pi} N^{1/2} \bigl(\llvert
d_{\pi} \rrvert - 1\bigr)^{1/2} \bigl(M(\theta_{\pi})
- m(\theta_{\pi}) \bigr)_{[\pi]} \stackrel{d} {\sim} \bigoplus
_{\pi\in\Pi} \frac{ 1}{\llvert  d_\pi\rrvert  \sqrt{\llvert  d_\pi\rrvert  +
1}} \Phi_{\pi};
\]
here $(\Phi_\pi)_{\pi\in\Pi}$ is a family of independent Gaussian
matrices, where each $\Phi_\pi$ is a $\llvert  \pi\rrvert  \times
\llvert  \pi\rrvert $ matrix whose covariance is given by
(\ref{GUEcovariances}).
\end{proposition}

\begin{pf}
The proof is a straightforward extension of that of
Proposition~\ref{propositionGaussiancase}, and we only indicate the
changes. For each argument $\theta_\pi$, we use Schur's complement
formula on the whole block $[\![1,r ]\!]$. Thus, instead of
(\ref{schurforGaussiancase}), we get
\begin{eqnarray*}
\hspace*{-4pt}&& N^{1/2} \bigl(\llvert d_\pi\rrvert - 1\bigr)^{1/2}
\bigl(M(\theta_\pi) - m(\theta _\pi) \bigr)
\\
\hspace*{-4pt}&&\qquad = d_\pi^{-2} N^{1/2} \bigl(\llvert
d_\pi\rrvert - 1\bigr)^{1/2} \bigl(-H_{[1 \cdots
r]} + W(
\theta_\pi) \bigr)
+ O \bigl(\varphi^C N^{-1/2}
\bigl(\llvert d_\pi\rrvert - 1\bigr)^{-3/2} \bigr).
\end{eqnarray*}
This gives
%
%e9.5 #&#
\begin{eqnarray}\label{HHWJ}
&& \bigoplus_{\pi\in\Pi} N^{1/2}
\bigl(\llvert d_{\pi} \rrvert - 1\bigr)^{1/2} \bigl(M(
\theta_{\pi}) - m(\theta_{\pi}) \bigr)_{[\pi]}
\nonumber\\[-8pt]\\[-8pt]
&&\qquad\stackrel{d} {\sim} \bigoplus_{\pi\in\Pi} d
_{\pi}^{-2} N^{1/2} \bigl(\llvert d_\pi
\rrvert - 1\bigr)^{1/2} \bigl(-H_{[1 \cdots
r]} + W(\theta_\pi) \bigr)_{[\pi]},\nonumber
\end{eqnarray}
which is the appropriate generalization of (\ref{HHW}). By definition,
$H_{[1\cdots r]}$ is independent of the family of matrices
$(W(\theta_\pi))_{\pi\in\Pi}$, and the submatrices $H_{[\pi]}$, $\pi
\in\Pi$, are obviously independent. We may now repeat verbatim the
proof of (\ref{convergenceofW}) to get
%
%e9.6 #&#
\begin{equation}
\label{convergenceofWJ} \bigoplus_{\pi\in\Pi}
N^{1/2} \bigl(\llvert d_\pi\rrvert - 1\bigr)^{1/2}
W_{[\pi]} (\theta_\pi) \stackrel{d} {\sim} \bigoplus
_{\pi\in\Pi}\frac{1}{ \sqrt {\llvert  d_\pi\rrvert  +
1}} \Phi_\pi.
\end{equation}
The claim now follows from (\ref{HHWJ}).
\end{pf}

Next, we consider the almost Gaussian case (generalization of
Section~\ref{secAGcase}).

%
%pr9.5 #&#
\begin{proposition} \label{propalmostgaussian2}
Let $\rho> 0$. Suppose that the Wigner matrix $H$ satisfies
%
%e9.7 #&#
\begin{equation}
\label{condforpropsition2} \max_{1 \leq l \leq r} \max\bigl\{\llvert
V_{il} \rrvert, \llvert V_{jl} \rrvert \bigr\} \leq
\varphi^{-\rho} \quad\Longrightarrow\quad h_{ij}\mbox{ is
Gaussian}.
\end{equation}
Define $\widetilde \Upsilon$ to be the matrix $\Upsilon$ without the
shift arising from $\mathcal S(V)$, that is, $\widetilde \Upsilon=
\bigoplus_{\pi\in \Pi} \widetilde \Upsilon{}^\pi$ with
\[
\widetilde \Upsilon{}^\pi:= \bigl(\llvert d_{\pi} \rrvert + 1\bigr)
\bigl(\llvert d_{\pi} \rrvert - 1\bigr)^{1/2} \biggl(
\frac{N^{1/2} V_\delta^* H V_\delta}{d_{\pi
}^2} \biggr)_{[\pi]}.
\]
Then for large enough $K$ we have
%
%e9.8 #&#
\begin{equation}
\label{hjTw} \bigoplus_{\pi\in\Pi} X^\pi
\stackrel{d} {\sim} \bigoplus_{\pi\in\Pi} \bigl(\widetilde \Upsilon{}^\pi+\Psi^\pi \bigr).
\end{equation}
\end{proposition}

\begin{pf}
We start exactly as in the proof of
Proposition~\ref{propalmostgaussian}. We repeat the steps up to
(\ref{M-msum}) verbatim on the family of $r \times r$ matrices $
(M(\theta_\pi) - m(\theta _\pi) )_{\pi\in\Pi}$, whereby all of the
reduction operations are performed simultaneously on each matrix
$M(\theta_\pi) - m(\theta_\pi)$. Note that these matrices only differ
in the argument $\theta_\pi$; hence all steps of the reduction (and in
particular the quantities $O$, $O_1$, $U$, $W$, $\widetilde W$, $A$,
$B$, $H_0$, $H_1$, $Z$, etc.) are the same for all matrices
$M(\theta_\pi) - m(\theta_\pi)$. We take over the notation from the
proof of Proposition~\ref{propalmostgaussian} without further comment.
Thus, we are led to the following generalization of (\ref{M-msum}):
%
%e9.9 #&#
\begin{equation}
\label{M-msum2} \bigoplus_{\pi\in\Pi} X^\pi
\stackrel{d} {\sim} \Theta_1 + \Theta_2 +
\Theta_3 + \Theta_4' +
\Theta_4'',
\end{equation}
where
\begin{eqnarray*}
\Theta_1 &:=& \bigoplus_{\pi\in\Pi}
\bigl(-N^{1/2} \bigl(\llvert d_{\pi} \rrvert - 1
\bigr)^{1/2} \bigl(\llvert d_{\pi
} \rrvert + 1\bigr)
d_\pi^{-2}
\bigl[{E_1^* \bigl(G_2(\theta_\pi) - m(\theta_\pi ) \bigr) E_1}\bigr]_{[\pi] } \bigr),
\\
\Theta_2 &:=& \bigoplus_{\pi\in\Pi}
\bigl(N^{1/2} \bigl(\llvert d_{\pi} \rrvert - 1
\bigr)^{1/2} \bigl(\llvert d_{\pi
} \rrvert + 1\bigr)
 d_\pi^{-2} \bigl[{U^* A U} \bigr] _{[\pi] }  \bigr),
\\
\Theta_3 &:=& \bigoplus_{\pi\in\Pi}
\bigl(N^{1/2} \bigl(\llvert d_{\pi} \rrvert - 1
\bigr)^{1/2} \bigl(\llvert d_{\pi
} \rrvert + 1\bigr)
 d_\pi^{-2} \bigl[{\widetilde W^* H_1 \widetilde W} \bigr]
_{[\pi] } \bigr),
\\
\Theta_4' &:=& \bigoplus_{\pi\in\Pi}
\bigl(N^{1/2} \bigl(\llvert d_{\pi} \rrvert - 1
\bigr)^{1/2} \bigl(\llvert d_{\pi
} \rrvert + 1\bigr)
\\
&&\hspace*{20pt}{}\times \bigl[{d_\pi^{-2} \bigl(U^* B^* W + W^* B U \bigr) +
d_\pi^{-3} \IE \bigl(U^* B^* B U \bigr)}
\bigr]_{[\pi] } \bigr),
\\
\Theta_4'' &:=& \bigoplus
_{\pi\in\Pi} \bigl(N^{1/2} \bigl(\llvert d_{\pi}
\rrvert - 1\bigr)^{1/2} \bigl(\llvert d_{\pi
} \rrvert + 1\bigr)
 d_\pi^{-3}
\\
&&\hspace*{20pt}{}\times  \bigl[{\widetilde U^* Z \widetilde W + \widetilde W^* Z^* \widetilde U + \IE
\bigl(\widetilde W Z^* Z \widetilde W\bigr)} \bigr]_{[\pi] } \bigr).
\end{eqnarray*}
(We deviate somewhat from the convention of Section~\ref{secAGcase} in
that, unlike there, we include the normalization factor, which depends
on $\pi$, in the definition of the variables $\Theta$.) By definition,
the random matrices $\Theta_1$, $\Theta_2$, $\Theta_3$, $\Theta_4'$
and $\Theta_4''$ are independent. They are all block diagonal, and we
sometimes use the notation $\Theta_1 = \bigoplus_{\pi\in\Pi}
\Theta_1^{\pi}$, etc., for their blocks. What remains is to identify
their individual asymptotic distributions.

The matrix is $\Theta_1$ is easy: from
Proposition~\ref{propositionGaussiancase2} we immediately get
\[
\Theta_1 \stackrel{d} {\sim} \bigoplus
_{\pi\in\Pi} \frac{ \sqrt{\llvert  d_\pi\rrvert  + 1}}{\llvert  d_\pi\rrvert ^3} \Phi_{\pi},
\]
where $(\Phi_\pi)_{\pi\in\Pi}$ is defined as in
Proposition~\ref{propositionGaussiancase2}. The matrix $\Theta_2$ is
dealt with in the same way as in the proof of
Proposition~\ref{propalmostgaussian}; we omit the details. By
definition, $\Theta_3$~is Gaussian with mean zero. A short computation
yields the covariance
\[
\mathbb{E}\bigl(\Theta^\pi_3\bigr)_{ij}
\bigl(\Theta^{\pi'}_3\bigr)_{kl}
= \biggl(\prod
_{p=\pi, \pi'}\frac{(\llvert  d_p \rrvert  - 1)^{1/2} (\llvert  d_p \rrvert + 1)}{d_p^2} \biggr) \mathcal
T_{ij,kl}\bigl(W^* W,W^* W\bigr)
\]
for $\pi,\pi' \in\Pi$, $i,j \in\pi$ and $k,l \in\pi'$. We may therefore
conclude that, similar to (\ref{sum3}) and (\ref{varPsi3}), we have
%
%e9.10 #&#
\begin{equation}
\label{sum3J} (\Theta_1 + \Theta_2 +
\Theta_3 ) \stackrel{d} {\sim } \bigoplus
_{\pi\in\Pi} \frac{ \sqrt{\llvert  d_\pi\rrvert  + 1}}{\llvert  d_\pi\rrvert ^3} \Phi_{\pi} + \bigoplus
_{\pi\in\Pi} \widetilde \Upsilon{}^\pi + \bigoplus
_{\pi\in\Pi} \Psi^\pi_3,
\end{equation}
where $\bigoplus_{\pi\in\Pi} \Psi^\pi_3$ is a block diagonal Gaussian
matrix with mean zero and covariance
%
%e9.11 #&#
%e9.12 #&#
\begin{eqnarray}\label{varPsi3J}
\qquad\mathbb{E}\bigl(\Psi_3^\pi \bigr)_{ij} \bigl(\Psi_3^{\pi'}\bigr)_{kl}
&=& \biggl(\prod_{p=\pi, \pi'}\frac{(\llvert  d_{p} \rrvert  - 1)^{1/2}(\llvert  d_{p} \rrvert + 1)}{ d_p^2} \biggr)\nonumber
\\
&&{} \times \bigl(2\mathcal T_{ij,kl}\bigl(U_\delta^*
U_\delta, \widehat U^*_\delta\widehat U_\delta\bigr) + \mathcal
T_{ij,kl}\bigl(\widehat U_\delta^* \widehat U_\delta,\widehat U_\delta^* \widehat U_\delta\bigr)
\\
&&\hspace*{110pt}{} + \mathcal T_{ij,kl}\bigl(W^*W,W^*W\bigr) \bigr)\nonumber
\end{eqnarray}
for $\pi,\pi' \in\Pi$, $i,j \in\pi$ and $k,l \in\pi'$.

Next, we deal with $\Theta_4'$ and $\Theta_4''$. By the central limit
theorem and the definition of $W$, as in the proof of
Proposition~\ref{propalmostgaussian}, both of these matrices are
asymptotically Gaussian (with mean zero). The variances may be computed
along the same lines as in the proof of
Proposition~\ref{propalmostgaussian}. The result is, for $\pi$, $\pi'
\in\Pi$, $i$, $j \in\pi$ and \mbox{$k$, $l \in\pi'$,}
\begin{eqnarray*}
&& \mathbb{E}\bigl(\Theta'_4\bigr)_{ij} \bigl(
\Theta'_4\bigr)_{kl}
\\
&&\qquad  = \biggl(\prod
_{p=\pi, \pi'} \frac{(\llvert  d_{p} \rrvert  -
1)^{1/2} (\llvert  d_{p} \rrvert + 1)}{d_p^2} \biggr)
\\
&&\quad\qquad{}\times \biggl( 2 T_{ij,kl}\bigl(U^* U,W^* W\bigr) +
\frac{1}{d_{\pi} d_{\pi' }} \bigl(\mathcal T_{ij,kl}\bigl(U^* U,U^* U\bigr)+
\mathcal R_{ij,kl}(U) \bigr)
\\
&&\hspace*{51pt}{}+ \frac{N^{-1/2}}{d_{\pi'}}\sum_{a,b} \bigl(
\widebar U_{ai} \widebar U_{ak} U_{al}
\mu^{(3)}_{ab} W_{bj} + \widebar W_{ia} \mu
^{(3)}_{ab} U_{bj} \widebar U_{bk}
U_{bl} \bigr)
\\
&&\hspace*{78pt}{}+\frac{N^{-1/2}}{d_{\pi}}\sum_{a,b} \bigl(
\widebar U_{ak} \widebar U_{ai} U_{aj}
\mu^{(3)}_{ab} W_{bl} + \widebar W_{ka} \mu
^{(3)}_{ab} U_{bl} \widebar U_{bi}
U_{bj} \bigr) \biggr)
\end{eqnarray*}
as well as
\begin{eqnarray*}
\mathbb{E}\bigl(\Theta_4''
\bigr)_{ij} \bigl(\Theta_4''
\bigr)_{kl} &=& \biggl(\prod_{p=\pi, \pi'} \frac{(\llvert  d_{p} \rrvert  -
1)^{1/2} (\llvert  d_{p} \rrvert + 1)}{d_p^3} \biggr)
\\
&&{}\times\bigl(2 \mathcal T_{ij,kl}\bigl(U^* U,W^* W
\bigr) + \mathcal T_{ij,kl}\bigl(W^* W,W^* W\bigr) \bigr).
\end{eqnarray*}

Putting everything together, we get
%
%e9.13 #&#
\begin{equation}\label{hjT}
\bigoplus_{\pi\in\Pi} X^\pi
\stackrel{d} {\sim} \bigoplus_{\pi \in\Pi} \widetilde \Upsilon{}^\pi +
\bigoplus_{\pi\in\Pi} \Psi^\pi_4,
\end{equation}
where $\bigoplus_{\pi\in\Pi} \Psi^\pi_4$ is a Gaussian block diagonal
matrix with mean zero that is independent of $H$, and whose
covariance is given by
\begin{eqnarray*}
&& \mathbb{E}\bigl(\Psi^\pi_4\bigr)_{ij} \bigl(\Psi^{\pi'}_4\bigr)_{kl}
\\
&&\qquad  = \frac{\llvert  d_{\pi} \rrvert  + 1 }{d_{\pi}^{2} }
\delta_{\pi\pi'}\Delta _{ij, kl} + \delta_{\pi\pi'}E_{ij, kl}
\\
&&\quad\qquad{} + \biggl(\prod_{p=\pi, \pi'}
\frac{(\llvert  d_{p} \rrvert  -
1)^{1/2} (\llvert  d_{p} \rrvert + 1)}{d_p^2} \biggr)
\\
&&\hspace*{44pt}{}\times \biggl( - \mathcal P_{ij,kl}\bigl(V_\delta^*
V_\delta\bigr) + \frac{1}{d_{{\pi}}  d_{{\pi'}} } \mathcal R_{ij,kl}(V)
\\
&&\hspace*{62pt}{}+ \frac{N^{-1/2}}{ d_{\pi'} }\sum_{a,b} \bigl(
\widebar U_{ai} \widebar U_{ak} U_{al}
\mu^{(3)}_{ab} W_{bj} + \widebar W_{ia} \mu
^{(3)}_{ab} U_{bj} \widebar U_{bk}
U_{bl} \bigr)
\\
&&\hspace*{62pt}{} + \frac{N^{-1/2}}{ d_\pi} \sum_{a,b}
\bigl( \widebar U_{ak} \widebar U_{ai} U_{aj}
\mu^{(3)}_{ab} W_{bl} + \widebar W_{ka} \mu
^{(3)}_{ab} U_{bl} \widebar U_{bi}
U_{bj} \bigr) \biggr).
\end{eqnarray*}
Similar to (\ref{errorincov}), we find using the definition of $U$
and $W$ that the two last lines are asymptotic to $\frac{\mathcal
W_{ij,kl}(V)}{d_{\pi'}} + \frac{\mathcal W_{kl, ij}(V)}{d_\pi}$. Thus,
we get
%
%e9.14 #&#
\begin{equation}
\bigoplus_{\pi\in\Pi} \Psi^\pi_4
\stackrel{d} {\sim} \bigoplus_{\pi\in
\Pi}
\Psi^\pi.
\end{equation}
This concludes the proof.
\end{pf}

In order to conclude the proof of
Proposition~\ref{propgeneraldistributionJ}, we finally consider the
general case (generalization of Section~\ref{secgencase}). As in
Proposition~\ref{propgeneraldistribution}, in the general case we get a
deterministic shift $\bigoplus_{\pi\in\Pi} \mathcal S^\pi$, where
%
%e9.15 #&#
\begin{equation}
\label{929} \mathcal S^\pi:= \frac{(\llvert  d_\pi\rrvert +1)(\llvert  d_{\pi} \rrvert -1)^{1/2}}{d_\pi^4} \mathcal
S_{[\pi]}(V).
\end{equation}

The proof is similar to those of
Lemma~\ref{lemmaGreenfunctioncompwithshift} and
Proposition~\ref{propgeneraldistribution}. We take over the setup and
notation from Section~\ref{secgencase} up to, but not including,
(\ref{definitionofx}). For each $\pi\in\Pi$, we define the spectral
parameter $z_\pi:=\theta_\pi+ \ii N^{-4}$ and the $\llvert \pi\rrvert
\times\llvert  \pi\rrvert $ matrix
%
%e9.16 #&#
\begin{equation}
\label{definitionofxJ} x_R^\pi:= N^{1/2}
\bigl(\llvert d_\pi\rrvert - 1\bigr)^{1/2} \bigl(V^*
R(z_\pi ) V - m(z_\pi) \bigr)_{[\pi]},
\end{equation}
we well as the $\llvert  [\Pi] \rrvert  \times\llvert  [\Pi] \rrvert $
block diagonal matrix $x_R:=\bigoplus_{\pi\in\Pi} x_R^\pi$. The
quantities $x_S$ and $x_T$ are defined analogously with $R$ replaced by
$S$ and $T$, respectively. The following is the main comparison
estimate, which generalizes
Lemma~\ref{lemmaGreenfunctioncompwithshift}.

%
%le9.6 #&#
\begin{lemma} \label{lemmaGreenfunctioncompwithshiftJ}
Provided $\rho$ is a large enough constant, the following holds. Let $f
\in C^3  (\mathbb{C}^{\llvert  [\Pi] \rrvert  \times\llvert [\Pi]
\rrvert } )$ be bounded with bounded derivatives and $q \equiv q_N$ be
an arbitrary deterministic sequence of $\llvert  [\Pi] \rrvert
\times\llvert  [\Pi] \rrvert $ matrices. Then
%
%e9.17 #&#
%e9.18 #&#
\begin{eqnarray}
\mathbb{E}f(x_T + q) & =&
\mathbb{E}f(x_R + q) + \sum_{i,j \in
[\Pi]}Z^{(ab)}_{ij} \mathbb{E}\frac{\partial f}{\partial x_{ij}}(x_R +q)
\nonumber\\[-8pt]\label{comparisonforTJ} \\[-8pt]
&&{}  + A_{ab} + O \bigl(\varphi^{-1} \mathcal E_{ab} \bigr),\nonumber
\\
\label{comparisonforSJ} \mathbb{E}f(x_S + q) & =&
\mathbb{E}f(x_R + q) + A_{ab} + O \bigl(
\varphi^{-1} \mathcal E_{ab} \bigr),
\end{eqnarray}
where $A_{ab}$ satisfies $\llvert  A_{ab} \rrvert  \leq\varphi^{-1}$,
the error term $\mathcal E_{ab}$ is defined in (\ref{defofcalE}), and
$Z^{(ab)}$ is the $\llvert  [\Pi] \rrvert  \times\llvert  [\Pi] \rrvert
$ block diagonal matrix $Z^{(ab)}:=\bigoplus_{\pi\in\Pi} Z^{(ab), \pi}$
with $\llvert  \pi\rrvert  \times\llvert  \pi\rrvert $ blocks
\begin{eqnarray}
Z^{(ab),\pi}_{ij}:= - N^{-1} \bigl(\llvert
d_\pi\rrvert - 1\bigr)^{1/2} \bigl(m(z_\pi)^4
\mu^{(3)}_{ab} \widebar V_{ai} V_{bj} +
m(z_\pi)^4 \mu ^{(3)}_{ba} \widebar V_{bi} V_{aj} \bigr) \nonumber
\\
\eqntext{(i,j \in\pi).}
\end{eqnarray}
\end{lemma}

\begin{pf}
The proof of Lemma~\ref{lemmaGreenfunctioncompwithshift} may be taken
over almost verbatim, following the proof of Lemma~7.13 of \cite{KY2}.
\end{pf}

The comparison estimate from
Lemma~\ref{lemmaGreenfunctioncompwithshiftJ} yields the shift described
by $\mathcal S$. The precise statement is given by the following
proposition, which generalizes
Proposition~\ref{propgeneraldistribution}.

%
%pr9.7 #&#
\begin{proposition} \label{propjointGFCconclusion}
For large enough $K$, we have
\[
\bigoplus_{\pi\in\Pi} X^\pi\stackrel{d} {\sim}
\bigoplus_{\pi\in\Pi} \bigl(\widetilde \Upsilon{}^\pi+
\Psi^\pi+ \mathcal S^\pi \bigr),
\]
where $\mathcal S^\pi$ was defined in (\ref{929}).
\end{proposition}

\begin{pf}
As in the proof of Proposition~\ref{propgeneraldistribution}, we follow
the proof of Theorem~2.14 in Section~7.4 of \cite{KY2}. The inputs are
Proposition~\ref{propalmostgaussian2} and
Lemma~\ref{lemmaGreenfunctioncompwithshiftJ}.
\end{pf}

Now Proposition~\ref{propgeneraldistributionJ} follows immediately from
Proposition~\ref{propjointGFCconclusion} using $\Upsilon^\pi=
\widetilde \Upsilon{}^\pi+ \mathcal S^\pi$. This concludes the proof of
Proposition~\ref{propgeneraldistributionJ}.

%s10 #&#
\begin{appendix}\label{appendix}
\section*{Appendix: Near-degenerate perturbations}
\setcounter{theorem}{0}
\setcounter{equation}{0}

In this appendix, we record some basic results on the perturbation of
near-degenerate spectra.

%
%pr10.1 #&#
\begin{proposition} \label{propperturbation}
Let $A$ and $B$ be nonzero Hermitian matrices on $\mathbb{C}^N$. Let
$n + m =
N$, so that $\mathbb{C}^N = \mathbb{C}^{n} \oplus\mathbb{C}^{m}$,
and assume that $A$ and $B$
are of the form
\[
A = \pmatrix{ A_{11} & 0 \vspace*{2pt}
\cr
0 & A_{22} },
\qquad B = \pmatrix{ 0 & B_{12} \vspace*{2pt}
\cr
B_{21} & 0
}
\]
(in self-explanatory notation). Define the spectral gap
\[
\Delta:= \operatorname{dist} \bigl(\sigma(A_{11}), \sigma(A_{22})
\bigr)
\]
and assume that $\Delta\geq3 \llVert  B \rrVert $.

Define the domain
\[
\mathcal D:= \bigl\{\mu\in\mathbb{C}\dvtx \operatorname{dist} \bigl(\mu,\sigma
(A_{11}) \bigr) < 2 \llVert B \rrVert \bigr\}.
\]
Then $A + B$ has exactly $n$ eigenvalues $\mu_1 \leq\cdots\leq\mu_n$ in
$\mathcal D$ (counted with multiplicity), which satisfy
\[
\bigl\llvert \mu_i - \lambda_i(A_{11})
\bigr\rrvert \leq\frac{\llVert
B \rrVert ^2}{\Delta
- 2 \llVert  B \rrVert } \qquad(i = 1, \ldots, n).
\]
\end{proposition}

\begin{pf}
The eigenvalue--eigenvector equation reads $(A + B) \mathbf{x} =
\mu\mathbf{x}$. Writing $\mathbf{x} = (\mathbf{x}_1, \mathbf{x}_2)
\in\mathbb{C}^n \oplus\mathbb{C}^m$ leads to the system
%
%e10.1 #&#
\begin{eqnarray}\label{systemforperturbation}
A_{11} \mathbf{x}_1 + B_{12} \mathbf{x}_2 &=& \mu\mathbf{x}_1,
\nonumber\\[-8pt]\\[-8pt]
A_{22} \mathbf{x}_2 + B_{21}\mathbf{x}_1 &=& \mu\mathbf{x}_2.\nonumber
\end{eqnarray}
By assumption, for $\mu\in\mathcal D$ we have
%
%e10.2 #&#
\begin{equation}
\label{muawayfromQ}
\operatorname{dist}\bigl(\mu, \sigma(A_{22})\bigr) \geq
\Delta- 2 \llVert B \rrVert.
\end{equation}
Since $\Delta- 2 \llVert  B \rrVert  \geq\llVert  B \rrVert  > 0$, we
find that (\ref{systemforperturbation}) is equivalent to the system
\[
\mathbf{x}_2 = - (A_{22} - \mu)^{-1}
B_{21} \mathbf{x}_1, \qquad A_{11}
\mathbf{x}_1 - \mu\mathbf{x}_1 - B_{12}
(A_{22} - \mu)^{-1} B_{21} \mathbf{x}_1
= 0.
\]
Replacing $B$ with $tB$ for $t \in[0,1]$, we conclude that for $\mu
\in\mathcal D$ we have the equivalence
\[
\mu\in\sigma(A + t B) \quad\Longleftrightarrow\quad f_t(\mu) = 0,
\]
where
\[
f_t(\mu):= \det \bigl(A_{11} - \mu- t^2
B_{12} (A_{22} - \mu)^{-1} B_{21}
\bigr).
\]
Moreover, from Lemma~\ref{lemmagrowthofspectrum} below we find that
$\mathcal D$ contains exactly $n$ eigenvalues of $A + t B$, for all $t
\in[0,1]$. It is well known that the eigenvalues $\mu_i(t)$ of $A + t
B$ are continuous in $t$. We now claim that each such continuous
$\mu_i(t)$ is in fact Lipschitz continuous with Lipschitz constant
\[
L:= \frac{\llVert  B \rrVert ^2}{\Delta- 2 \llVert  B \rrVert
}.
\]
Assuming this is proved, the claim immediately follows from $\llvert
\mu_i - \lambda_i \rrvert  = \llvert  \mu_i(1) - \mu_i(0) \rrvert
\leq L$.

In order to prove the Lipschitz continuity of $\mu_i(t)$, note that
$\mu_i(t)$ is an eigenvalue of the matrix
\[
X_i(t):= A_{11} - t^2 B_{12}
\bigl(A_{22} - \mu_i(t)\bigr)^{-1}
B_{21}.
\]
Then the Lipschitz continuity of $\mu_i(t)$ follows readily from
Lemma~\ref{lemmagrowthofspectrum} below and the estimate
\[
\bigl\llVert B_{12} \bigl(A_{22} - \mu_i(t)
\bigr)^{-1} B_{21} \bigr\rrVert \leq L
\]
as follows from (\ref{muawayfromQ}), the fact that $\mu_i(t) \in
\mathcal D$ for all $t \in[0,1]$, and the fact that $A_{22}$ is
Hermitian.
\end{pf}

%
%le10.2 #&#
\begin{lemma} \label{lemmagrowthofspectrum}
Let $A$ and $B$ be square matrices, with $A$ Hermitian. Then the
spectrum of $A+B$ is contained in the closed $\llVert  B \rrVert
$-neighborhood
of the spectrum of $A$.
\end{lemma}

\begin{pf}
Using the identity $(A + B - z)^{-1} = (A - z)^{-1}  (1 + B (A -
z)^{-1} )^{-1}$ we conclude that if $\operatorname{dist}(z, \sigma(A))
> \llVert  B \rrVert $ then $z \notin\sigma(A + B)$.
\end{pf}
\end{appendix}


%suskaldyti doi

% imsref loaded by linak, 2013-10-23 11:02:53

\printaddresses

\end{document}